\RequirePackage{fix-cm}
\documentclass[twocolumn,showkeys]{revtex4}          
\usepackage{graphicx}
\usepackage{amssymb}
\usepackage{amsbsy}
\usepackage{subfigure}
\usepackage{color}
\usepackage{float}

\newcommand{\modif}[1]{#1}
\newcommand{\modifnew}[1]{#1}

\newcommand{\bm}[1]{{\bf #1}}
\newcommand{\R}{\mathbb{R}}

\newtheorem{algo}{Algorithm}

\begin{document}

\title[On stochastic generation of RVEs for analysis of composites within the framework 
of homogenization]{On efficient and reliable stochastic generation of RVEs for analysis of composites within the framework 
of homogenization
}

\author{Vladimir Salnikov} 
\author{Daniel Cho\"i}
\author{Philippe Karamian-Surville}

\affiliation{Nicolas Oresme Mathematics Laboratory \\ University of Caen Lower Normandy\\ 
  CS 14032, Bd. Mar\'echal Juin,  BP 5186\\
  14032, Caen Cedex,  France \\
  \email{vladimir.salnikov@unicaen.fr, daniel.choi@unicaen.fr, philippe.karamian@unicaen.fr}
}


\date{18 March 2014}

\begin{abstract}
In this paper we describe efficient methods of generation of representative volume 
elements (RVEs) suitable for producing the samples for analysis of effective properties 
of composite materials via and for stochastic homogenization. We are interested in 
composites reinforced by a mixture of spherical and cylindrical inclusions. For 
these geometries we give explicit conditions of intersection in a
convenient form for verification. Based on those conditions we present two methods to 
generate RVEs: one is based on the Random Sequential Adsorption scheme, the 
other one on the time driven Molecular Dynamics. We test the efficiency of these methods 
and show that the first one is extremely powerful for low volume fraction of inclusions, 
while the second one allows us to construct \modif{denser} configurations.   
All the algorithms are given explicitly so they can be implemented directly. 

\end{abstract}

\maketitle

\section{Introduction/motivation}
\label{sec:intro}

In this paper we describe some approaches to generate the 
representative volume elements (RVE) in order to estimate
the effective properties of composite materials within the 
framework of stochastic homogenization.

Our consideration of this problem is motivated by direct applications, namely 
the estimation of mechanical, thermal and electrical properties 
of composites reinforced by spherical inclusions, microtubes as well as
inclusions of irregular shapes. 
Our main approach to estimate  these properties is based on 
 homogenization techniques, so we need an efficient algorithm of generation 
of stochastic RVEs; in particular it should be sufficiently fast 
and require minimal interaction with the user. The nature of the considered 
materials defines the geometries that we need to generate, namely 
we focus our attention on a mixture of spherical and cylindrical inclusions
that are not allowed to intersect or to overlap.

There has been a number of works on the RVE generation 
for various geometries of inclusions: spheres (\cite{segu}, \cite{levesque_sp}), 
ellipsoids (\cite{man} -- \cite{levesque_ell}), 
``spherocylinders'' (cylinders with half-spheres attached to the ends -- 
\cite{williams}, \cite{zhao}). The common point of all these geometries is that 
the relative positions of two figures is characterized by simple 
algebraic conditions. For the applications that interest us we would need 
also the inclusions of true cylindrical form, since we are going to 
use them as elementary ``building blocks'' for designing and optimizing
more complex geometries. The purpose of this paper is thus to present 
appropriate methods of dealing with such geometries \modifnew{in different situations
occuring in pratice}, give an estimation of their efficiency and implementation details.

To generate an RVE one can adapt several approaches that are more or less efficient for 
different geometries. The most natural and probably historically 
the first one is random sequential adsorption (RSA, see for example \cite{rsa1,rsa}) 
-- random generation of the parameters of the geometry 
and verification if these parameters satisfy the imposed conditions, like the intersection one.
In practice one considers an empty RVE and starts generating the 
inclusions one after the other \modif{and} rejecting those that do not verify 
the conditions. This process usually solves the problem when the 
volume fraction of inclusions is sufficiently small, otherwise the 
generation process can take a long time or even \modif{get} stuck while the 
RVE is still far from the theoretical volume fraction of inclusions. 

Another family of approaches is inspired by molecular dynamics (MD): basically the 
generated inclusions are allowed to interact and \modif{change} until the desired 
configuration is constructed.
Here, one should 
distinguish two qualitatively different techniques: ``event driven'' and ``time driven''
simulation. For the ``event driven'' MD one is not interested in dynamics 
itself but only in a particular configuration when some \emph{event} occurs: 
that can be a collision between inclusions or an interaction 
of an inclusion with the boundary of the considered region. For every such 
event, the parameters (coordinates, velocities, angular velocities, sizes...)
of the inclusions are updated and the process is repeated until the desired configuration 
is achieved. \modif{A nice description of such a process with an example of rigid disks in a plane
can be found in \cite{LSA}; because of this publication the process is often 
called Lubachevsky--Stillinger algorithm.}
The difference of the ``time driven'' MD is that the 
parameters are updated at each time step. The former method is more efficient 
provided that there is an easy way to compute the time of the next event.
In practice, however, already for simple geometric shapes like ellipsoids, 
it is not an easy task (\cite{levesque_ell}): one needs to predict the time 
of collision of moving inclusions and provide a consistent model of interaction itself.
In this paper, we explain how the above mentioned methods 
can be adapted in order to generate the RVEs with cylindrical and spherical inclusions. 
We describe in details the random generation of non-intersecting inclusions,
as well as the relaxation procedure allowing us to produce non-intersecting configurations from the 
intersecting ones. 

The paper is organized as follows. In the next section we discuss the 
conditions of intersection of spheres with cylinders and of cylinders 
between themselves under convenient form for verification. 
This is already sufficient to implement the random generation strategy
described in section \ref{sec:MC}. 
As an alternative to this strategy in section  \ref{sec:dynamic_geometry} we present a relaxation procedure: we describe a model 
for dynamics of intersecting inclusions which is then used in the time driven MD 
simulation. We comment on the 
mechanics behind this model as well as on the computational and implementation details. 
The short section \ref{sec:compare} is devoted to comparison of these methods. 
We observe that for low volume fractions both strategies are acceptable, while 
for higher fraction (which we may need in the applications)
only the latter one produces satisfactory results in reasonable time.
We also study the dependence of time needed to generate the RVE samples
depending on the geometric parameters of the configuration: 
relative fraction of spheres and cylinders, size and number of them.

\section{Intersection conditions}
\label{sec:static_geometry}

In this section we describe the preliminaries related to 
the conditions of intersection of the geometric shapes 
that will be used as inclusions in the RVE generation. 
We present these conditions in the form of small algorithms 
to make their implementation transparent. 
In each algorithm we not only detect the intersection 
but also specify the type of it (notations of the form $sc1, sc2, cc1\ldots$).
The reason for this is that the application of these algorithms is twofold: 
in section \ref{sec:MC} we use them as a parts of RSA-type methods
(only the detection itself is important), but we also use them 
to define the MD interaction laws in section \ref{sec:dynamic_geometry}
(one does need to distinguish various types). 
These types of intersections \modif{will be described after each algorithm; an illustration of each of them can also} 
be seen in the appendix 1 (figures \ref{fig:sp-cyl1} --
\ref{fig:cyl-cyl9}).

\subsection{Sphere with cylinder}
Let us start with a simpler case of intersection of spheres 
with cylinders. Throughout this paper, we will characterize 
a sphere $S$ by its central point $\bm p_s \in \R^3$ and its radius
$r_s \in \R$, and denote it $S(\bm p_s, r_s)$ The parameters for a cylinder 
$C(\bm p_c, r_c, \bm l_c)$ will 
be its central point $\bm p_c \in \R^3$, its radius $r_c \in \R$
and the direction of the axis of symmetry $\bm l_c \in \R^3$ -- 
a vector which is not normalized and thus encodes also the 
information about the (half of the)length of a the cylinder. 
In what follows we will use the notation $a$
for the aspect ratio of a cylinder $a \equiv \|\bm l_c\|/r_c$.

\begin{figure}[ht]  
\centering
\subfigure[\,\modif{3D view}]{    
\includegraphics[trim = 0cm 0cm 0cm 17cm, clip, width=0.9\linewidth]{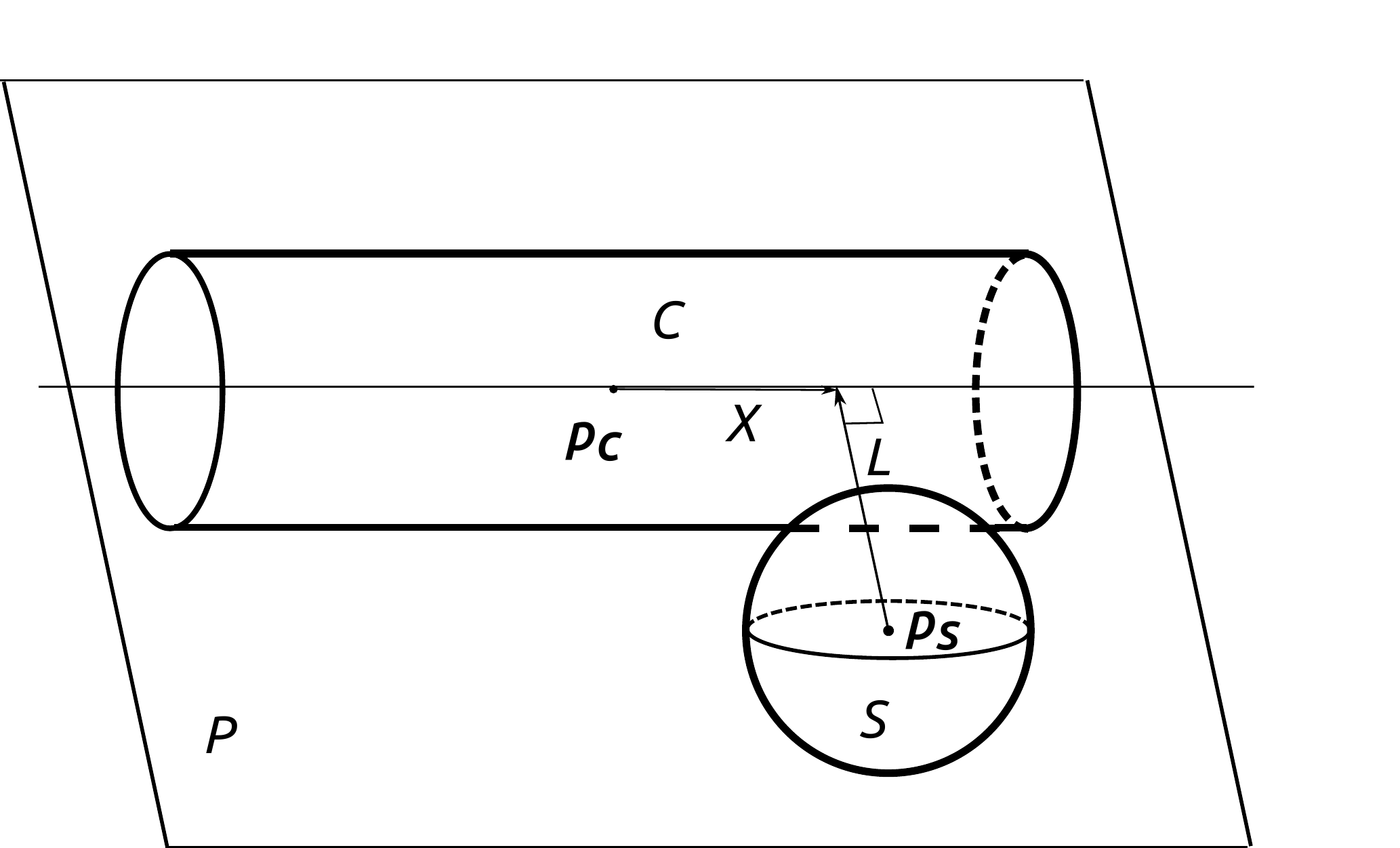}
  
}
\subfigure[\, \label{fig:sp-cyl-inter_b} The 
plane $P$ passing through the axis of symmetry of the cylinder $l_c$
and the center of the sphere $\bm p_s$. $\bm p_s$ in the waved/dotted regions 
corresponds to intersection of type $sc2$ and $sc3-sc4$ respectively]{
   \includegraphics[width=0.9\linewidth]{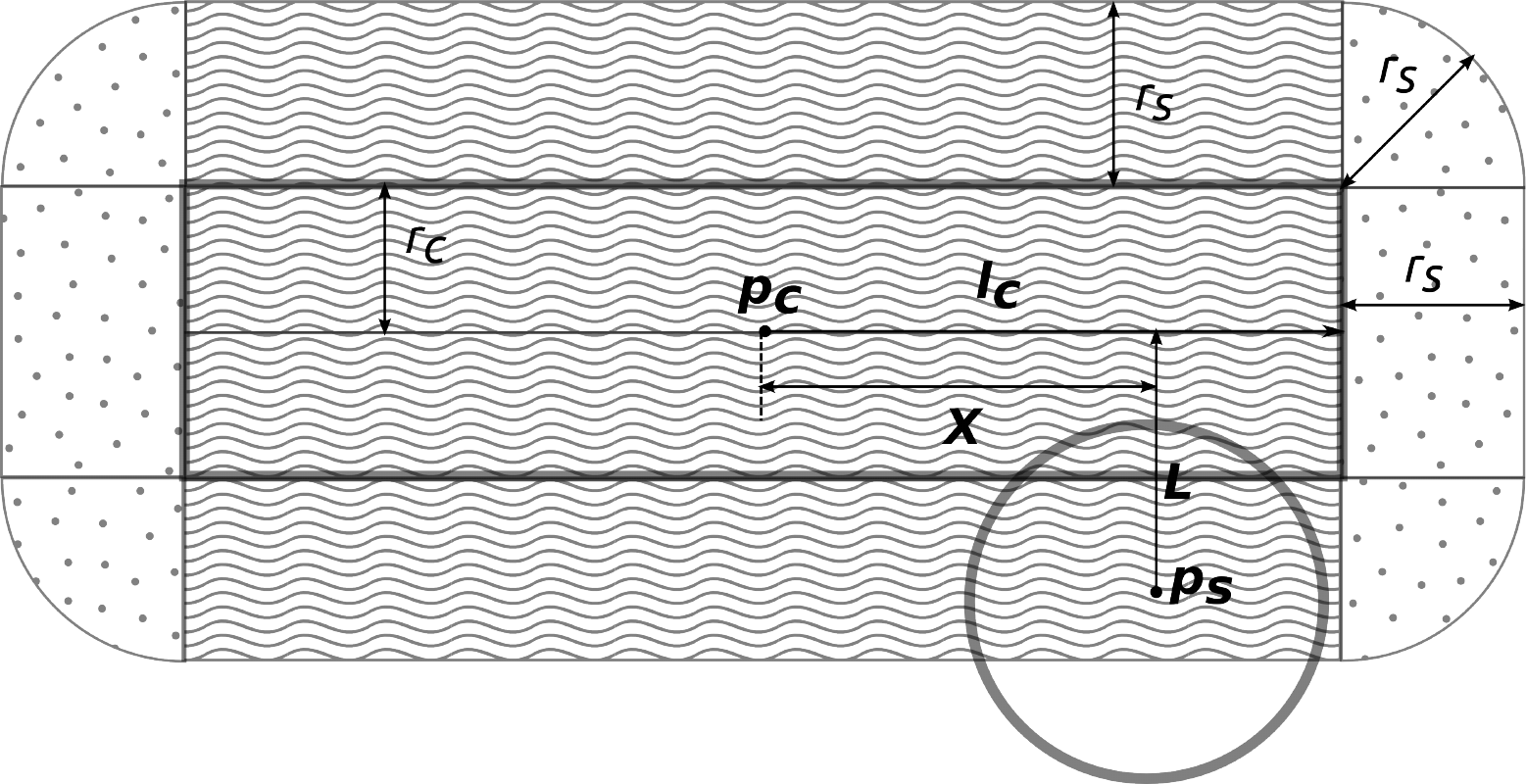}
 }

\caption{Intersection of a sphere and a cylinder}
\label{fig:sp-cyl-inter}
\end{figure}

To guarantee that a sphere $S$ does not intersect with a cylinder $C$,
one needs to verify two conditions: that there is no intersection 
neither with the cylindrical (curved) face nor with any of the two 
bases (extreme disks in the orthogonal sections). Two important 
quantities in the process are $L$ --- the distance from $\bm p_s$
to the symmetry axis of $C$ and $X$ --- the distance from $\bm p_c$ to the 
orthogonal projection of $\bm p_s$ \modif{to the axis}, depending on them we 
check if $\bm p_s$ is in the domain of intersection 
(\modif{waved or dotted regions of} fig \ref{fig:sp-cyl-inter_b}).

The following algorithm describes the computation.
\begin{samepage}
 \begin{algo} \label{alg:sp-cyl} Input: $S(\bm p_s, r_s), C(\bm p_c, r_c, \bm l_c)$. \newline
 \begin{tabular}{c|l}   & \parbox{0.95\linewidth}{  
\begin{enumerate}
  \item   Compute $X = (\bm p_s - \bm p_c)\cdot \bm l_c/\|\bm l_c\|$
    
  \item if $( |X|>\|\bm l_c\| + r_s )$  \newline
      \begin{tabular}{c|l}   & \parbox{0.8\linewidth}{

    \begin{enumerate}
       \item if $(\|\bm p_c - \bm p_s \|<r_s)$ $\to$ \emph{intersection of type sc1, stop.} 
	      \newline
	     else $\to$ \emph{no intersection, stop}
    \end{enumerate}
    
    }  \end{tabular}
    
  \item if $( |X|<\|\bm l_c\| )$ \newline
  \begin{tabular}{c|l}   & \parbox{0.8\linewidth}{

    \begin{enumerate}
       \item Compute $L = \sqrt{ (\bm p_s - \bm p_c)^2 - X^2 }$
       \item if $(L < r_s + r_c)$ $\to$ \emph{intersection of type sc2, stop.} 
	      \newline
	     else $\to$ \emph{no intersection, stop}
    \end{enumerate}
    }  \end{tabular}

  \item if $(\|\bm l_c\| \le |X| \le \|\bm l_c\| + r_s )$  
  \newline \begin{tabular}{c|l}   & \parbox{0.95\linewidth}{
     \begin{enumerate}
         \item       Compute $L = \sqrt{ (\bm p_s - \bm p_c)^2 - X^2 }$
         \item if $( L < \sqrt{r_s^2 - (|X| - \| \bm l_c \|)^2} + r_c)$ 
	    \newline \begin{tabular}{c|l}   & \parbox{0.8\linewidth}{
              \begin{enumerate}
                   \item if$(|L| < r_c)$ $\to$ \emph{intersection of type sc3, stop.} 
                   \item if$(|L| \ge r_c)$ $\to$ \emph{intersection of type sc4, stop.} 
		 \end{enumerate}
	      }  \end{tabular} \newline
		else $\to$ \emph{no intersection, stop.}  
	 
     \end{enumerate}
   } \end{tabular} \newline
     
\end{enumerate}
 }\end{tabular}    
\end{algo}
\end{samepage}
\modif{The intersection of type $sc2$ is the most natural, when the sphere touches the cylindrical 
face. Types $sc3$ and $sc4$ correspond to the sphere intersecting the disk base of the cylinder; that 
depends on whether or not the center of the sphere is inside the infinite cylinder. Type $sc1$ corresponds 
to a degenerate situation when the cylinder is inside the sphere.}

\subsection{Two cylinders}
Let us now turn to the intersection of two cylinders $C_1$ and $C_2$.
In contrast to the previous situation, we have potentially 
four geometries of intersection, that is all combinations 
of cylindrical faces or base disks. 

\begin{figure}[ht]  
\centering
\includegraphics[width=0.9\linewidth]{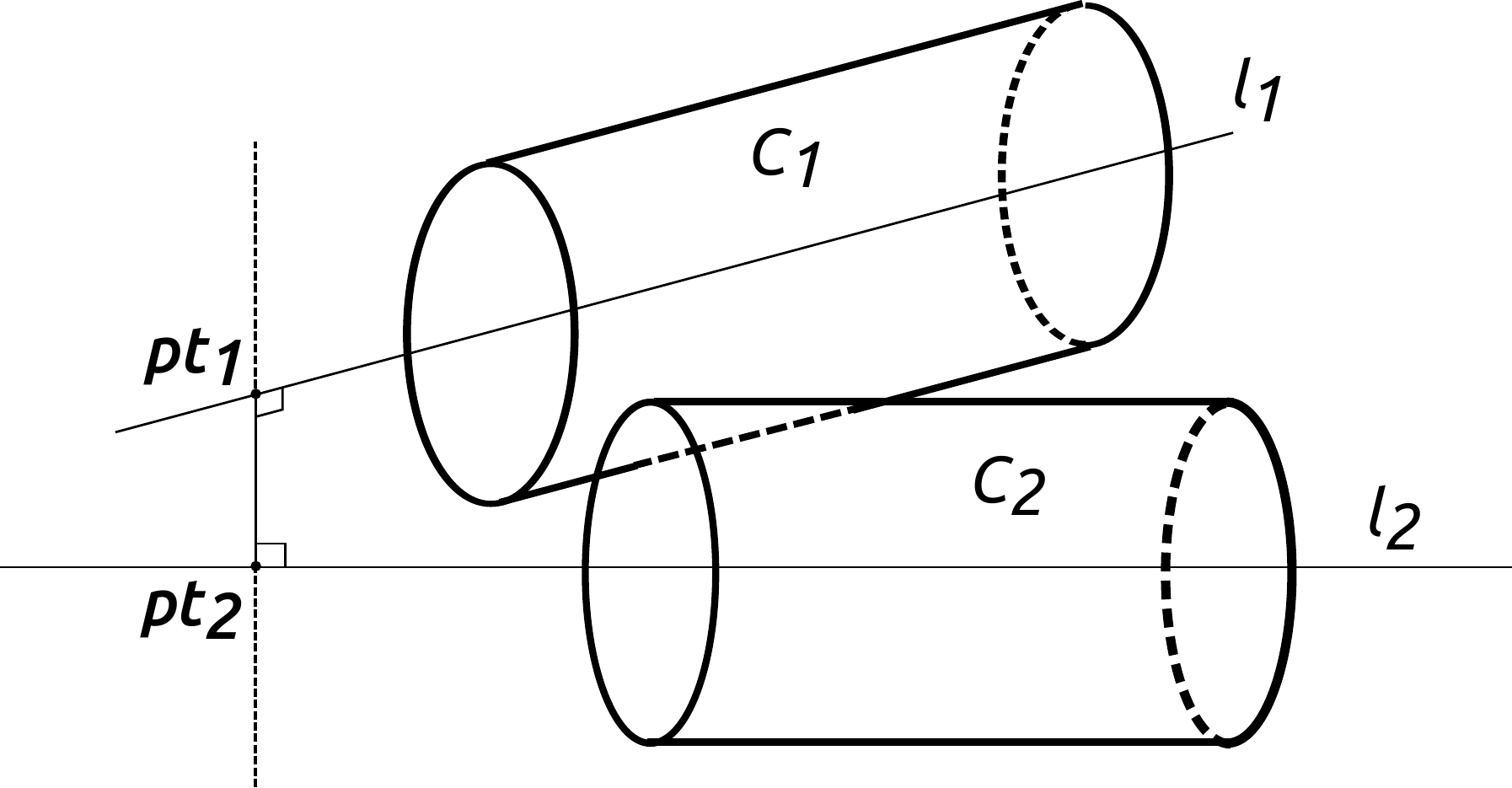}  
\caption{Geometry of intersecting cylinders}
\label{fig:skew-lines}
\end{figure}

An important (and not obvious from the first impression) observation 
is that one can 
distinguish these cases in terms of simple geometry 
of skew (i.e. not coplanar) straight lines. More precisely, consider the 
symmetry axes $l_1$ and $l_2$ of the cylinders, let $\bm{pt}_1$ and $\bm{pt}_2$
be the intersection points of the respective axes with the orthogonal
line realizing the distance between them (fig. \ref{fig:skew-lines}). 
If the cylinders are intersecting non-trivially (i.e. one is not inside the 
other one) but $\bm{pt}_1$ is not inside 
$C_1$ and $\bm{pt}_2$ is not inside $C_2$ 
then at least one of the disk bases of one cylinder intersects with the other cylinder.
Let us give an idea of the proof of this statement. If there is no intersection 
by the cylindrical face of one of the cylinders the result is automatic. 
If the symmetry axes are coplanar the statement is also trivial and follows from 
the intersection of rectangles in this plane.
Else, consider the intersection of infinite cylinders with the 
generating axes $l_1$ and $l_2$. They are quadratic surfaces in 
$\R^3$ and their intersection is given by continuous curves. 
Since the distance between the axes is less than the sum of radii 
of the cylinders there are points belonging to these curves 
lying in the plane orthogonal to $l_1$ that contains the segment $[\bm{pt}_1, \bm{pt}_2]$.
By assumption this plane is outside $C_1$. The similar picture is valid for the other cylinder. 
But since the cylinders do intersect, 
there is also a point on the surfaces of $C_1$ and $C_2$ belonging to one of these curves. 
Thus, by continuity, the curve has to intersect at least one circular boundary of one of the 
bases.

The above observation simplifies a lot the computation
process, since for detecting the intersection it is now sufficient 
to verify explicit algebraic conditions. \modif{A simple counting of floating point operations 
clearly shows that this} is much faster than 
solving a minimization problem to find the distance between two convex bodies.
Moreover it allows one to distinguish different types of intersection 
geometries. We give details of these conditions in the following algorithms.

\begin{figure}[ht]  
\centering
\includegraphics[trim = 0cm 0cm 9cm 19cm, clip,width=0.9\linewidth]{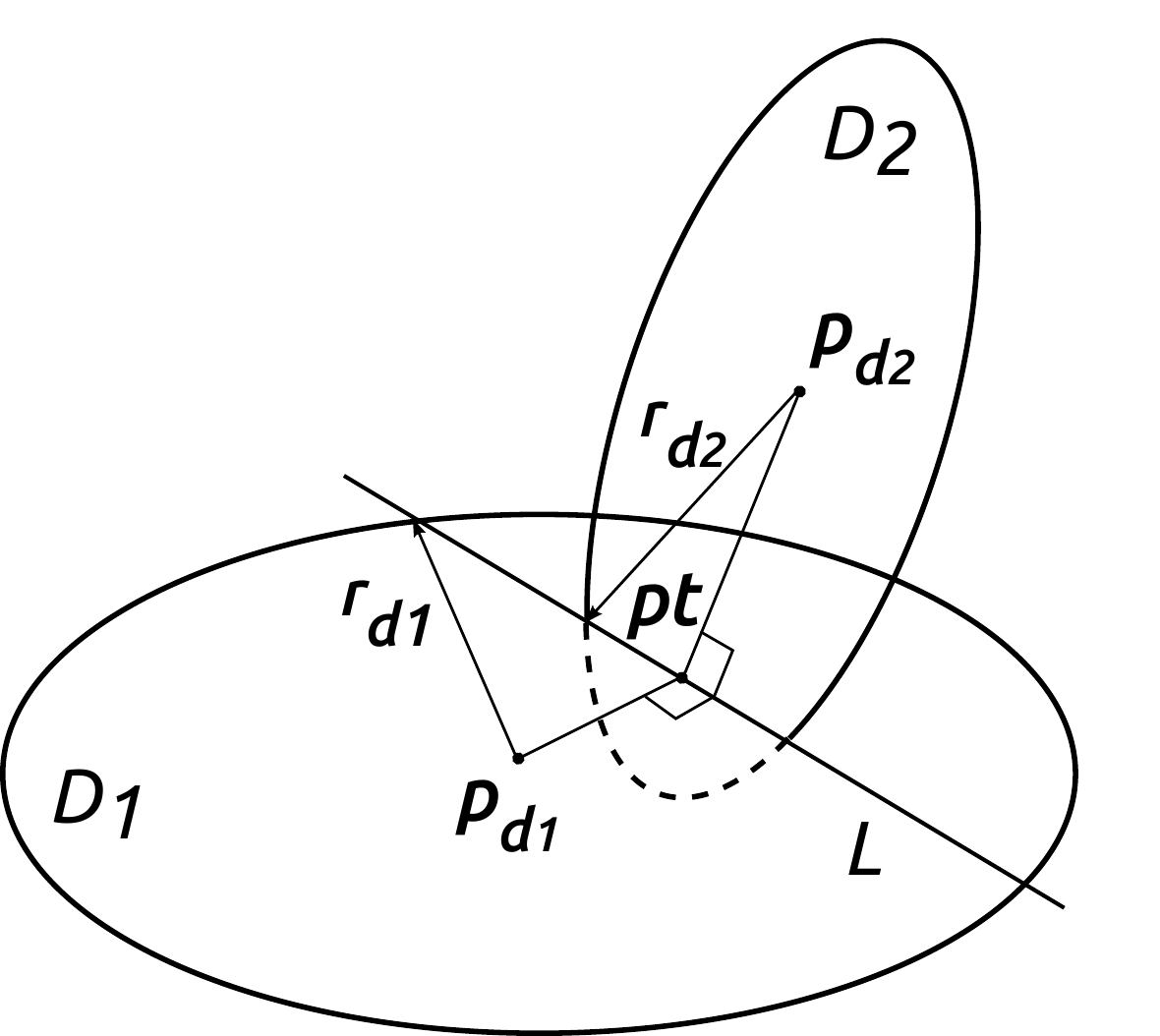}  
\caption{\modif{Geometry of intersecting disks. Type $d1$ depicted -- for type $d2$ exchange the disks.}}
\label{fig:disk-disk}
\end{figure}

First of all, let us consider the intersection between two disks in $\R^3$.
A disk $D$ is characterized by its center $\bm p_d$, radius $r_d$ and normal vector 
$\bm n_d$, and denoted by $D(\bm p_d, r_d, \bm n_d)$. 
To check the intersection we characterize the common line 
of the planes containing the disks and compare the distance from 
it to the centers of the disks with the respective radii.
\begin{samepage}
\begin{algo} \label{alg:d-d} Input: $D_1(\bm p_{d_1}, r_{d_1}, \bm n_{d_1}), D_2(\bm p_{d_2}, r_{d_2}, \bm n_{d_2})$
\newline \begin{tabular}{c|l}   & \parbox{0.95\linewidth}{  
\begin{enumerate}
  \item Compute the direction vector of the line $L$ of intersection of the 
  planes containing $D_1$ and $D_2$: $\bm n = \bm n_{d_1} \times \bm n_{d_2}$.
  \item Compute the (common) projection of the disk centers to $L$: 
  $\bm {pt} = \bm p_{d_1} + t \bm{v}$, where $\bm v = \bm n \times \bm n_{d_1}$, 
  $t = \bm n_{d_2} \cdot (\bm p_{d_1} - \bm p_{d_2})/(\bm n_{d_2} \cdot \bm v)$.
  \item if $(\|\bm{pt} -\bm p_{d_1}\| \le r_{d_1}) \& (\|\bm{pt} -\bm p_{d_2}\| \le r_{d_2} )  $ \newline
  \begin{tabular}{c|l}
       &     if $( r_{d_1}^2 - \|\bm{pt} -\bm p_{d_1}\|^2  >  r_{d_2}^2 - \|\bm{pt} -\bm p_{d_2}\|^2 )$ \\
        &   		  $\to$ \emph{intersection of type d1, stop} \\
          &      else $\to$ \emph{intersection of type d2, stop} \\
  \end{tabular}  \newline
      else $\to$ \emph{no intersection, stop}
\end{enumerate}
}\end{tabular}
\end{algo}
\end{samepage}
\modif{In the generic situation only one of the disks encounters the boundary circle of the other one 
-- this is the difference between the cases $d1$ and $d2$ (see also figure \ref{fig:disk-disk}). }

Now turn to the intersection of a disk with the cylindrical face of a cylinder. 
Here we will basically consider the point of intersection of the axis of the 
cylinder with the plane containing the disk and check if it is inside, 
close to, or far away outside the disk. This algorithm will treat only the cases 
that are not covered by the previous one, i.e. the intersection is considered 
only if the disk does not hit the circular boundary of the cylindrical face. 
\begin{samepage}
\begin{algo} \label{alg:cyl-d} Input: $C(\bm p_{c}, r_{c}, \bm l_{c}), D(\bm p_{d}, r_{d}, \bm n_{d})$
\newline \begin{tabular}{c|l}   & \parbox{0.95\linewidth}{  
\begin{enumerate}
  \item Compute the intersection point of the axis of the cylinder with the plane of the disk:
   $\bm {a} = \bm p_c + t \bm{l_c}$, where  
   $t = \bm n_{d}\cdot(\bm p_{d} - \bm p_{c})/(\bm n_{d} \cdot \bm l_c)$.
  \item Compute the point of the boundary of the disk which is the closest to $\bm a$:
   $\bm{pt_c} = \bm p_d + r_d\frac{\bm a - \bm p_d}{\|\bm a - \bm p_d\|}$.
  \item Compute the distance between the center of the cylinder and the projection $\bm b$ 
  of $\bm{pt_c}$ on it's axis: $X = (\bm{pt}_c - \bm p_c)\cdot\bm l_c/ \| \bm l_c\|$, 
  $\bm b = \bm p_c + X\cdot\bm l_c/ \| \bm l_c\|$.
  \item if $(|X| < \| \bm l_c \| )$
  \newline \begin{tabular}{c|l}   & \parbox{0.8\linewidth}{
            \begin{enumerate}
             \item if $(\|\bm p_d - \bm a \| > r_d) \& (\|\bm b - \bm{pt_c} \| < r_c)$
                \newline $\to$ \emph{intersection of type cd1, stop}, 
              \item if $(\|\bm p_d - \bm a \| < r_d) \& (\|\bm b - \bm{pt} \| < r_c)$
                \newline $\to$ \emph{intersection of type cd2, stop}, 
            \end{enumerate}
   } \end{tabular} \newline
   \item if $(\|\bm a - \bm p_d \| < r_d)$ $\to$ \emph{intersection of type cd3, stop} 
  
   \item if no intersection of type cd1, cd2 or cd3 \newline
   $\to$    \emph{no intersection, or disk-disk intersection, stop.}

   \end{enumerate}

} \end{tabular}
\end{algo}
\end{samepage}
\modif{The choice of the type of intersection depends on how far the disk penetrates into the 
cylindrical face: $cd1$ -- the disk just encounters the surface, $cd2$ -- it intersects the axis of symmetry, 
$cd3$ -- the whole cylinder goes through the disk.}

And finally, let us present the algorithm of verification of intersection 
between two cylinders in its whole generality. Here we use the statement from the beginning of the 
section allowing us to distinguish the cases of intersection by the base disks
and by the cylindrical surfaces. \modif{To the intersection types described above we add $cc1$
corresponding to both cylindrical surfaces intersecting.} 
\begin{samepage}
\begin{algo} \label{alg:cyl-cyl} Input: $C_1(\bm p_{c_1}, r_{c_1}, \bm l_{c_1}), C_2(\bm p_{c_2}, r_{c_2}, \bm l_{c_2})$  
\newline \begin{tabular}{c|l}   & \parbox{0.95\linewidth}{  
  \begin{enumerate}
   \item  Compute the vector parallel to the common normal to the symmetry axes of the 
   cylinders: \newline
   $\bm n = \frac{\bm l_{c_1} \times \bm l_{c_2}}{\|\bm l_{c_1} \times \bm l_{c_2}\|}$
   \item Compute the distance between the symmetry axes of the cylinders: 
   $\rho = | (\bm p_{c_1} - \bm p_{c_2})\cdot \bm n |$
   \item if $\rho > r_{c_1} + r_{c_2}$ $\to$ \emph{no intersection, stop.} \newline
	 else 
	    \newline \begin{tabular}{c|l}   & \parbox{0.9\linewidth}{
	      \begin{enumerate}
		\item Compute the normals to the planes containing $\bm n$ and the axes of the cylinders
			  respectively: 
			$\bm n_1 = \bm n \times \bm l_{c_1}$, $\bm n_2 = \bm n \times \bm l_{c_2}$.
		\item Compute the points realizing the distance between the axes: 
		    $\bm{pt}_1 = \bm p_{c_1} + t_1 \bm l_{c_1}$, $\bm{pt}_2 = \bm p_{c_2} + t_2 \bm l_{c_2}$, 
		    where \newline
		    $t_1 = (\bm p_{c_2} - \bm p_{c_1})\cdot\bm n_2/(\bm l_{c_1}\cdot \bm n_2)$, \newline
		    $t_2 = (\bm p_{c_1} - \bm p_{c_2})\cdot\bm n_1/(\bm l_{c_2}\cdot \bm n_1)$.
		\item if $(|t_1| \le 1) \& (|t_2| \le 1)$ $\to$ \emph{intersection of type cc1, stop}
		      \newline
		      else 
		      	    \newline \begin{tabular}{c|l}   & \parbox{0.9\linewidth}{
		      \begin{enumerate}
		       \item Using the algorithm \ref{alg:cyl-d},
			      check intersection of $C_{1}$ with the disks 
			      $D((\bm p_{c_2} + \bm l_{c_2}), r_{c_2}, \bm l_{c_2})$ and \newline
			      $D((\bm p_{c_2} - \bm l_{c_2}), r_{c_2}, - \bm l_{c_2})$
      		       \item Using the algorithm \ref{alg:cyl-d},
			      check intersection of $C_{2}$ with the disks 
			      $D((\bm p_{c_1} + \bm l_{c_1}), r_{c_1}, \bm l_{c_1})$ and \newline
			      $D((\bm p_{c_1} - \bm l_{c_1}), r_{c_1}, - \bm l_{c_1})$
		       \item Using the algorithm \ref{alg:d-d},
			      check intersection of \newline
			      $D((\bm p_{c_1} + \bm l_{c_1}), r_{c_1}, \bm l_{c_1})$ with
			      $D((\bm p_{c_2} + \bm l_{c_2}), r_{c_2}, \bm l_{c_2})$, \newline
			      $D((\bm p_{c_1} + \bm l_{c_1}), r_{c_1}, \bm l_{c_1})$ with
			      $D((\bm p_{c_2} - \bm l_{c_2}), r_{c_2}, - \bm l_{c_2})$, \newline
			      $D((\bm p_{c_1} - \bm l_{c_1}), r_{c_1}, -\bm l_{c_1})$ with
			      $D((\bm p_{c_2} + \bm l_{c_2}), r_{c_2}, + \bm l_{c_2})$, \newline
			      $D((\bm p_{c_1} - \bm l_{c_1}), r_{c_1}, -\bm l_{c_1})$ with
			      $D((\bm p_{c_2} - \bm l_{c_2}), r_{c_2}, - \bm l_{c_2})$.
			\item if there is no intersection in the three points above 
			\newline$\to$ \emph{no intersection, stop.}
		      \end{enumerate}
		      	   } \end{tabular} \newline
	      \end{enumerate}
	   } \end{tabular} \newline
  \end{enumerate}
}\end{tabular}
\end{algo}
\end{samepage}

Let us make several remarks about the above algorithms. 
First, in their implementation, one should be careful about degenerate cases, 
like the symmetry axes of the cylinders close to intersecting or being parallel, 
or centers of the figures close to coinciding. Such situations might lead to 
some norms of vectors being close to vanishing. We didn't include this 
detail to the algorithms since they correspond to very 
explicit geometric configurations and their description makes the exposition 
too technical without giving any significant input. 
Second, having some extra information one can optimize a little the algorithms 
excluding some of the particular cases. For example in what 
follows we will discuss the generation of identical spheres and cylinders, that means that the point  
5 of the algorithm \ref{alg:cyl-d} never takes place. 
Let us also note that the algorithms presented in this section were constructed in 
such a way that one treats first the most frequent configuration, that is, in most of the 
cases the algorithm stops after very few operations. This remark 
is however valid with one exception of the point 
3.c. of the algorithm \ref{alg:cyl-cyl} where we consider all possible intersections
involving disks. If one is interested in intersection condition only (like in section 
\ref{sec:MC}) this can be optimized in an obvious way, but in section \ref{sec:dynamic_geometry} we will 
need all this information.

   \clearpage
  \newpage

\section{Random generation}
\label{sec:MC}

The algorithms presented in the previous section allows one to formulate 
explicitly the method of generating the RVE containing non-intersecting 
spherical and cylindrical inclusions. 
As we have agreed, a sphere is characterized by
its central point $\bm p_s$ and its radius
by $r_s$, the cylinder by its central point $\bm p_c$, its radius $r_c$
and the direction of the axis of symmetry $\bm l_c$. 
Suppose that the volume fractions  $f_s$, $f_c$ as well as the number 
$n_s$, $n_c$ of spheres and cylinders are given, let us also fix the aspect ratio
$a = \|\bm l_c\|/r_c$ of all the cylinders. This defines the size of all geometric shapes.
Thus, creating a random sphere is just generating three real numbers to form $\bm p_s$; 
\modif{for a cylinder one needs to provide} six numbers: for $\bm p_c$ and $\bm l_c$, and 
\modif{then} rescale $\bm l_c$
to fit the aspect ratio.
The natural RSA-type algorithm is then the following:
\begin{samepage}
\begin{algo} \label{alg:mc_gen} Input: $f_s, f_c, n_s, n_c, a$
\newline \begin{tabular}{c|l}   & \parbox{0.95\linewidth}{
\begin{enumerate}
 \item  Compute the radius of cylinders $r_c = \sqrt[3]{\frac{f_c}{2 \pi a n_c}}$
  \item  Compute the radius of spheres $r_s = \sqrt[3]{\frac{3f_s}{4 \pi n_s}}$
 \item  NumOfGenSpheres = 0, \newline NumOfGenCylinders = 0
 \item  while $(NumOfGenCylinders < n_c)$
      \newline \begin{tabular}{c|l}   & \parbox{0.8\linewidth}{
      \begin{enumerate}
       \item Generate a new cylinder \label{gencyl} 
       \item Using the algorithm \ref{alg:cyl-cyl} check (taking periodicity into account) if it intersects with any cylinder generated before
       \item if yes go back to \ref{gencyl}
       \item if no increase NumOfGenCylinders
       \end{enumerate}
        }\end{tabular}
 \item  while $(NumOfGenSpheres < n_s)$
      \newline \begin{tabular}{c|l}   & \parbox{0.95\linewidth}{
           \begin{enumerate}
       \item Generate a new sphere \label{gensp} 
       \item By comparing the distance between the centers with the sum of the radii, 
       check (taking periodicity into account) if it intersects with any sphere generated before
       \item if yes go back to \ref{gensp}
       \item Using the algorithm \ref{alg:sp-cyl} check (taking periodicity into account) if it intersects with any cylinder generated before
       \item if yes go back to \ref{gensp}
       \item if no increase NumOfGenSpheres
       \end{enumerate}
        }\end{tabular}
  
\end{enumerate}
}\end{tabular}
 
\end{algo}
\end{samepage}
As one sees from the above algorithm there is an issue of periodicity to deal with. 
The reason for this is that the concept of RVE has to cover the cases when the 
inclusions penetrate the boundary of the considered volume. The convention that is 
often used is that the part of the inclusion that exits the volume is mapped periodically 
to the other side of it. From the point of view of implementation it means that 
each generated object should have an attribute of intersecting the boundary 
of the RVE, which is assigned depending on its geometric properties.
(We leave the formulation of explicit conditions of this intersection 
as a simple exercise for a curious reader). 
If this attribute is present, the algorithm \ref{alg:sp-cyl} or \ref{alg:cyl-cyl}
should be applied to all the couples of the objects or their periodic images. 
One should not forget, that  one object can potentially intersect several 
RVE's boundaries, e.g. a sphere centered in the corner of the RVE should be 
considered as eight objects. 

\begin{figure}[ht] 
\centering
\includegraphics[trim = 1.7cm 1.5cm 6.5cm 19cm, clip, width=0.98\linewidth]{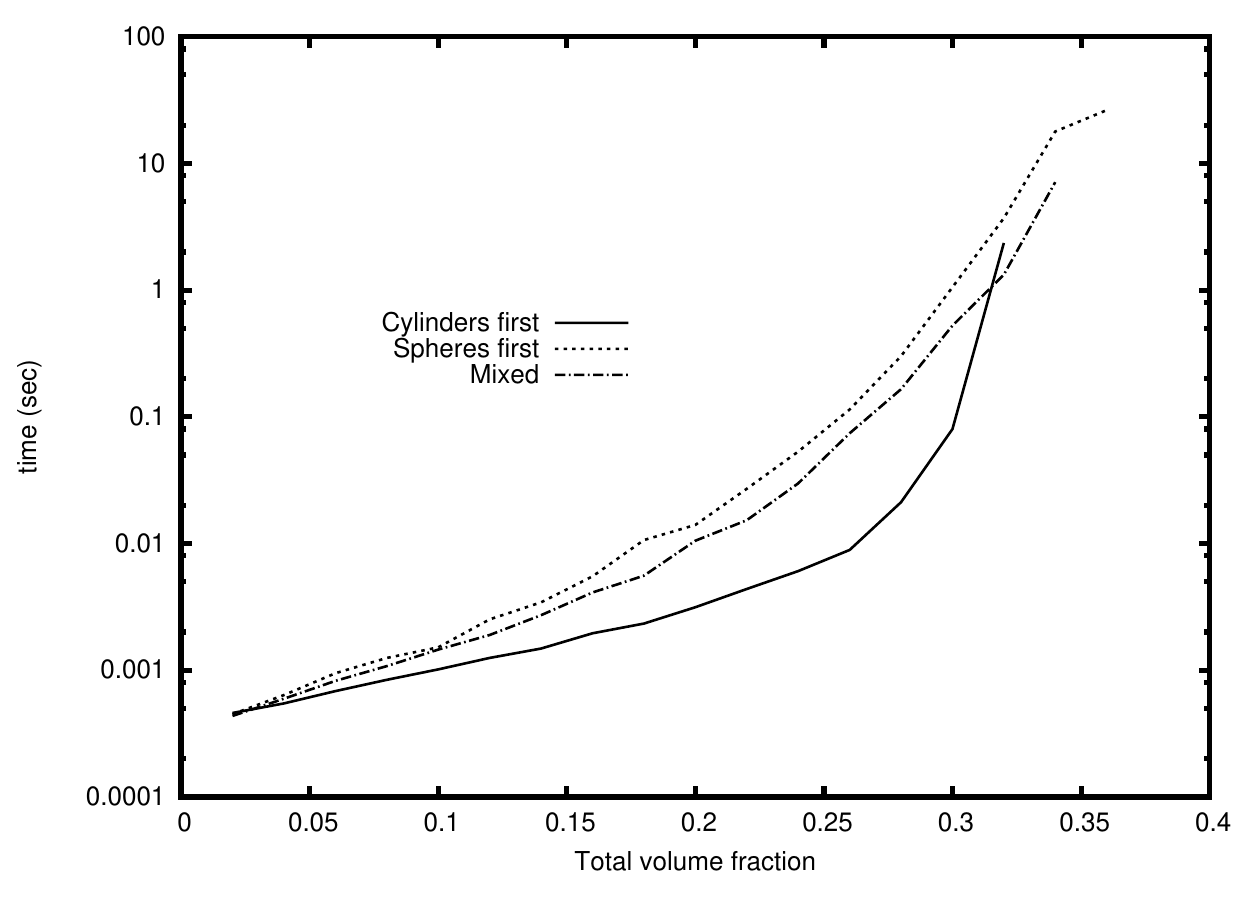}  
\caption{RSA: dependence of efficiency on the order choice strategy. 
Volume fraction distributed equally between 30 spheres and 30 cylinders of aspect ratio $a=5$, 
time estimation in seconds) averaged over 20 runs, \modif{volume fraction varies with the step $0.01$}.}
\label{fig:gen_order}
\end{figure}

Let us also note that one can suggest several versions of the above algorithm, 
all of the RSA type. Namely there is a freedom in the choice of the order of generation 
of spheres and cylinders: first all the cylinders, first all the spheres, or some mixed order.
We have made several tests of efficiency depending on the 
generation strategy. The results are shown on the figure \ref{fig:gen_order}\footnote{All the algorithms 
presented in this paper have been implemented using C++; the efficiency tests have been carried out on 
an Intel\textsuperscript{\textregistered} Core{\texttrademark} i7 960 3.2GHz machine running Kubuntu 13.10 
with the GNU compiler version 4.8.1.}: 
there is no particular optimal strategy, since generating first the cylinders gives
better performance on small volume fractions, but generating the spheres first 
permits us to achieve higher volume fractions, so the choice should apparently be 
made empirically depending on concrete applications.
In what follows we have chosen ``cylinders first''
strategy for comparison, since for higher volume fraction we will suggest another algorithm. 

We have performed several tests in order to study the dependence of the generation time 
on various parameters of the RVEs. The typical generation time is of order 
$10^{-3}$ seconds for small total volume fraction of inclusions (like $10\%$) 
and it may reach several seconds for relatively high one ($35\%$). 

We clearly observe that 
in many cases the same total volume fractions are easier achieved by generating spheres 
than cylinders, which is reasonable, since the cylinder intersection verification 
is longer. There is however a little saturation effect in the extreme case: 
when there are a lot of cylinders, representing a relatively low volume in total 
-- this apparently is influenced by our choice of the ``cylinders first''
strategy. 
It is also quite predictable that the same volume fraction is achieved faster with the 
smaller number of inclusions, outside the saturation effect certainly.
We also see that assembling cylinders with a higher aspect ratio (ratio between its length 
and diameter) is more difficult 
than with a lower one, which is also intuitively understandable. 
The details of these tests can be found in tables \ref{tab:RSA10} -- \ref{tab:RSA50} in the 
appendix 2.

During the tests we have fixed the maximal permitted generation time to $50$ seconds.
With this limitation we were able to reach the total volume fraction of $35\%$ with almost any 
distribution of this volume between spheres and cylinders. 
Certainly if one is sufficiently patient it is possible to generate the RVEs with higher 
volume fraction (permitted by geometry), with however a necessity to eventually relaunch the program
in the case of stagnating configuration.
As a global conclusion, we see that the method is very efficient for low volume fractions, and the generation time increases rapidly for medium values of order $30$ -- $35\%$.

\section{Time driven molecular dynamics}
\label{sec:dynamic_geometry}
We have seen in the previous section that as expected, the RSA-type 
algorithms of generation the RVEs are efficient when the desired volume fraction is 
small enough. We have also observed that for higher volume fractions 
the time of generation grows substantially and one cannot guarantee that the
desired number of inclusions will be reached even if the volume fraction 
is quite far from the theoretical value for a dense packing. This effect is 
perfectly explainable since it is rather easy to generate examples 
of ``bad'' geometries even with low number of inclusions. 

In the introduction we mentioned other strategies of generation of RVEs
based on improvements of the RSA-type algorithms or  
molecular dynamics. We have chosen to implement the so-called 
time driven version of the latter one. The scheme is the following: 
all the inclusions are generated not taking intersections in the account, 
then the interaction force is assigned to each couple of intersecting inclusions, 
the dynamics governed by these forces is described by a system of ODEs, that are solved 
numerically in order to achieve a relaxed configuration. In their dynamics the 
inclusions are also affected by \modif{damping} forces that are supposed to slow the system down and 
make it stay in the relaxed configuration.

Before going into details let us explain the choice of this strategy in comparison to 
a couple of other possible ones. We could have applied the pattern presented in 
\cite{levesque_ell} in the context of ellipsoids. The major problem 
there would be to formulate an efficient algorithm of collision time computation. 
We have seen in the section \ref{sec:static_geometry} that the intersection  
condition for cylinders can be verified algebraically, however,
in contrast to ellipsoids, there is no ``nice'' way to characterize it in terms of 
zeros of some function (to our knowledge at least). This means that 
using the rigid colliding cylinders for event-driven MD  would lead to solving a 
complicated minimization problem at each step. Moreover the collision problem 
for two rigid cylinders is also more involved. That is we doubt this 
approach to be more efficient than just integrating the ODEs. The other possible option 
would be to modify the RSA algorithm using some random moves of the inclusions, 
like in \cite{rsa+random}. There, the difficulty would be that because of the 
complexity of geometry under consideration one random move will not be enough to 
exclude all the intersections. A possible modification of this method is to consider 
several random moves to reorder the configuration -- this boils down to 
a rather classical approach to minimization of the functional using the Monte Carlo 
techniques (\cite{Monte-Carlo}) and needs the functional characterizing the intersection to be defined. 
We expect this approach to be very close to the direct time-driven MD both in the 
local relaxation trajectory and efficiency. 

Let us now turn to the description of the interaction of inclusions. 
As we have outlined before the idea is to introduce the forces 
if two inclusions intersect. In our model the forces will be 
of linear elastic nature, that is the value of the force is proportional 
to the depth of overlapping domain. The direction will be certainly chosen
to make the forces repulsive. 
For example for two spheres $S_1(\bm p_{s_1}, r_{s_1})$, $S_2(\bm p_{s_2}, r_{s_2})$: 
if $\| \bm p_{s_1} - \bm p_{s_2}\| < r_{s_1} + r_{s_2}$, the force acting 
on $S_1$ is $\bm F_{s_1} = - (\| \bm p_{s_1} - \bm p_{s_2}\| - r_{s_1} - r_{s_2}) 
\frac{ \bm p_{s_1} - \bm p_{s_2} }{\| \bm p_{s_1} - \bm p_{s_2}\|}$. 
The force acting on the second sphere is certainly opposite: $\bm F_{s_2} = - \bm F_{s_1}$.

For the forces acting on cylinders it is  necessary to specify the point of application,
since their motion also includes rotation. The table \ref{tab:forces} recapitulates all the 
forces that can act on a cylinder depending on the type of intersection 
from the algorithms \ref{alg:sp-cyl} -- \ref{alg:cyl-cyl}, using the same notations as in section \ref{sec:static_geometry}. The described forces always 
act on the first cylinder or the only one in the corresponding algorithm, the 
force acting on the other object is  opposite and has the same application point. 
\begin{table*}[ht] \centering
\begin{tabular}{|c|c|c|}
\hline
   & force $\bm F_c$ & application point $\bm{pt}_c$ \\
\hline \hline
 sc1 & $2\bm l_c \frac{\bm p_c - \bm p_s}{\|\bm p_c - \bm p_s\|} $ & $\bm p_c$ \\
\hline sc2 & $ - \sqrt{ (r_s + r_c)^2 - L^2} \frac{ \bm p_c - \bm p_s + X\frac{\bm l_c}{\|\bm l_c\|} }{ \| \bm p_c - \bm p_s + X\frac{\bm l_c}{\|\bm l_c\|}  \|}$ &
      $ \bm p_s + (r_s - \sqrt{(r_s + r_c)^2 - L^2})\frac{\bm F_c}{\|\bm F_c \|} $ \\
 \hline
 sc3 & $  -(\| \bm l_c \| + r_s - |X|)\frac{X}{|X|}\frac{\bm l_c}{ \| \bm l_c \| }$ &
   $ \bm p_c + X\frac{\bm l_c}{ \| \bm l_c \| }  - L\frac{\bm p_c - \bm p_s + X\frac{\bm l_c}{ \| \bm l_c \| }}{\| \bm p_c - \bm p_s + X\frac{\bm l_c}{ \| \bm l_c \| } \|}$ \\
 \hline
 sc4 & $  (\sqrt{r_s^2 - (|X| - \| \bm l_c \|)^2} + r_c - L) \frac{\bm{pt}_c - \bm p_s}{ \| \bm{pt}_c - \bm p_s \| }$ & 
 $  \bm p_c + X\frac{\bm l_c}{ \| \bm l_c \| }  - r_c\frac{\bm p_c - \bm p_s + X\frac{\bm l_c}{ \| \bm l_c \| }}{\| \bm p_c - \bm p_s + X\frac{\bm l_c}{ \| \bm l_c \| } \|}$ 
 \\
\hline \hline
 cc1 &$ (r_{c_1} + r_{c_2} - \rho) \frac{\bm{pt}_1 - \bm{pt}_2}{\| \bm{pt}_1 - \bm{pt}_2 \|}$ &		    
 $   (\bm{pt}_1 + \bm{pt}_2)/2$ \\
\hline 
 cd1 & $ -(\|\bm b - \bm{pt_c} \| - r_c) \frac{\bm b - \bm{pt_c} }{\|\bm b - \bm{pt_c} \|}$ &
 $ \bm p_d + r_d\frac{\bm a - \bm p_d}{\|\bm a - \bm p_d\|}$ \\
 \hline
 cd2 & $ -(\|\bm b - \bm{pt_c} \| - 2r_c) \frac{\bm{pt_c} - \bm b }{\|\bm b - \bm{pt_c} \|}$ &
 $ \bm p_d + r_d\frac{\bm a - \bm p_d}{\|\bm a - \bm p_d\|}$\\
 \hline
 cd3 & $ 2 r_c \frac{\bm p_c - \bm p_d}{\| \bm p_c - \bm p_d\|} $ & 
 $ \bm a$  \\
 \hline
 d1 & $ -(r_{d_2} - \|\bm{pt}_c -\bm p_{d_2}\|) \frac{\bm n_{d_1}}{\| \bm n_{d_1} \|} $ & 
 $ \bm p_1 + t \bm{v}$ \\
 \hline
 d2 & $ -(r_{d_1} - \|\bm{pt}_c -\bm p_{d_1}\|) \frac{\bm n_{d_2}}{\| \bm n_{d_2} \|} $& 
 $ \bm p_1 + t \bm{v}$\\
\hline 
\end{tabular} 
\caption{\label{tab:forces} Values of forces depending on the intersection type}
\end{table*}
To compute the force  between each couple of objects one needs to apply the 
algorithm \ref{alg:sp-cyl} or \ref{alg:cyl-cyl} and choose the corresponding expression from the table. 
Continuing the remark from the previous section let us note the intersections of types
cd1 -- cd3 and d1 -- d2 can be combined, it means that one may need to consider simultaneously 
the forces coming from these types of intersection.

These forces will enter directly in the differential equations governing the 
motion of the spheres described by the position of the center of mass $\bm p_s$ and 
its velocity $\bm v_s$.
As for the motion of the cylinders we represent it as a composition of a translation
characterized at every moment by the velocity $\bm v_c$ of the center $\bm p_c$ and 
a rotation around it characterized by the angular velocity $\boldsymbol{\omega}_c$
(according to \cite{arnold}, this can always be done).
If one knows these quantities at every moment 
this induces a relatively simple kinematic law of motion 
for the variables defining the position of the cylinder: 
\begin{equation}
\left\{\begin{array}{l}
  \dot \bm p_c = \bm v_c, \\[0.5em]
  \dot \bm l_c = \boldsymbol{\omega}_c \times \bm l_c.
   \end{array} \right.
\end{equation}
For constant $\bm v_c$ and $\boldsymbol{\omega}_c$ these 
equations  define the dynamics of a rigid body (cylinder) in the absence 
of external forces. The interaction with other inclusions results in 
additional evolution equations (that can be deduced using the techniques 
of Lagrangian mechanics \cite{arnold}):
\begin{equation}
\left\{\begin{array}{l}
  m \dot \bm v_c = \sum_{j} \bm F_j, \\[0.5em]
   \dot \bm{M} = \sum_{j} (\bm p_j - \bm p_c) \times \bm F_j, 
  \end{array} \right.
\end{equation}
where $m$ is the mass of the cylinder, $\bm M$ -- its angular momentum, 
and $(\bm p_j - \bm p_c)\times \bm F_j$ is the moment of the external force $\bm F_j$
applied at a point $\bm p_j$. At each moment there exists a linear 
operator relating $\bm M$ and $\boldsymbol{\omega}_c$, so in principle by inverting 
it, one can recover $\boldsymbol{\omega}_c$. However, this operation is not that much explicit.
The total angular momentum of the rigid 
body is $\bm M = \int_C \rho(\bm q) (\bm q \times (\boldsymbol{\omega}_c  \times \bm q)) d \bm q$, 
$\bm q$ being the spatial variable running along the whole body;
thus, the linear operator in general does depend on time, except for 
some specific choice of moving coordinate system that will have to 
be chosen independently for each body in the system at each timestep. 
To avoid this we make a couple of simplifications 
in the model. First, we assume that the mass in a cylinder 
is concentrated along its axis, i.e. dynamically the cylinder becomes 
a thin rod affected however by external forces 
applied to the whole volume. This simplification is legitimate since 
anyway we disregard the rotation of a cylinder around its symmetry axis. 
Already, this allows us to simplify the angular momentum to 
$\frac{1}{3}m \|\bm l\|^2 \sin^2(\alpha) \boldsymbol{\omega}_c$, where $\alpha$ is the angle between 
the angular velocity vector and the axis of the cylinder, and $\|\bm l\|$ 
is the (semi)length of the cylinder that remains constant.  
Note that this factor of $\sin^2(\alpha)$ defines the relation between the 
translation and rotation acceleration. As we are not interested in 
the precise dynamics of the system provided that it is qualitatively 
acceptable, we can simplify the equations even further 
by replacing this factor by its spatial mean value of $\frac{1}{2}$.
So the mechanical equations that we are finally solving read: 
\begin{equation} \label{mech}
\left\{\begin{array}{l}
\dot \bm p_c = \bm v_c, \quad \dot \bm l_c = \boldsymbol{\omega}_c \times \bm l_c, \quad  
m \dot \bm v_c = \sum_{j} \bm F_j,  \\[0.5em]
\frac{1}{6}m\|\bm l\|^2 \dot {\boldsymbol{\omega}}_c = \sum_{j} (\bm p_j -\bm p_c) \times \bm F_j. 
\end{array} \right.
\end{equation}
We will see that even this simplified model governs rather well the desired dynamics. 
But before discussing this, let us introduce the last ingredient of the model -- 
the \modif{damping} forces that we have already mentioned. 

It is clear that without dissipation there is no reason for the system to stay in the 
relaxed configuration: even non-intersecting the inclusions will have non-zero velocity 
and can collide again. To deal with this fact we introduce dissipative forces to the 
system. We will consider two \modif{damping} models: the usual viscous one and the 
so-called mechanical thermostats. The first one amounts simply to adding 
a force proportional to the velocity or the angular velocity 
of the body with a negative constant prefactor $-\beta$. The second one is characterized
by a non-linear \modif{damping} force which is worth being commented on in more details. 
The idea of introducing mechanical thermostat comes directly from molecular 
dynamics or more specifically from the simulation of molecular systems 
at constant temperature. The key point is to introduce a \modif{damping} force of the 
form $- \gamma \bm v$, where in contrast to ordinary viscous \modif{damping} the coefficient  
$\gamma$ depends on the temperature (i.e. the kinetic energy) of the whole system.
In the Berendsen thermostat (\cite{ber}) the prefactor is computed explicitly using the formula
\begin{equation}
  \gamma_{Ber} =  \alpha_{Ber}(E_{kin} - \frac{1}{2} N k_B T),
\end{equation}
where $\alpha$ is a generally small constant coefficient, $E_{kin}$ -- the kinetic energy 
of the system, $N$ -- number of degrees of freedom, $k_B$ -- Boltzmann's constant, 
$T$ -- desired temperature; the expression in brackets corresponds thus to the difference 
between the actual kinetic energy and its value corresponding to the temperature $T$.
The idea is that when the energy is high, $\gamma$ is positive and the force 
slows (cools) down the system, if on the contrary the  energy is low, 
$\gamma$ is negative and the force accelerates (heats up) the system. 
In the Nos\'e--Hoover thermostat (\cite{nose,hoover}) the logic is rather similar, but 
the coefficient is defined by the differential equation:
\begin{equation}
  \dot \gamma_{NH} =  \alpha_{NH}(E_{kin} - \frac{1}{2} N k_B T).
\end{equation}
It has been shown (\cite{golo-shaitan1}) 
that the Berendsen thermostat indeed brings the system to the desired temperature 
\modif{(i.e. $E_{kin}$ approaches $\frac{1}{2} N k_B T$)}, and 
does it exponentially fast. But it fails to reproduce the correct energy distribution 
against degrees of freedom (\cite{golo-shaitan2}), namely the collective (corresponding to 
global translation or rotation as a rigid system) 
degrees of freedom are overheated, while the others are frozen.
In the Nos\'e--Hoover model energy oscillations and resonance effects have been 
observed (\cite{gss}). All these effects make it difficult to apply these thermostats 
in realistic molecular simulations\modif{;} we can however profit from them for our purposes. 
Indeed, heating up the collective degrees of freedom by Berendsen thermostat will precisely 
mean that the inclusions tend not to intersect, otherwise the distance 
between them will oscillate. 
And oscillatory regimes of the Nos\'e--Hoover thermostat can help to spontaneously 
heat up the system to reshuffle it. Or even simpler one can consider these models 
at zero temperature 
to freeze the system fast close to the relaxed configuration. 
It turns out that a combination of these approaches can lead to a more efficient algorithm
of relaxation.
The full model thus includes the viscous \modif{damping}, the Berendsen thermostat that are 
always present and the Nos\'e--Hoover one which is eventually ``switched on'' when the 
relaxation stagnates. The system of ODEs governing the model reads:
\begin{equation} \label{system}
\left\{\begin{array}{lll}
\dot \bm  p_{s_i}&=& \bm  v_{s_i}, \quad \dot \bm p_{c_k} = \bm v_{c_k},  \quad
 \dot \bm l_{c_k}= \boldsymbol{\omega}_{c_k} \times \bm l_{c_k}, \nonumber \\[0.5em]
\dot \bm  v_{s_i}&=& \sum\limits_j \bm F_{ji} - \beta \bm  v_{s_i} - (\gamma_{Ber} + \gamma_{NH}) \bm  v_{s_i}, \nonumber \\[0.5em]
\dot \bm v_{c_k}&=& \sum\limits_{j} \bm F_{jk} - \beta \bm  v_{c_k} - (\gamma_{Ber} + \gamma_{NH}) \bm  v_{c_k},  \nonumber \\[0.5em]
 \dot {\boldsymbol{\omega}}_{c_k} &=& \frac{6}{ \|\bm l_{c_k}\|^2 }\sum\limits_{j} (\bm p_{jk} -\bm p_{c_k}) \times \bm F_{j_k} - \nonumber \\[0.5em]
 & & \quad - \beta \boldsymbol{\omega}_{c_k} - (\gamma_{Ber} + \gamma_{NH}) \boldsymbol{\omega}_{c_k},    
\end{array} \right.
\end{equation}
where $\bm F_{ji}$  is the $j$-th force acting on the $i$-th sphere, and 
$\bm F_{jk}$  is the $j$-th force acting on the $k$-th cylinder;
for the system of $n_s$ spheres and $n_c$ cylinders the number of 
degrees of freedom in the definition of thermostats $N = 3n_s +5n_c$.
We have put all the masses equal to $1$, although one can choose another 
convention, for instance the mass and the \modif{damping} coefficient $\beta$ can depend on 
the volume of the inclusion.

\begin{figure}[ht] 
\centering
\subfigure[\, Dependence of the overlap volume and the introduced potential energy on time]{    
\includegraphics[trim = 2cm 1.5cm 6.5cm 19cm, clip, width=0.97\linewidth]{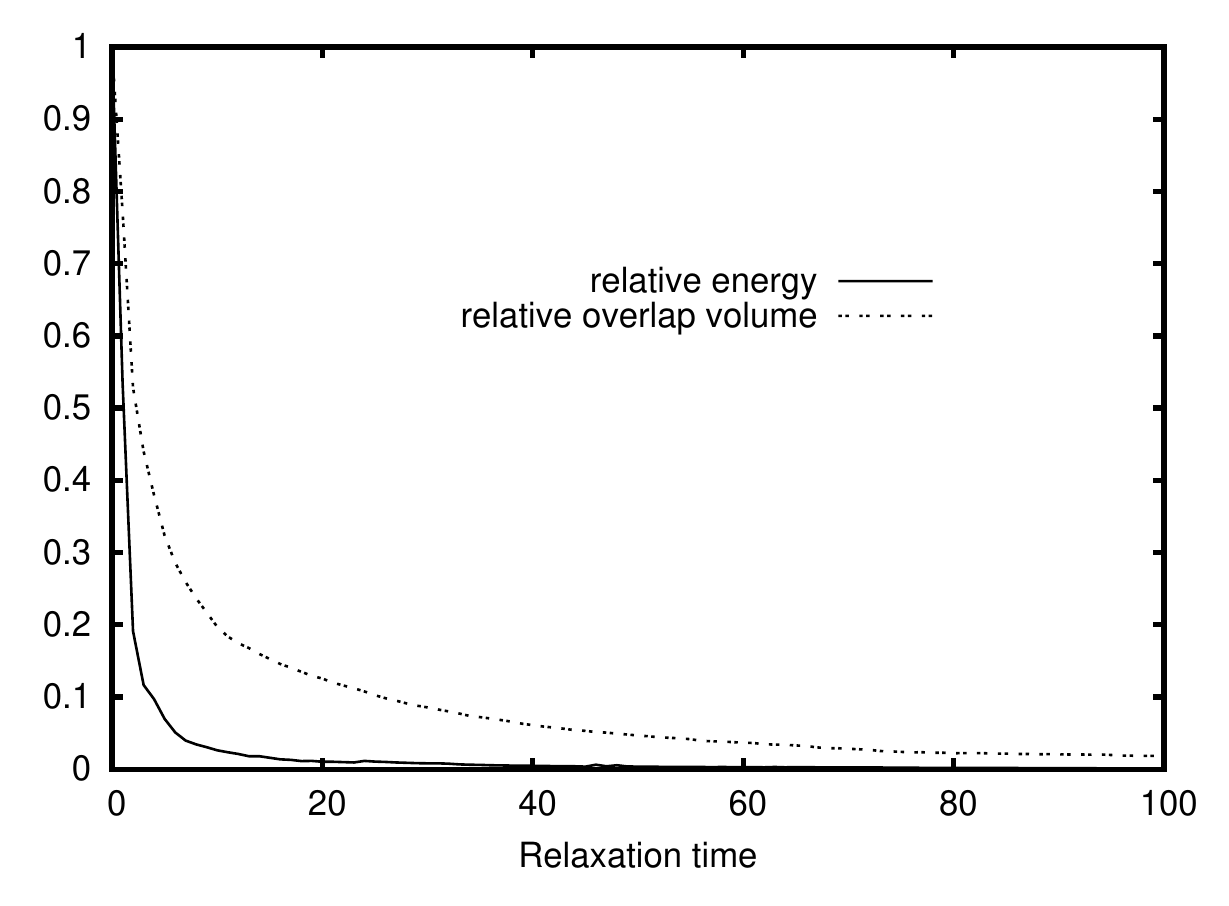}  
}
\subfigure[\, Relation between the (normalized) overlap volume and the defined energy]{
    \includegraphics[trim = 1.5cm 1.5cm 6.5cm 19cm, clip, width=0.97\linewidth]{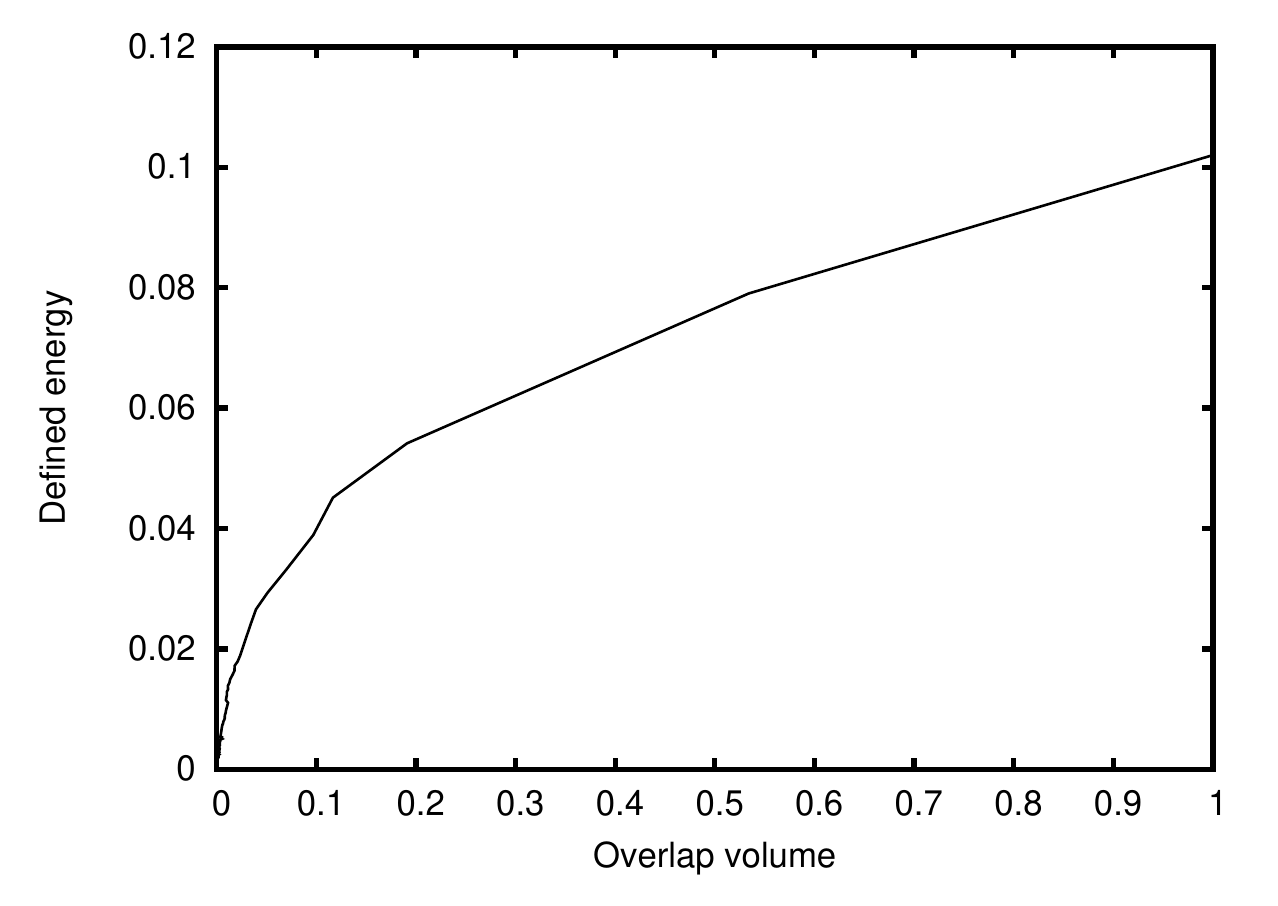}  
 }
\caption{The defined energy reflects well the overlap volume.}
\label{fig:stop-cond}
\end{figure}

The last but not the least point in description of this method is the condition
to stop the simulation. 
It is clear that a totally relaxed configuration corresponds to 
the volume occupied by all the inclusions equal to the desired one.
Computation of the overlap volume at each integration step is rather long. 
To overcome this difficulty we introduce the notion of potential energy 
which is basically the sum of squares of norms of the internal forces
(table \ref{tab:forces}), that is it uses the quantities computed at each step anyway.
The simulation is stopped when this energy is lower than some critical value, that 
we advise to determine for each family of simulations by performing several 
test runs.
The figure \ref{fig:stop-cond} shows that this energy is perfectly correlated with the 
real value of the overlap volume.

Having now the full model we have performed the above mentioned validation campaign
in order to make sure that the dynamics governed by the simplified mechanical equations
(\ref{mech}) is acceptable. 
The figures (\ref{fig:sp-sp} -- \ref{fig:cyl-cyl9}) in the Appendix 1 
show 
typical intersections of couples of spheres and cylinders and their 
configuration after the relaxation process governed by the 
equations (\ref{system}). As one sees the dynamics is precisely as expected 
from the physics of the model (see the figure captions for more details in each 
case).

Let us make several remarks concerning the implementation of the above relaxation method. 
As we have already mentioned we are not interested in high precision for the relaxation 
trajectory, as long as it does lead to the relaxed configuration. There is thus no need 
in involving advanced algorithms of numerical integration of ODEs: Leapfrog or velocity Verlet, 
or even trapezoid method already do the job rather well. 
Like in the previous section, we should not forget about the boundary conditions 
for the dynamics. To take them into account, an attribute of being close to 
the boundary is assigned to inclusions, and depending on it the periodic images 
are also taken into account to compute interaction forces. This attribute now has to be updated 
at each integration step.

\begin{figure}[ht] 
\centering
\includegraphics[trim = 1.5cm 1.5cm 6.5cm 19cm, clip, width=0.97\linewidth]{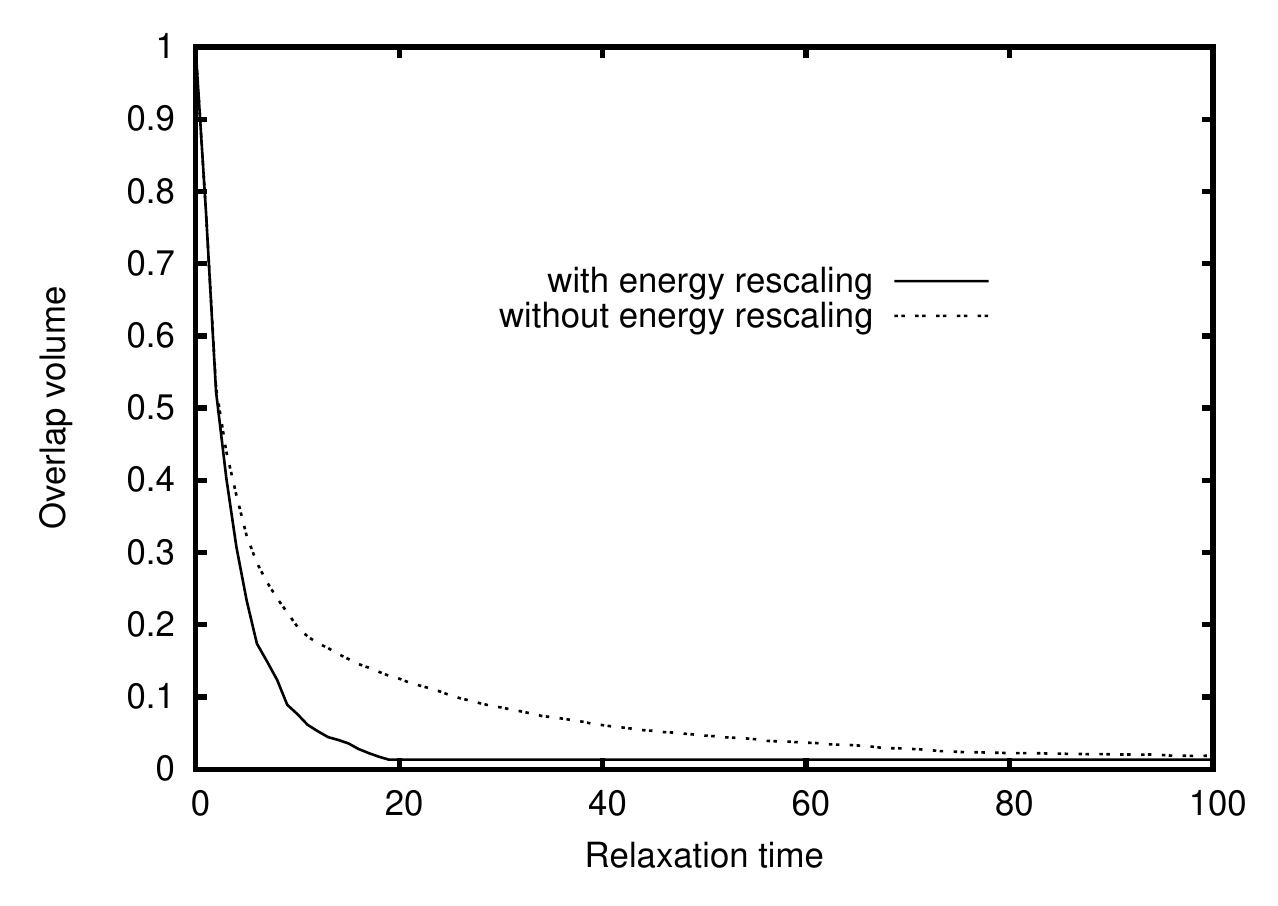}  
\caption{Acceleration of convergence by introducing energy rescaling.}
\label{fig:force-factor}
\end{figure}

Before turning to concrete computational results, let us also note that 
the described process allows one to introduce various numerical tricks, like 
reshuffling by the Nos\'e--Hoover thermostat, that we have already mentioned. 
Another one, which is useful for applications, is related to the fact that 
close to the relaxed configurations, the introduced forces necessarily 
become small, that slows down the process. A natural way to speed it up 
is to rescale the forces when the energy decreases. 
\modif{In practice, we introduce a global prefactor for all the forces, which 
is increased every time when the energy decreases for an example by a factor of two.
This allows one to have reasonable gradients in the beginning of the relaxation process 
and terminate the process rather fast.}
The figure 
\ref{fig:force-factor} shows the acceleration of the relaxation process due to this 
idea. 

As in the previous section, we have performed the tests to study the efficiency
of the suggested method. The time of generation at low volume fractions is of the order
$10^{-1}$ seconds, which is slower than with the RSA-type algorithm. But the MD method 
allows us to achieve volume fractions of $50\%$ and even more within several seconds, 
while for such values the RSA method almost never produces a result. 
The tendency is also roughly the same: spheres are easier to generate than cylinders, 
and it is more difficult to treat higher aspect ratio. We should note however, 
that there is no pronounced saturation effect with the MD method, i.e. 
one can achieve high volume fractions, conjecturally up to the theoretical 
limits. The tables \ref{tab:MD10} -- \ref{tab:MD50} from the appendix 3 show these 
dependencies in more details.

\newpage
\clearpage
\section{Efficiency comparison} \label{sec:compare}
From the previous two sections one can already get the idea that 
the RSA algorithm is much more efficient for low volume fraction but the 
time driven MD allows one to achieve higher fractions, for which RSA stagnates:
it is sufficient to compare the upper-left corners of tables \ref{tab:MD10} -- 
\ref{tab:MD50} with \ref{tab:RSA10} -- \ref{tab:RSA50}. 
\modifnew{This effect is perfectly understandable: every elementary step of the RSA method 
consists of the algorithm \ref{alg:sp-cyl} or \ref{alg:cyl-cyl} to be applied to the newly generated 
object and all the objects memorized before, while for MD the same algorithms are to 
be applied to all the \emph{couples} of objects at each timestep. For relatively low volume fraction 
the collisions in the RSA method appear rather rarely, while the MD method still needs to make 
several steps to converge. When the volume fraction increases the number of rejected inclusions for 
the RSA grows fast and at some point the total computation time becomes greater than for the MD.}
The figures below show this comparison for some cases in more details.

\begin{figure}[H] 
\centering
\includegraphics[trim = 1.5cm 1.5cm 6.5cm 19cm, clip, width=0.97\linewidth]{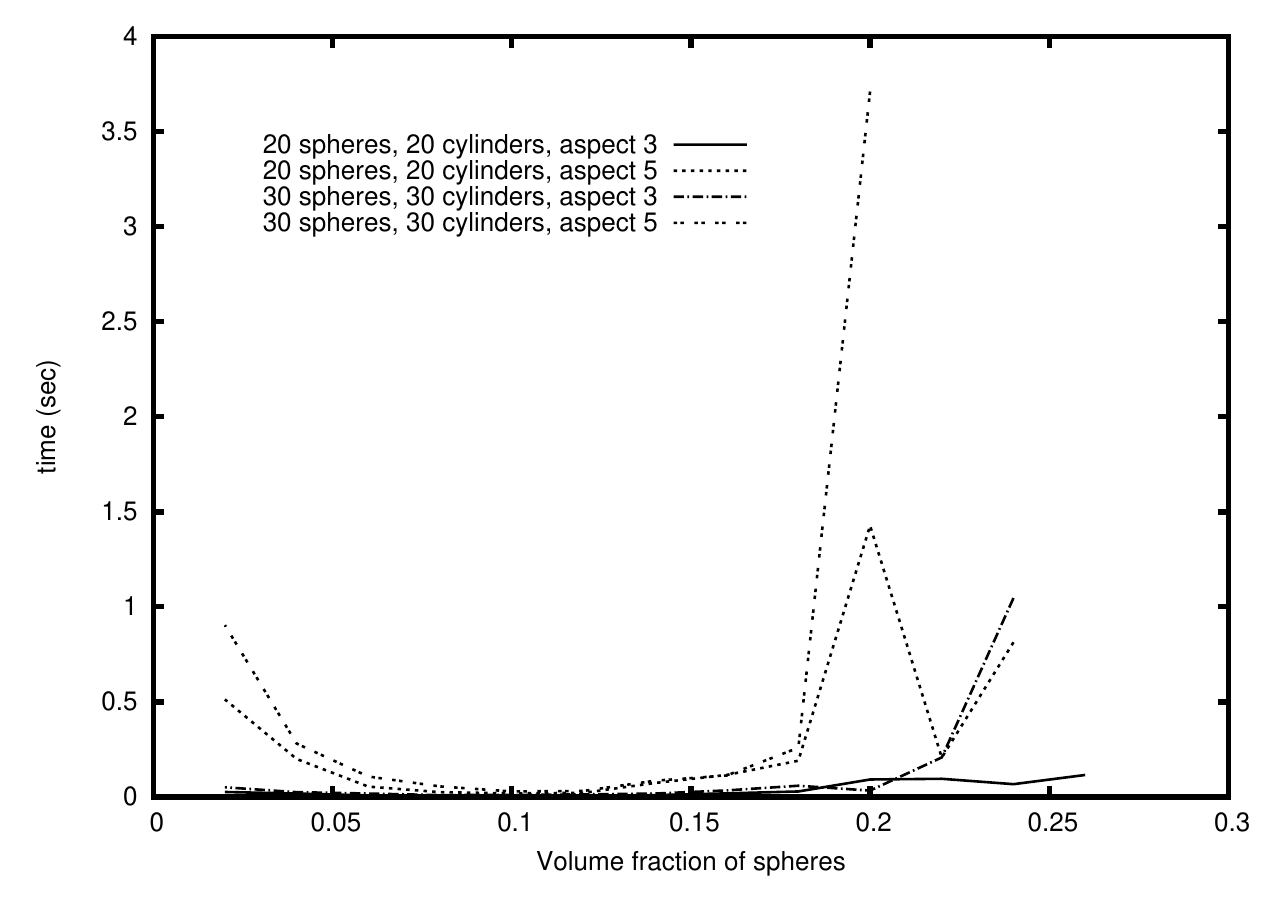}  
\caption{\label{fig:RSA030} RSA: \modif{total volume fraction fixed at $0.3$, 
$f_{sp}$ varies with the step of $0.01$, $f_{cyl} = 0.3 - f_{sp}$}.}
\end{figure}

\begin{figure}[H] 
\centering
\includegraphics[trim = 1.5cm 1.5cm 6.5cm 19cm, clip, width=0.97\linewidth]{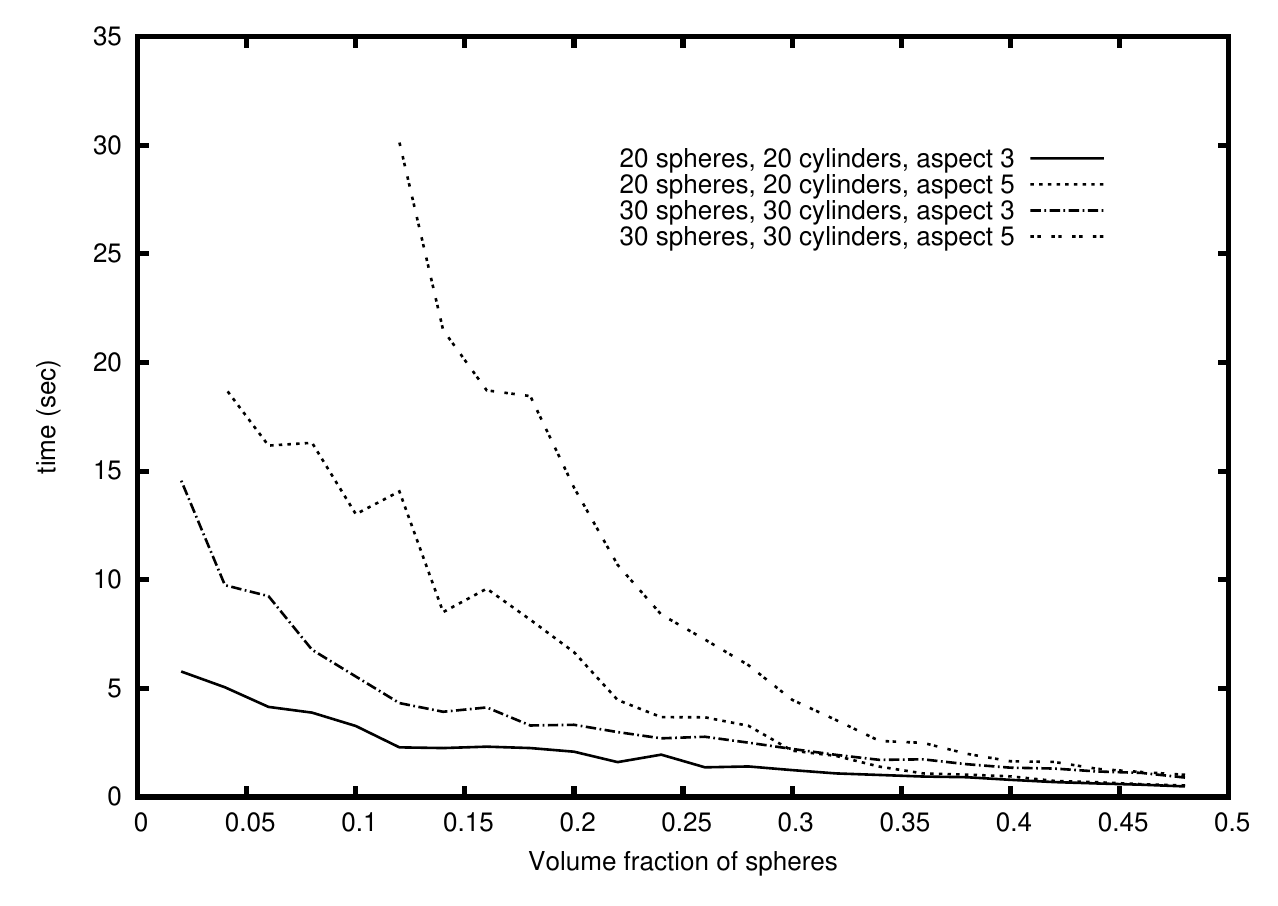}  
\caption{\label{fig:MD050} MD: \modif{total volume fraction fixed at $0.5$, 
$f_{sp}$ varies with the step of $0.01$, $f_{cyl} = 0.5 - f_{sp}$}.}
\end{figure}

First let us fix the volume fraction ($0.3$ for RSA and $0.5$ for MD) and 
see how the time of generation depends on the distribution of this fraction 
between spheres ($f_{sp}$) and cylinders ($f_{cyl}$). The figures (\ref{fig:RSA030}) and (\ref{fig:MD050})
show that for MD it is clearly more difficult to treat the cylinders, and for RSA 
the same phenomenon takes place for high proportion of cylinders.
The same saturation effect as before is observed for high proportion of spheres, 
which can be handled by changing the generation strategy.

\begin{figure}[H] 
\centering
\includegraphics[trim = 1.5cm 1.5cm 6.5cm 19cm, clip, width=0.97\linewidth]{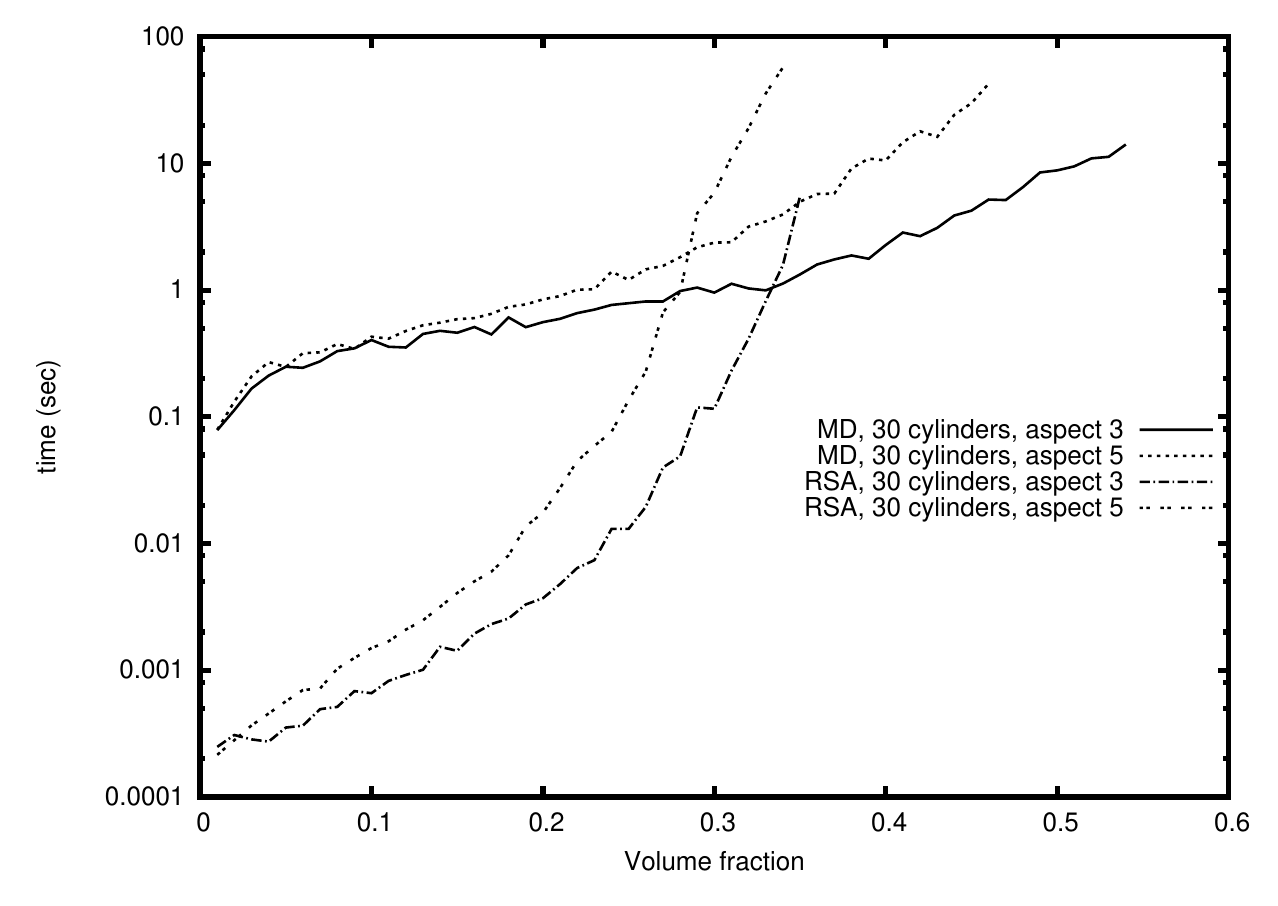}  
\caption{\label{fig:cyls_only} Cylinders only:  \modif{$f_{cyl}$
 varies with the step of $0.01$}}
\end{figure}

\begin{figure}[H] 
\centering
\includegraphics[trim = 1.5cm 1.5cm 6.5cm 19cm, clip, width=0.97\linewidth]{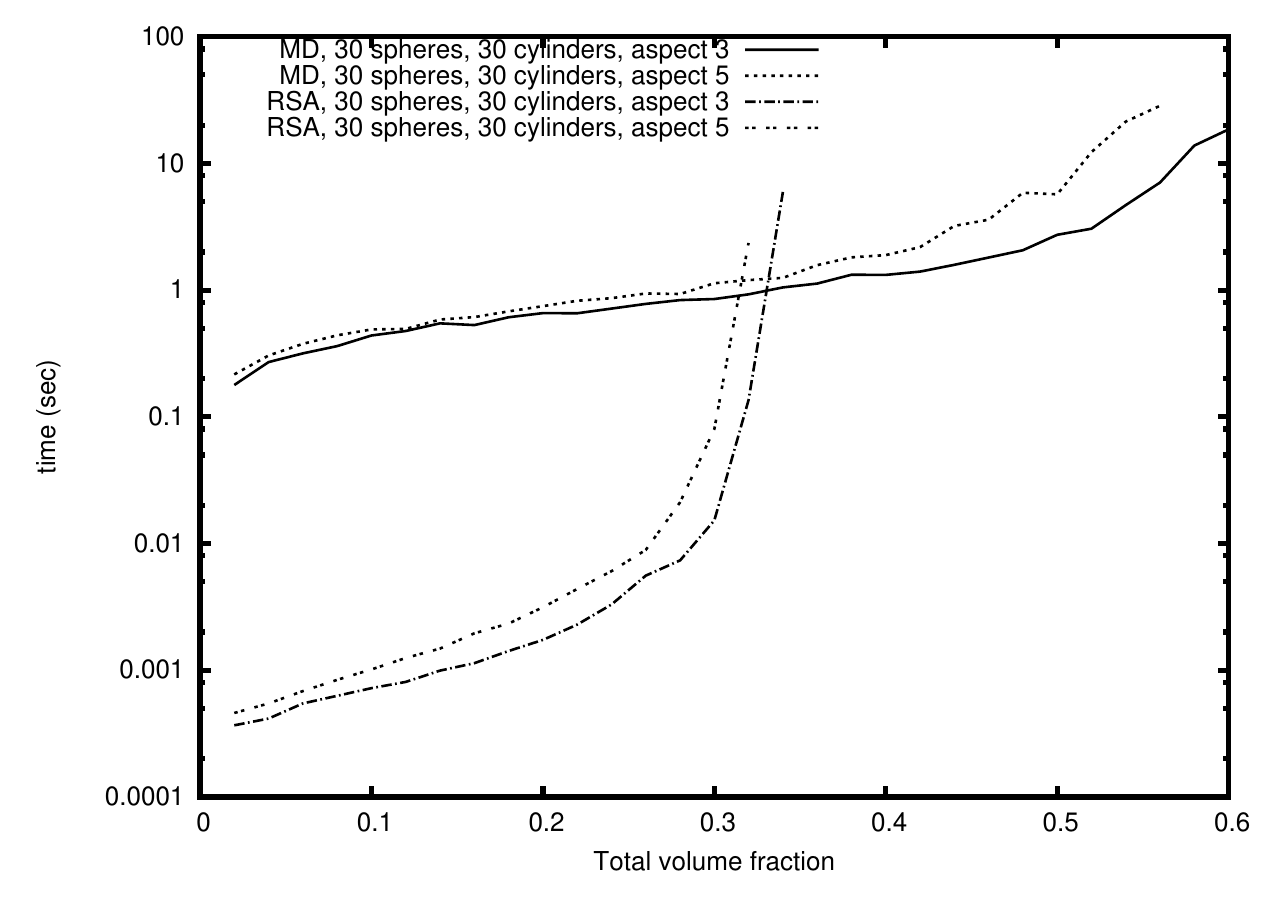}  
\caption{\label{fig:cyl_sp} Cylinders and spheres: \modif{$f_{cyl}$ and $f_{sp}$
 vary with the step of $0.01$}}
\end{figure}

\begin{figure}[H] 
\centering
\includegraphics[trim = 1.5cm 1.5cm 6.5cm 19cm, clip, width=0.97\linewidth]{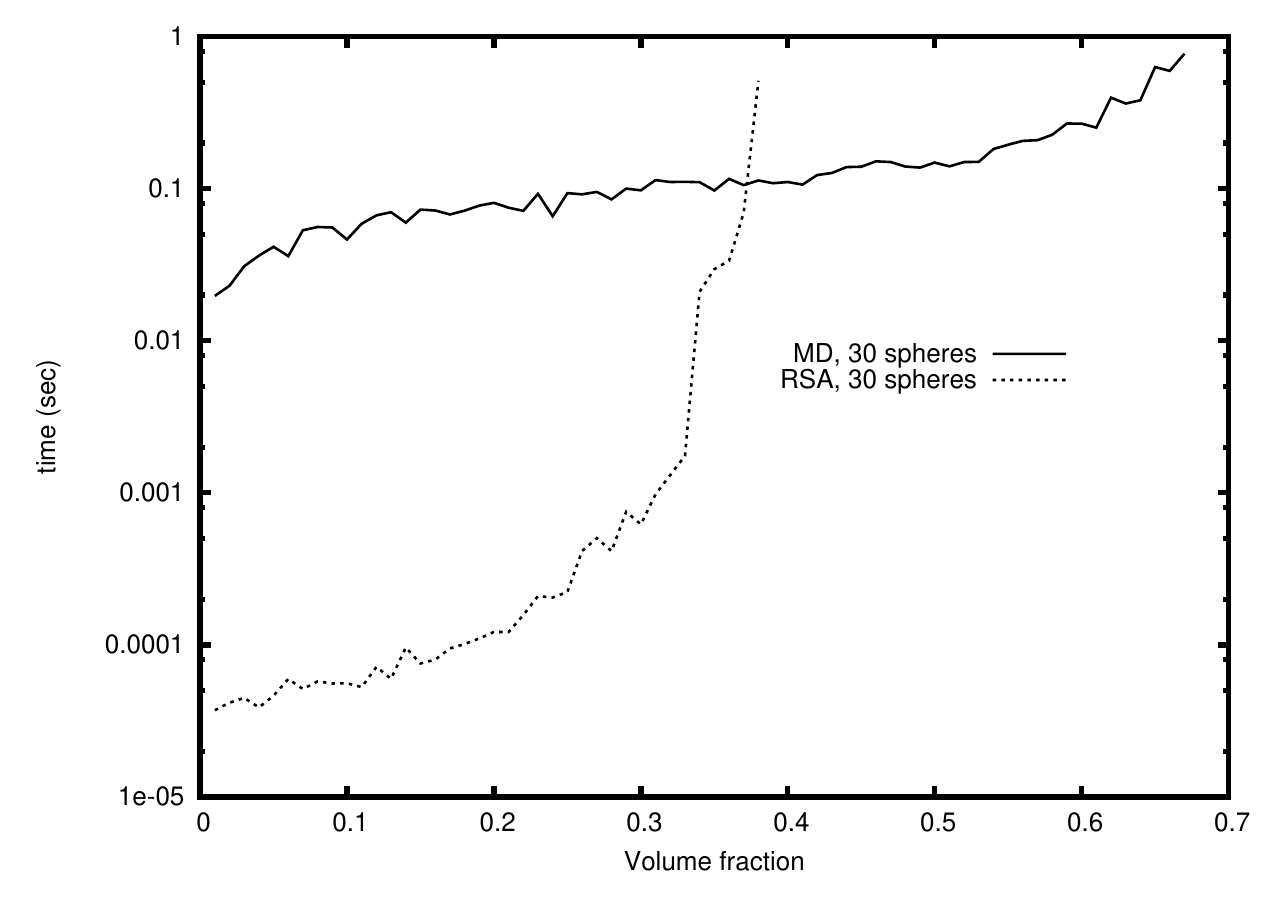}  
\caption{\label{fig:sp_only} Spheres only, \modif{$f_{sp}$
 varies with the step of $0.01$}}
\end{figure}

The following pictures (\ref{fig:cyls_only} -- \ref{fig:sp_only}) show the 
dependence of the generation time (as usual averaged over 20 runs) on the total 
volume fraction for spheres, cylinders, and their mixture. From them we clearly see the advantage 
of the RSA algorithm for small volume fractions as well as the capabilities of the MD method 
for higher ones. We have mentioned before that generating cylinders with higher aspect ratio
$a$ is more complicated for both algorithms. It is pretty obvious that for 
very high aspect ratio it is difficult even for rather low volume fraction of inclusions -- 
we indeed observe this effect in the tests. To give an example, consider the  generation 
of the RVE by the RSA technique with a mixture of spheres and cylinders at total volume fraction of 
$10 \%$. At the aspect ratio $a = 5$ it takes about $2 \cdot 10^{-3}$ sec.,  
for $a = 25$ the time is close to $10^{-2}$ sec., for $a = 50$ it approaches $0.1$ sec., and
for $a = 100$ it exceeds $10$ sec.
We should mention that the reason of the effect is twofold. First,  
with long cylinders it is easier to construct non-acceptable geometries of RVEs, 
especially in combination with a few large spheres. This makes the RSA algorithm get 
stuck and the MD to perform more verification steps. Second, and even more important, 
there is a limitation on the length of cylinders: we find it natural to consider  inclusions 
that are smaller than the studied RVE, in particular an inclusion should not intersect 
its periodic image. This forces the lower bound on the number of cylinders:
$N_{cyl} \ge \frac{4}{\pi}f_{cyl}a^2$. 
One thus needs to apply the algorithm to a much bigger number of 
inclusions ($150$ in the above example for $a = 100$) or review the notion of RVE for high aspect ratio.
 

%

\section{Conclusion/outlook}
In this paper we have presented the methods of generation of RVEs and performed 
the tests of their efficiency. We have seen that the RSA type algorithm allows us 
to generate the configuration extremely fast when the volume fraction is rather small. 
The MD-based algorithm does not have this limitation and permits us to achieve configurations 
even close to theoretical maximum of volume fraction.

\begin{figure}[H]  
\centering
\subfigure{    
\includegraphics[trim = 4cm 7cm 4cm 7cm, clip, width=0.9\linewidth]{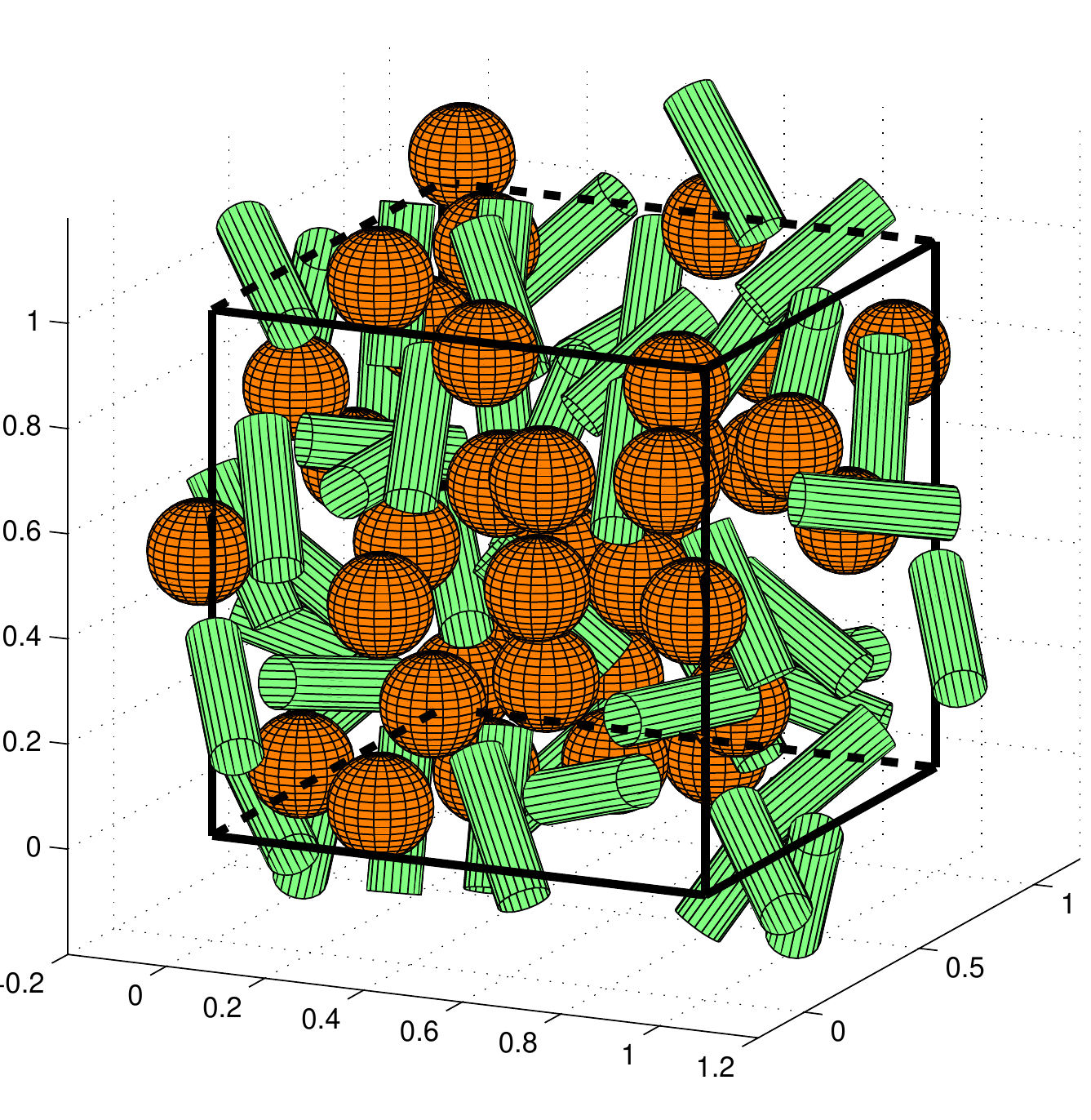}  
}
\subfigure{
      \includegraphics[width=0.45\linewidth]{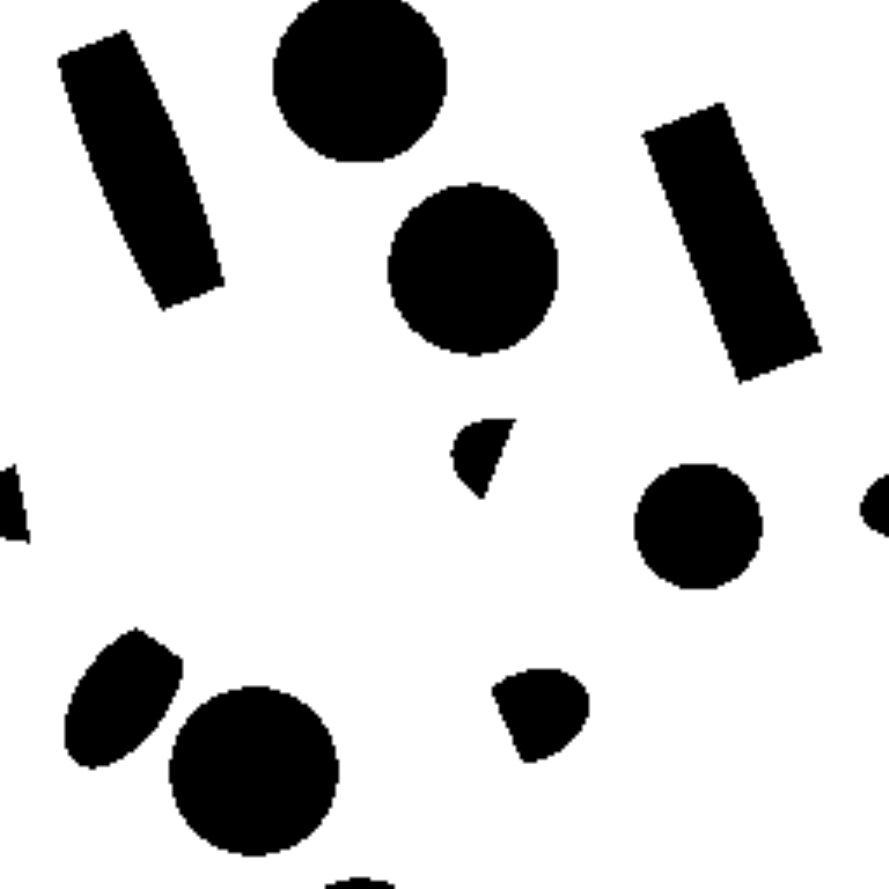}   }
\subfigure{
      \includegraphics[width=0.45\linewidth]{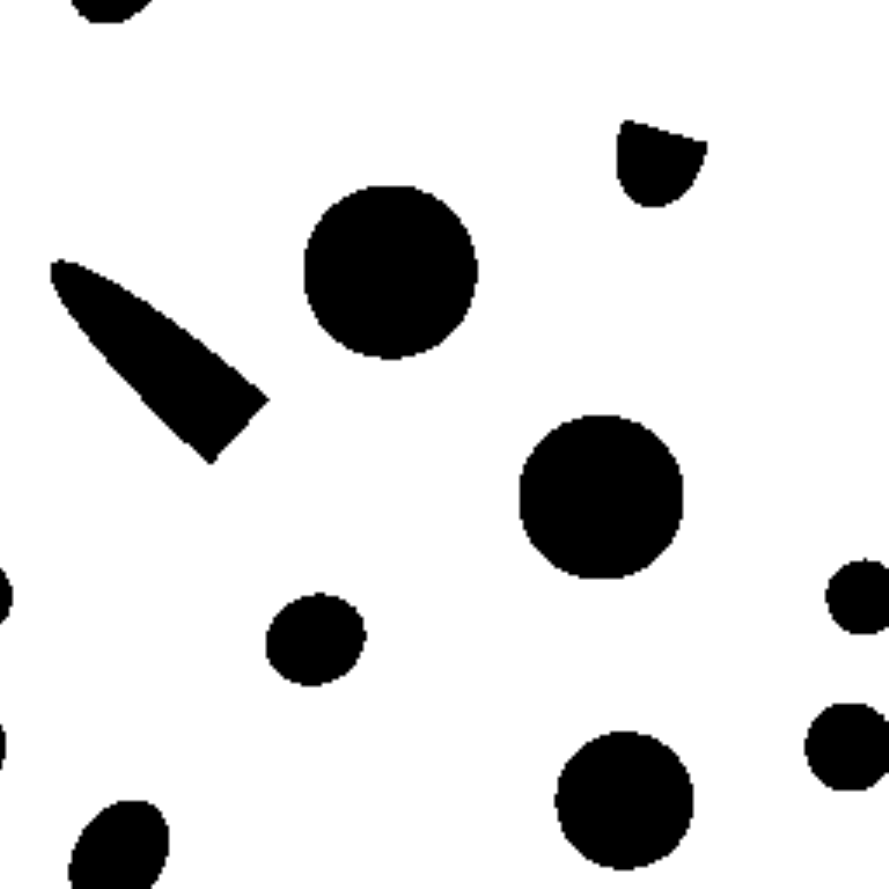}   }
      
\caption{Example of RVE generated by the RSA algorithm at total volume fraction of $20 \%$: 3D view and a couple
of typical slices of the pixelized 3D image}
\label{fig:RVE02}
\end{figure}

\begin{figure}[H]  
\centering
\includegraphics[width=0.8\linewidth]{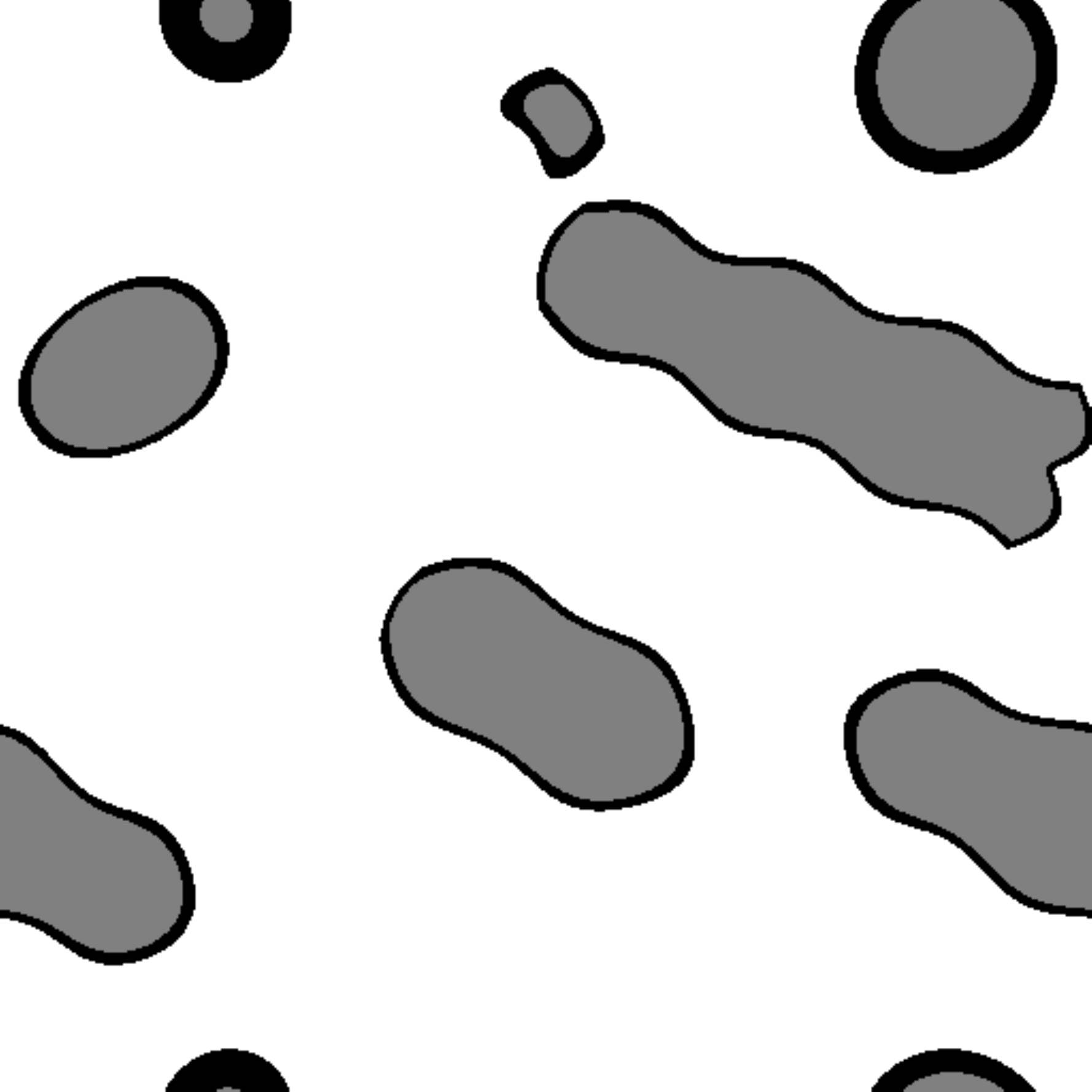}  
\caption{\modif{Example of a section of an RVE with waved inclusions, the interphase regions are highlighted.}}
\label{fig:waved}
\end{figure}

\begin{figure}[ht]  
\centering
\subfigure{    
\includegraphics[trim = 4cm 7cm 4cm 7cm, clip, width=0.9\linewidth]{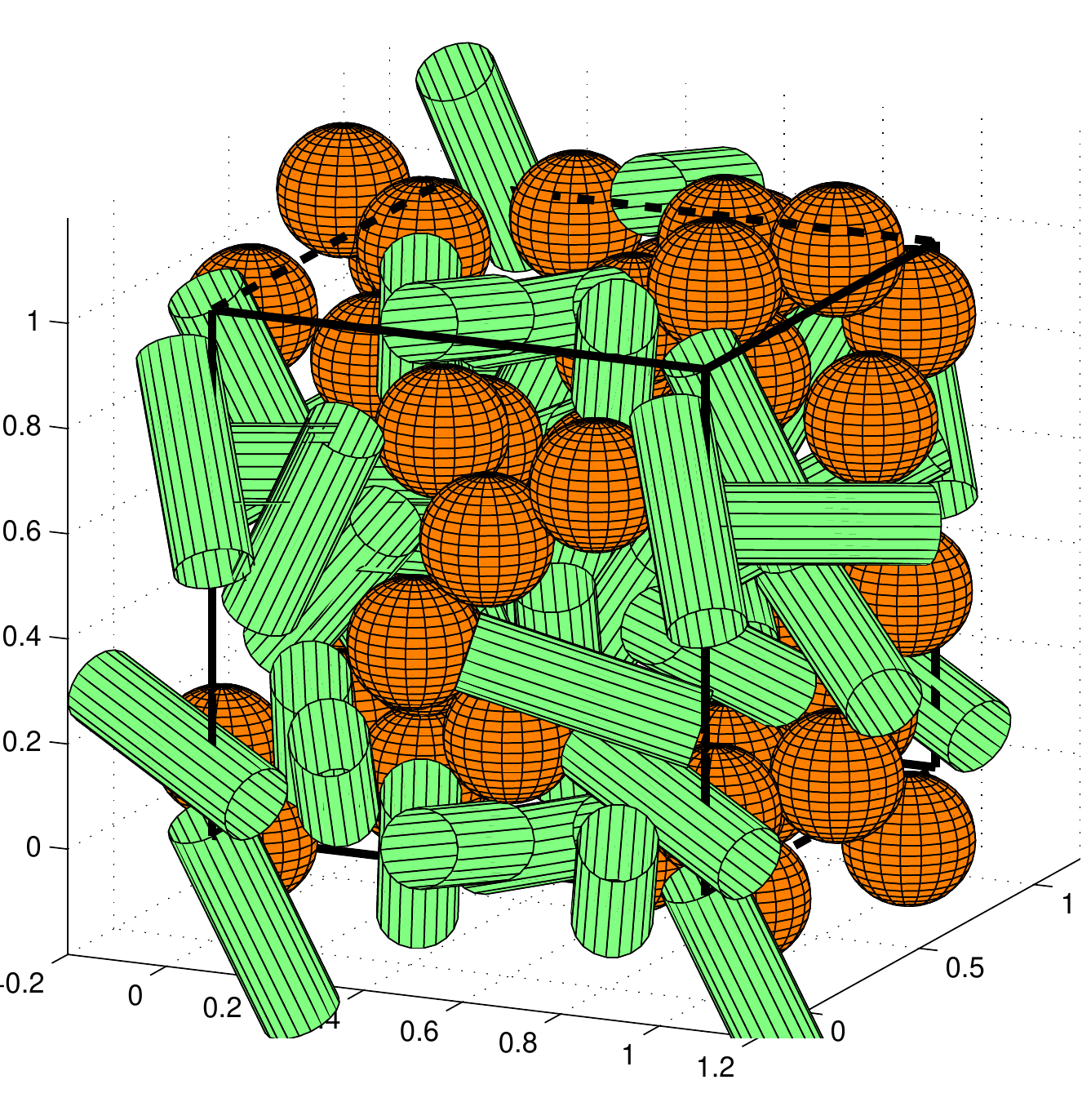}  
}
\subfigure{
      \includegraphics[width=0.45\linewidth]{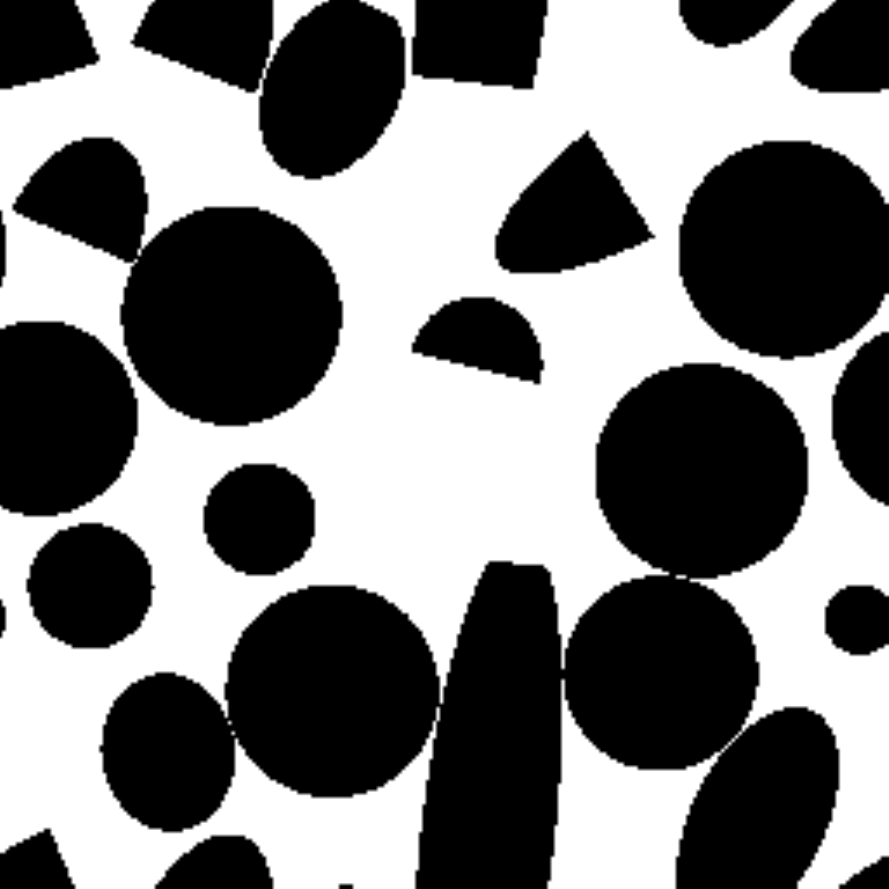}   }
\subfigure{
      \includegraphics[width=0.45\linewidth]{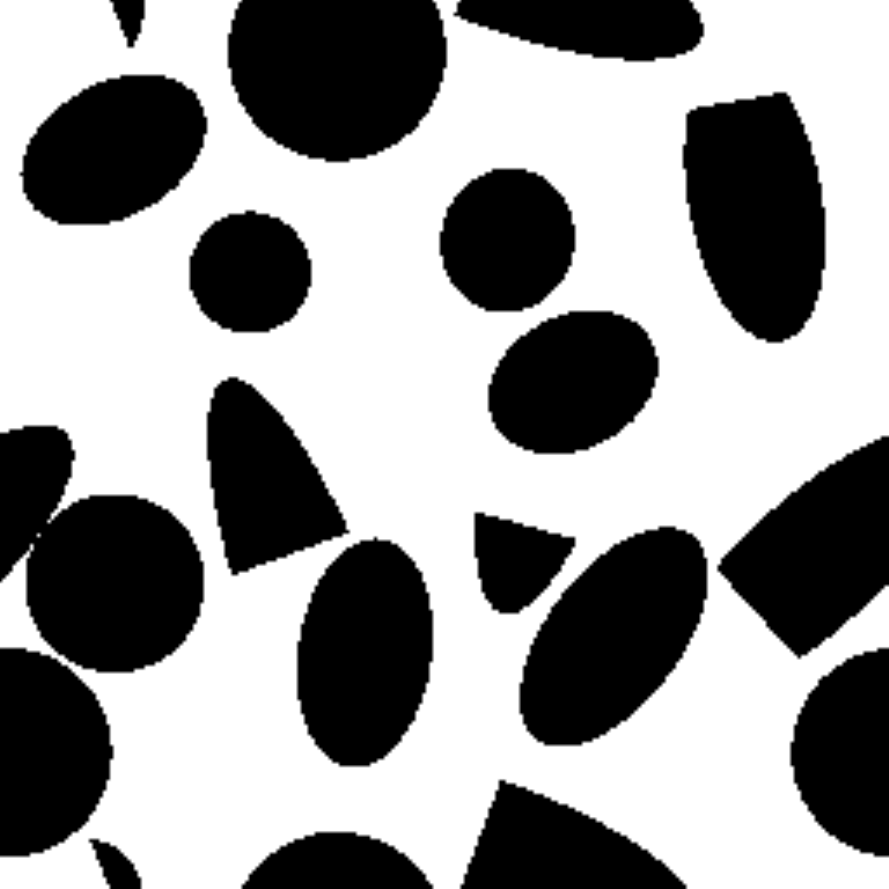}   }
\caption{Example of RVE generated by the MD algorithm at total volume fraction of $40 \%$: 3D view and a couple
of typical slices of the pixelized 3D image}
\label{fig:RVE04}
\end{figure}

As it follows from the title the main motivation for this work is to construct the 
RVE samples for computation of effective properties of composite materials. 
The computational methods that we use in the current work in progress 
are 
based on homogenization techniques using the Fast Fourier Transform\footnote{For the 
implementation of the methods we were inspired by the work \cite{levesque_sp} 
profiting from recent results of \cite{moul-suq}, \cite{monchiet}} or finite element method.
Here we should mention that the methods described in this paper are very well adapted to 
such applications. 
\modif{First, the algorithms are very flexible towards the modifications of the input data. 
For example, instead of using the identical spheres and cylinders as we did for efficiency tests, 
one can introduce any deterministic or probabilistic law for the parameters of their 
geometry. The computation that we perform in order to study the effective properties 
of composites shows that such modifications do not spoil the efficiency of the methods.}
\modif{Second}, the format of the output of the algorithms 
is very convenient for further usage in the computations, namely having all the information 
about the RVE encoded in the concise vector form we can on the one hand pixelize it to 
have a natural discretization of the analyzed sample, and on the other hand keep track of 
orientations of the inclusions (see figures \ref{fig:RVE02}, \ref{fig:RVE04}). 
\modif{In the context of composite analysis, pixelization is also a powerful tool. 
One can for instance introduce deformations of inclusions, detect the boundary and 
assign various properties to it, add defects etc. The figure \ref{fig:waved} shows the section 
of an RVE with waved cylinders and spheres and a highlighted interphase between the matrix and the inclusions.}
Moreover 
from the pixelized RVE we can easily construct a mesh for validation of the computation 
by Finite Elements method (taking the voxels as elements).
\modif{Third}, the main criteria of efficiency that we studied was the time of generation, which is reasonable 
since the methods are not memory consuming. 
To give an idea, let us mention that one computation of the homogenized stiffness 
tensor for a 3D sample
discretized at the resolution $256\times256\times256$ can take several hours for high contrast 
between the properties of the matrix and the inclusions. 
It means that the suggested methods are indeed efficient as the time of generation 
of the sample is negligible in comparison to the time of computation.

The presented methods are very flexible towards various fine-tuning procedures. For example
in the description of the methods we potentially authorized the tangent contacts 
between the inclusions (which is reasonable for our applications), but it is possible 
to avoid them by just modifying the effective interaction distances. 
One can also construct more complex geometries by authorizing some 
types of intersections, which is easy to check since all the information 
of the geometry of the sample is encoded in a concise vector form.
More precisely in the MD-based methods one can impose 
fixed distances or angles between some inclusions to 
produce interesting figures. Because of the natural form of the evolution 
equations (\ref{system}) these restrictions can be implemented 
using a well-developed formalism of mechanical systems with 
constraints and  Lagrange multipliers (\cite{lagrange}).

Let us also note that the presented methods can be useful for a purely mathematical 
purpose of studying the dense packings of simple geometric objects. For example 
with MD we managed to generate a configuration which is very close to the 
periodic one described by Gauss to realize the maximal volume fraction occupied
by identical spheres.

\acknowledgments
We would like to thank Martin L\'evesque for valuable bibliographic information
as well as Elias Ghossein for providing supplementary material related 
to \cite{levesque_sp}.  \newline
This work has been supported by the ACCEA project selected by the ``Fonds Unique Interminist\'eriel
(FUI) 15 (18/03/2013)'' program.

\section*{Appendix 1. Validation of dynamics} \label{sec:MD-validation}
\begin{figure}[H] 
\centering
\subfigure[\, Intersecting]{    
\includegraphics[trim = 5cm 7cm 4.5cm 7cm, clip, width=0.45\linewidth]{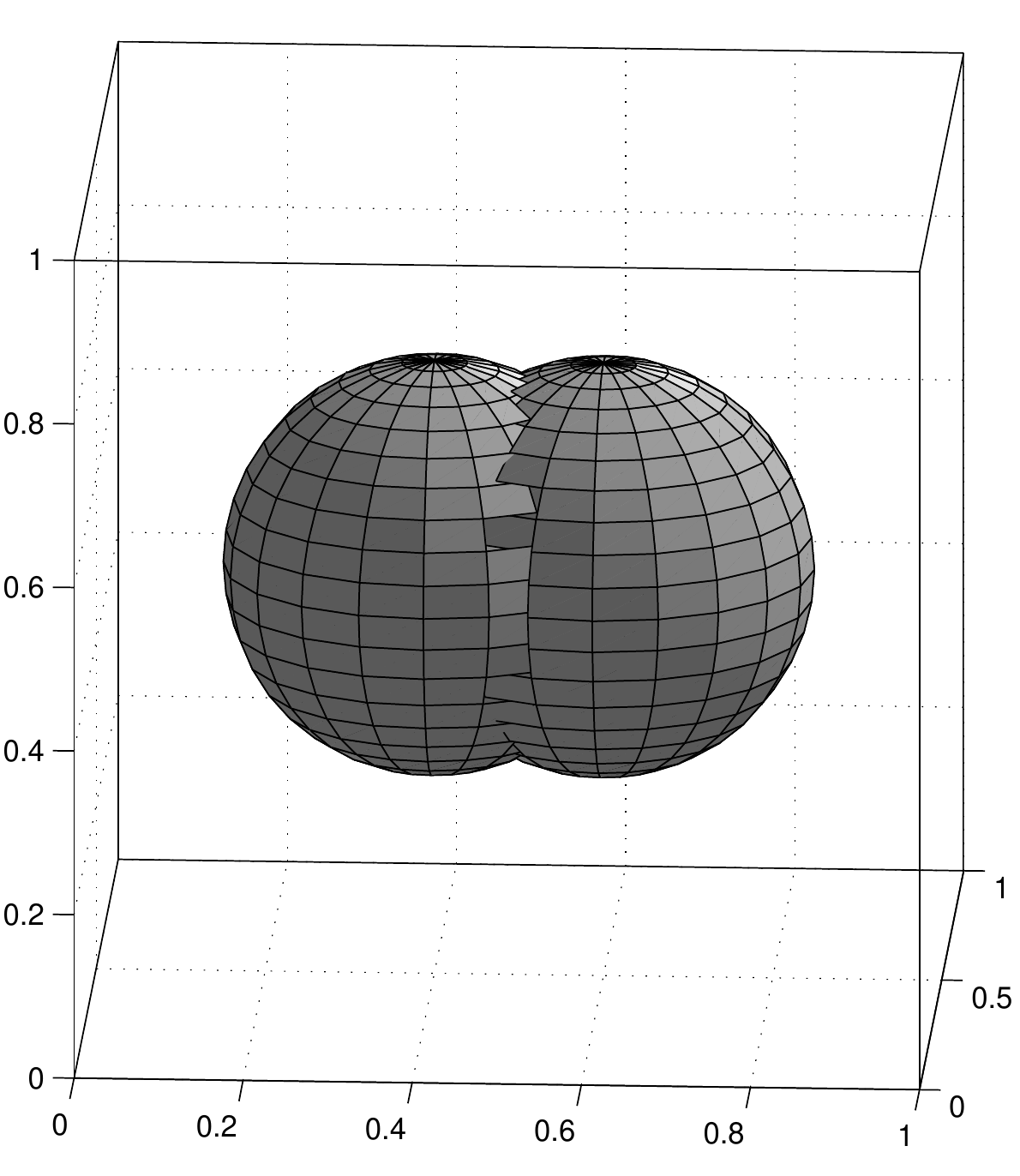}  
}
\subfigure[\, After relaxation]{
    \includegraphics[trim = 5cm 7cm 4.5cm 7cm, clip, width=0.45\linewidth]{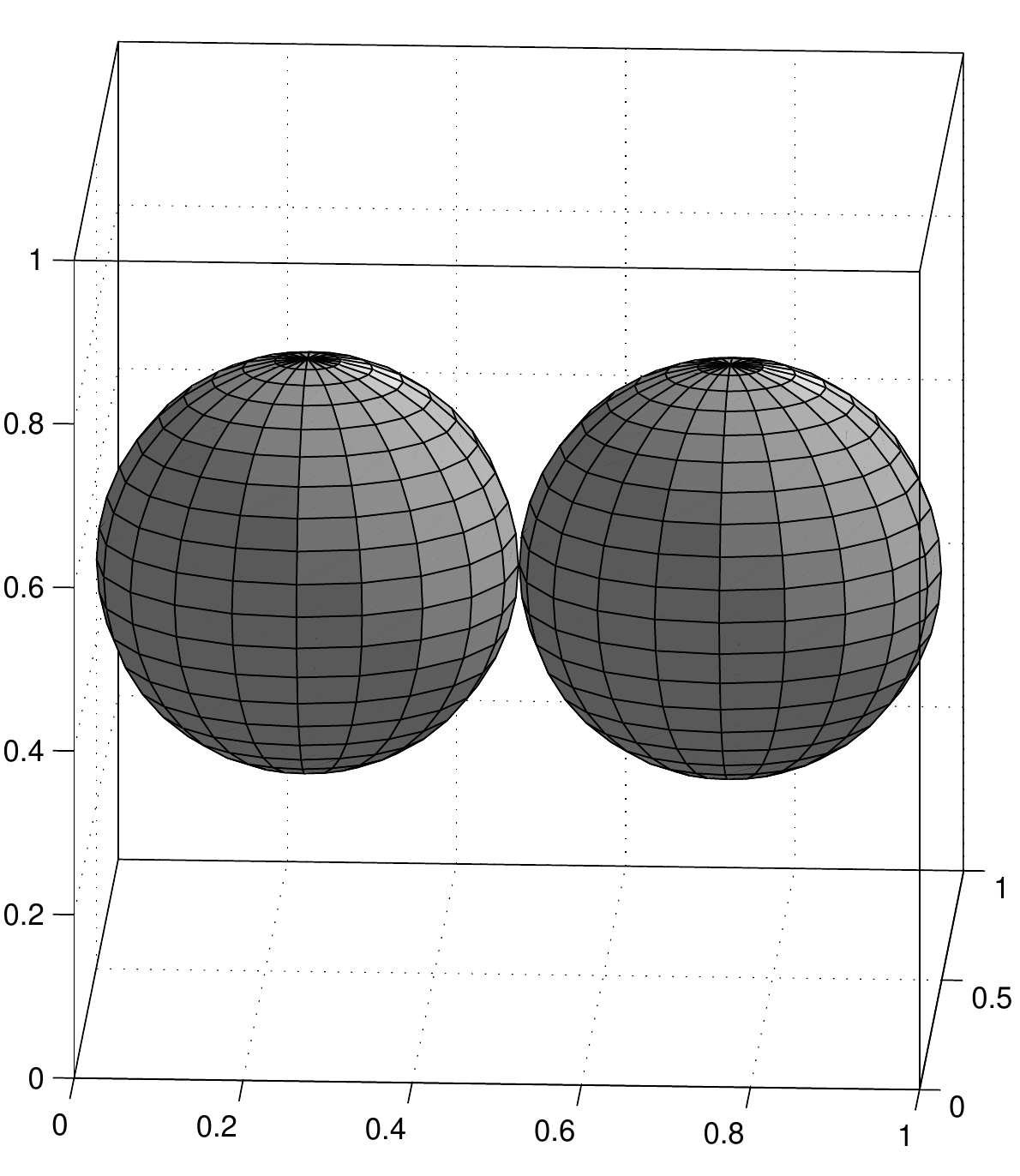}  
 }
\caption{Two spheres.}
\label{fig:sp-sp}
\end{figure}

\begin{figure}[H] 
\centering
\subfigure[\, Intersecting]{    
\includegraphics[trim = 5cm 7cm 4.5cm 7cm, clip, width=0.45\linewidth]{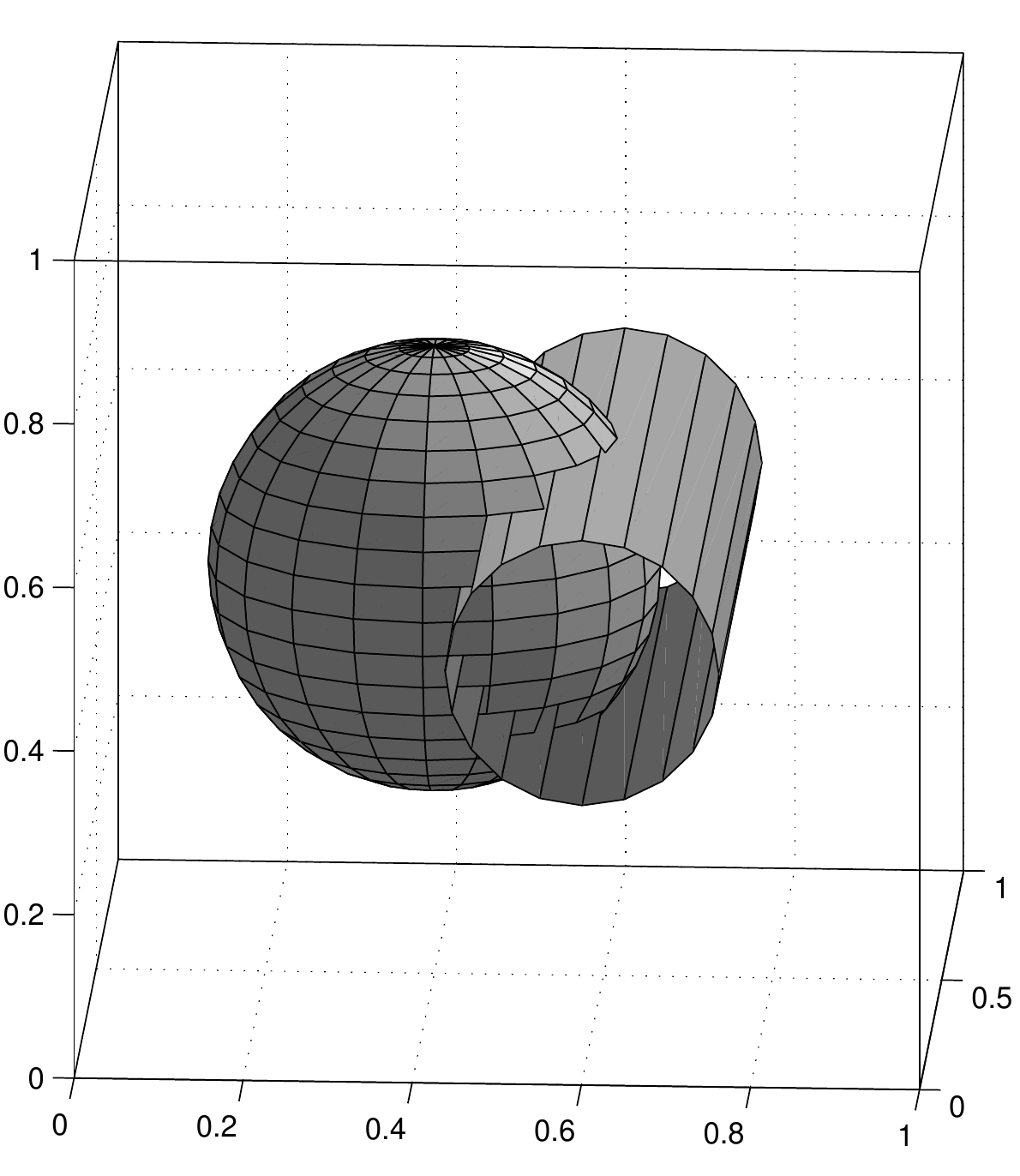}  
}
\subfigure[\, After relaxation]{
    \includegraphics[trim = 5cm 7cm 4.5cm 7cm, clip, width=0.45\linewidth]{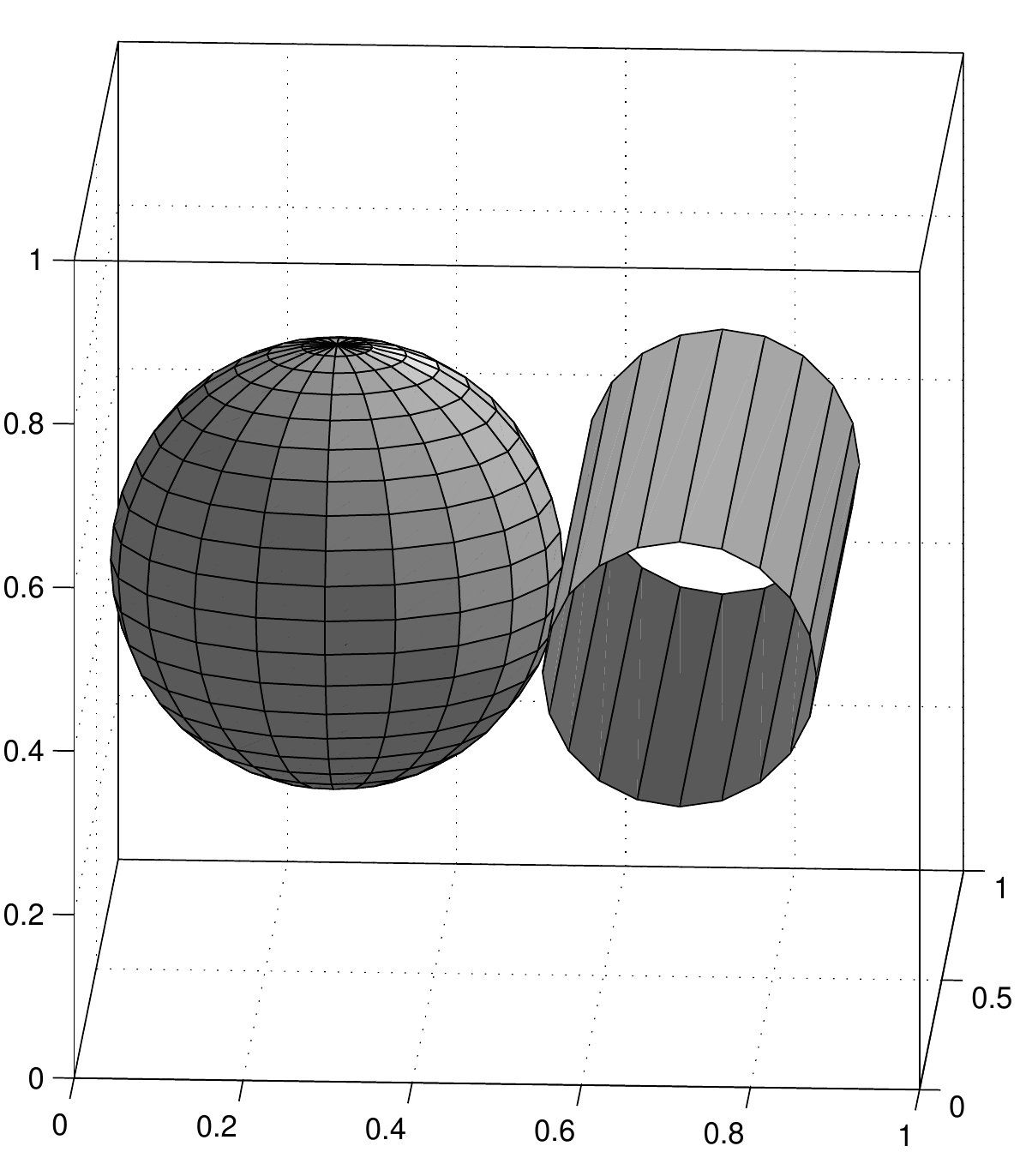}  
 }
\subfigure[\, Intersecting, top view]{    
\includegraphics[trim = 4.5cm 7cm 3.5cm 7cm, clip, width=0.45\linewidth]{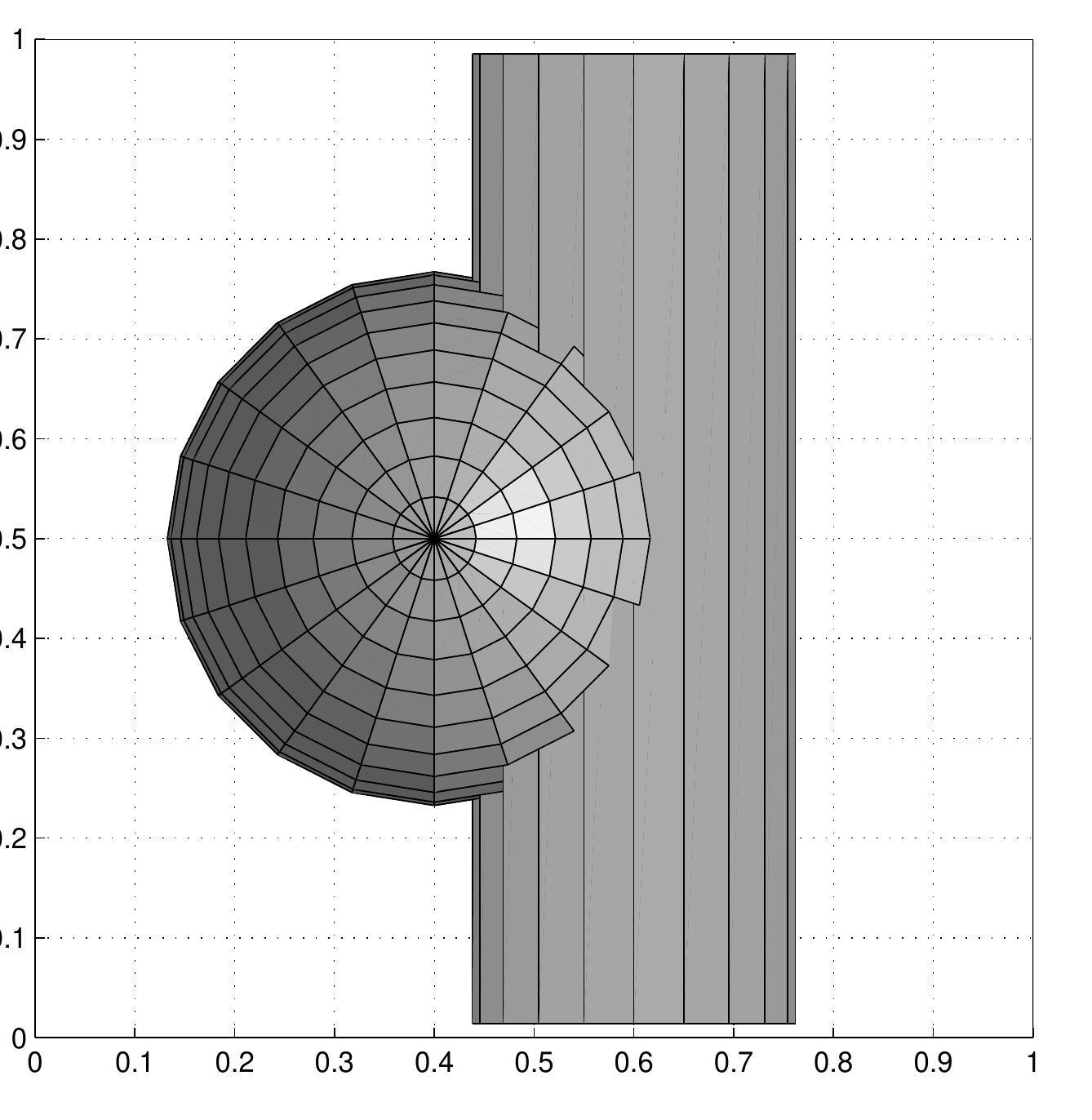}  
}
\subfigure[\, After relaxation, top view]{
    \includegraphics[trim = 4.5cm 7cm 3.5cm 7cm, clip, width=0.45\linewidth]{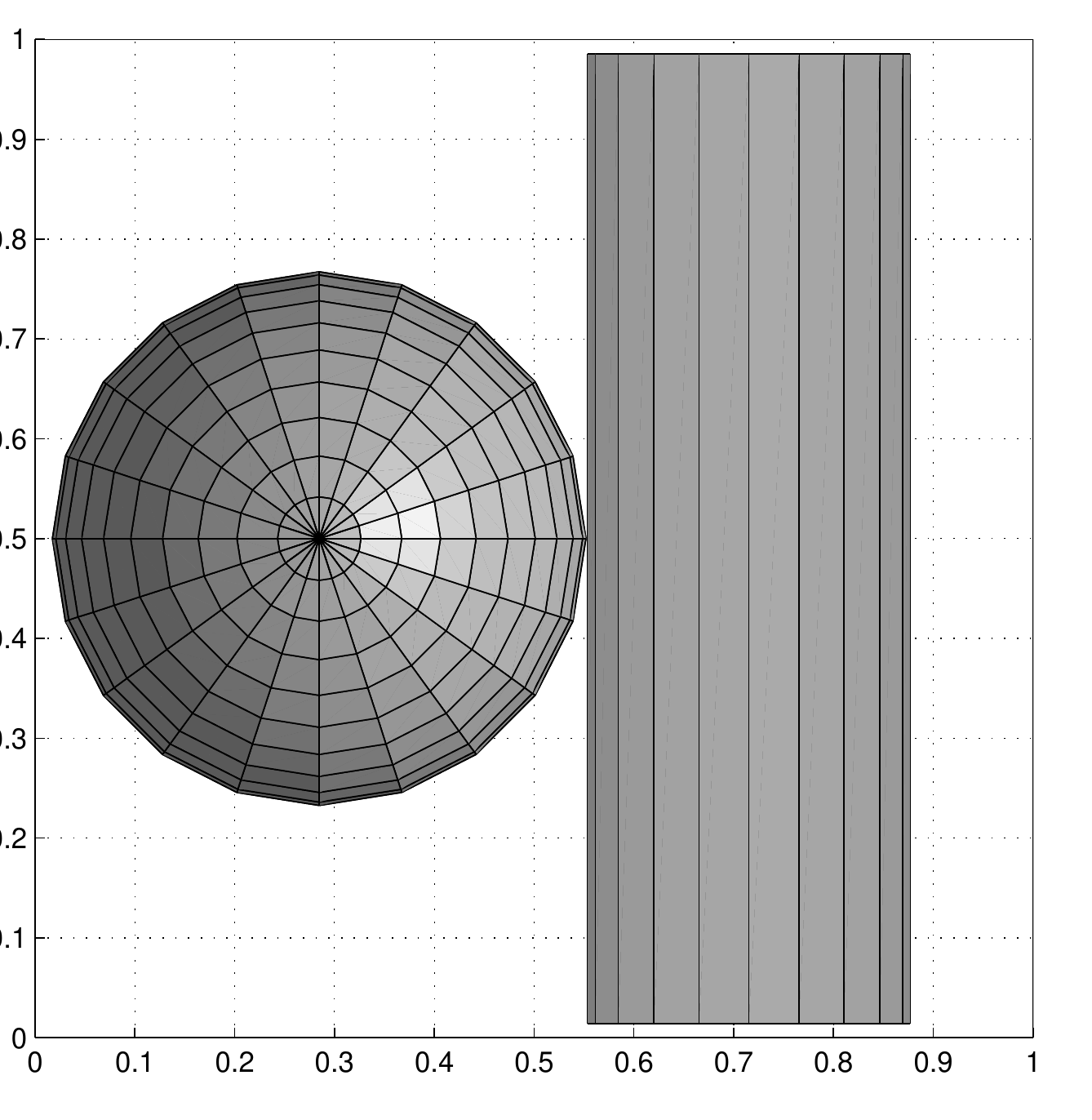}  
 }
\caption{A sphere and a cylinder, symmetric intersection, type $sc2$.}
\label{fig:sp-cyl1}
\end{figure}

\begin{figure}[H] 
\centering
\subfigure[\, Intersecting]{    
\includegraphics[trim = 5cm 7cm 4.5cm 7cm, clip, width=0.45\linewidth]{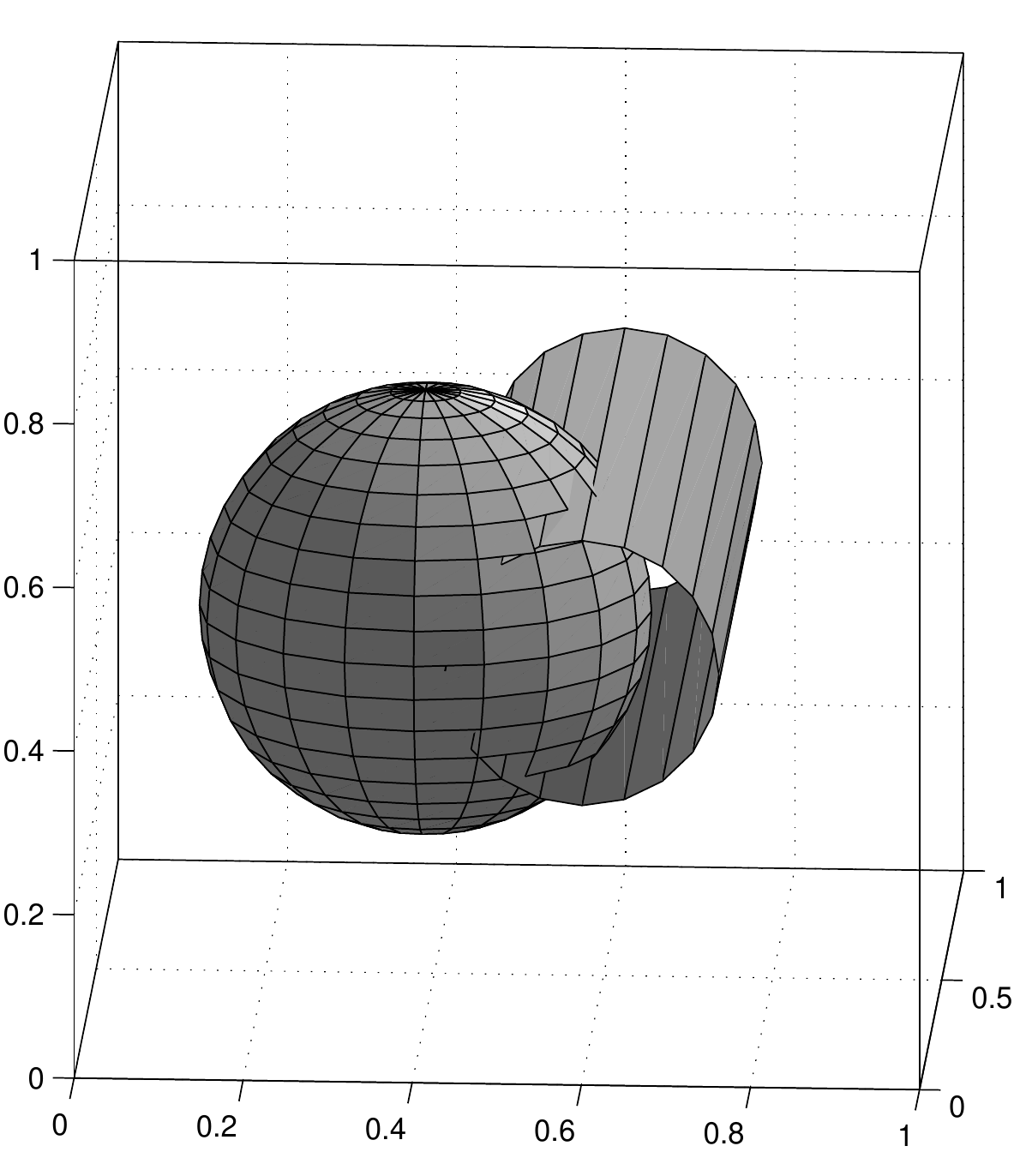}  
}
\subfigure[\, After relaxation]{
    \includegraphics[trim = 5cm 7cm 4.5cm 7cm, clip, width=0.45\linewidth]{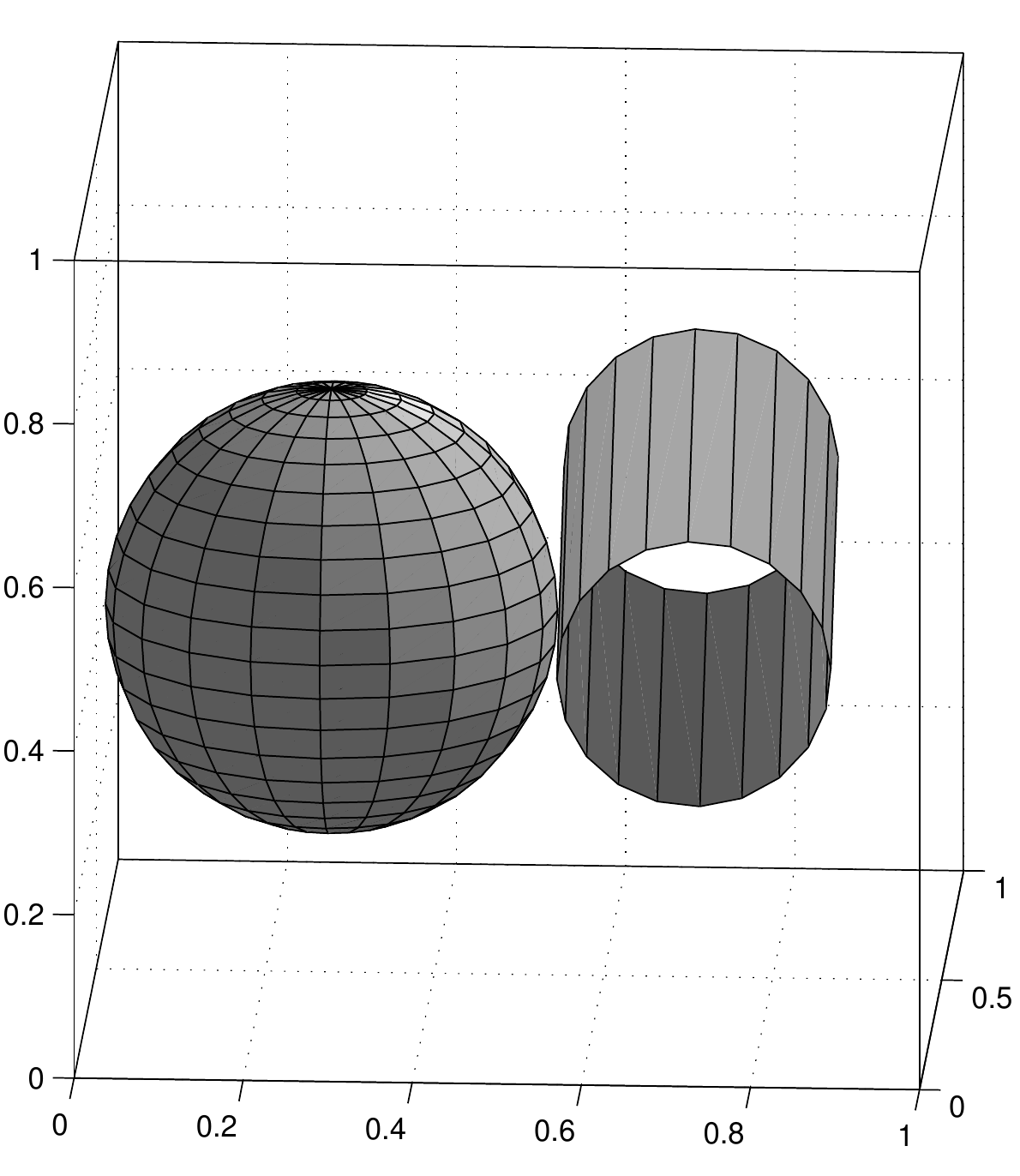}  
 }
\subfigure[\, Intersecting, top view]{    
\includegraphics[trim = 4.5cm 7cm 3.5cm 7cm, clip, width=0.45\linewidth]{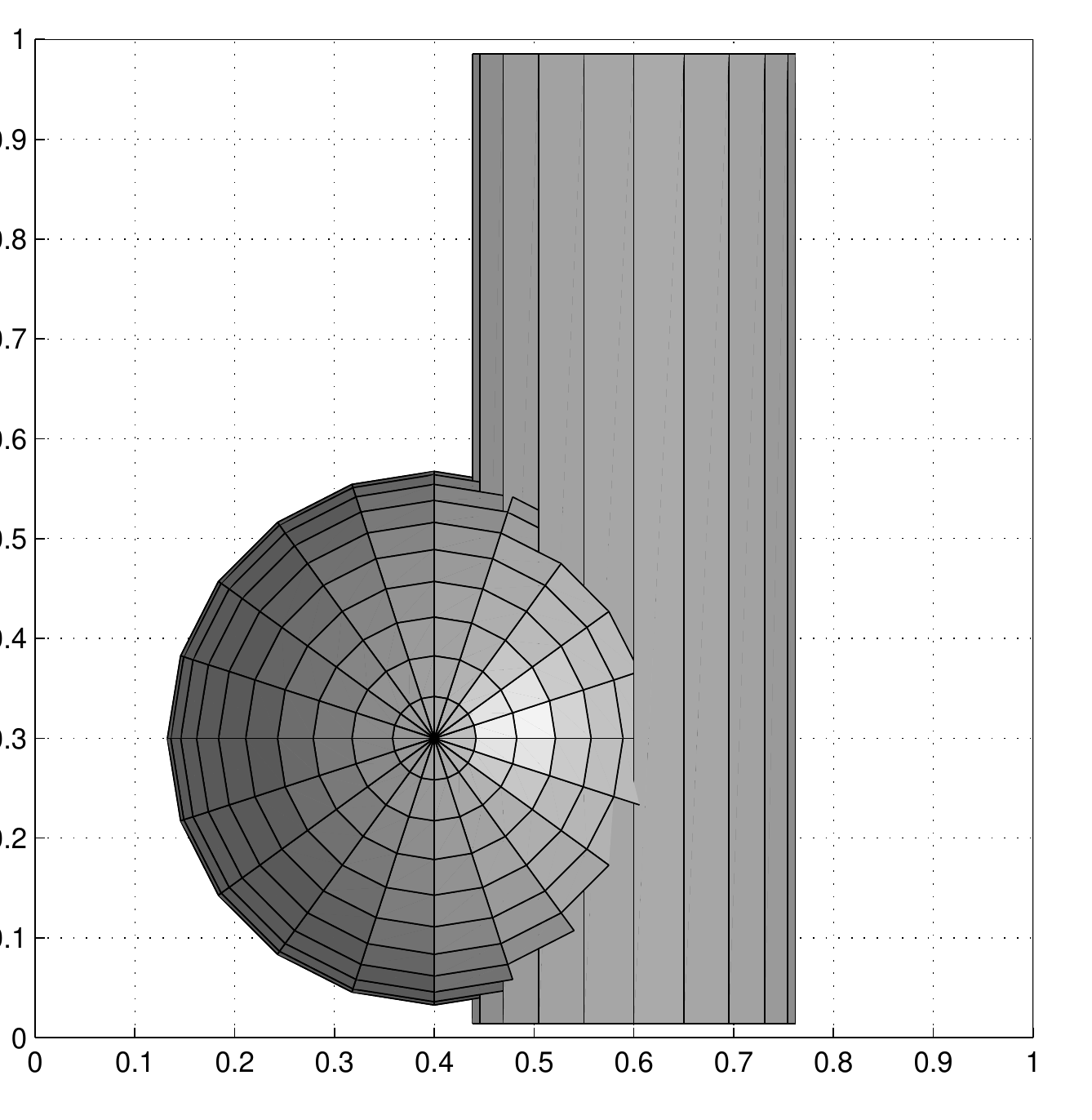}  
}
\subfigure[\, After relaxation, top view]{
    \includegraphics[trim = 4.5cm 7cm 3.5cm 7cm, clip, width=0.45\linewidth]{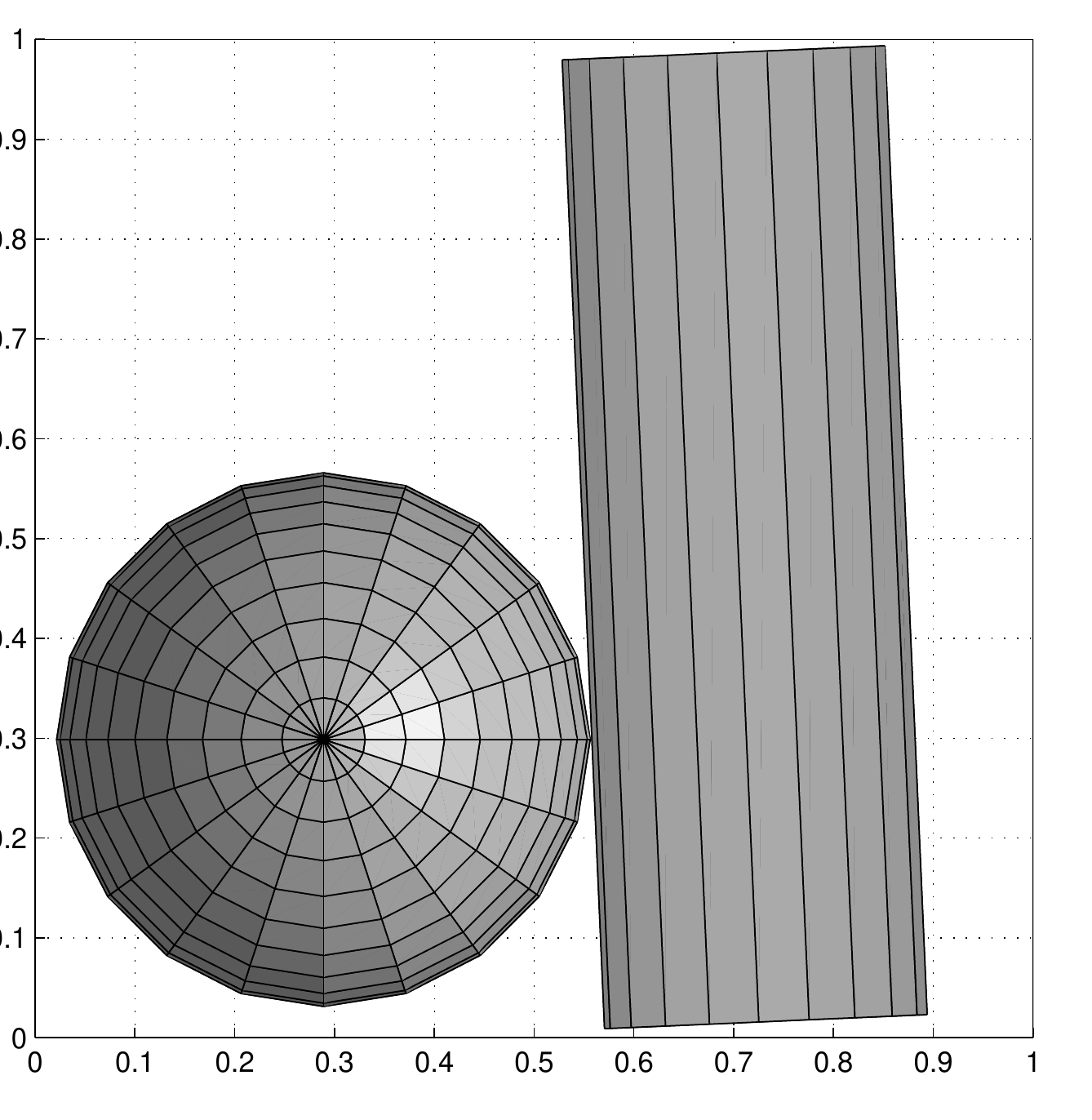}  
 }
\caption{A sphere and a cylinder, not symmetric intersection, type $sc2$  -- the cylinder is turning (compare with fig. \ref{fig:sp-cyl1}).}
\label{fig:sp-cyl2}
\end{figure}

\begin{figure}[H] 
\centering
\subfigure[\, Intersecting]{    
\includegraphics[trim = 5cm 7cm 4.5cm 7cm, clip, width=0.45\linewidth]{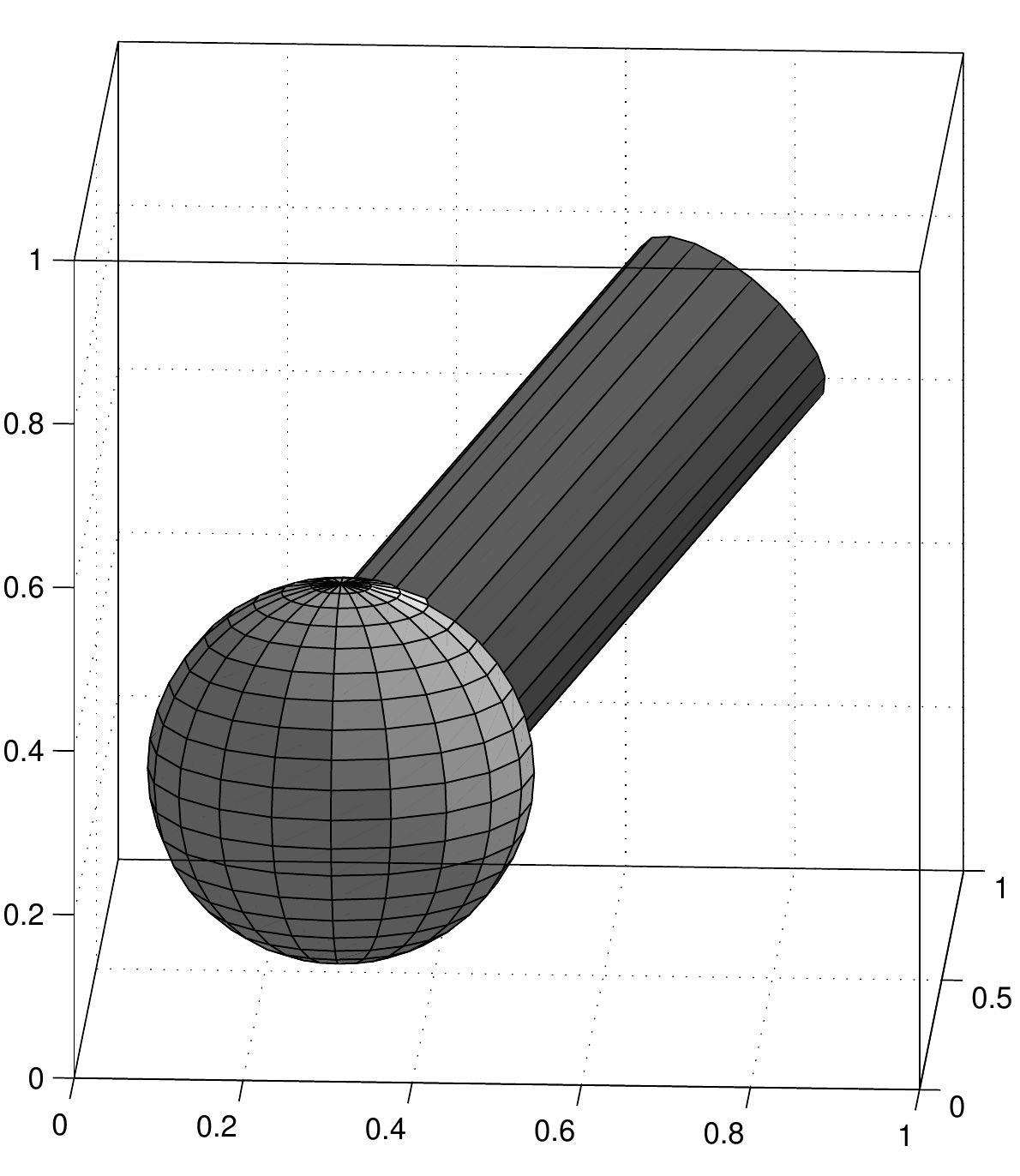}  
}
\subfigure[\, After relaxation]{
    \includegraphics[trim = 5cm 7cm 4.5cm 7cm, clip, width=0.45\linewidth]{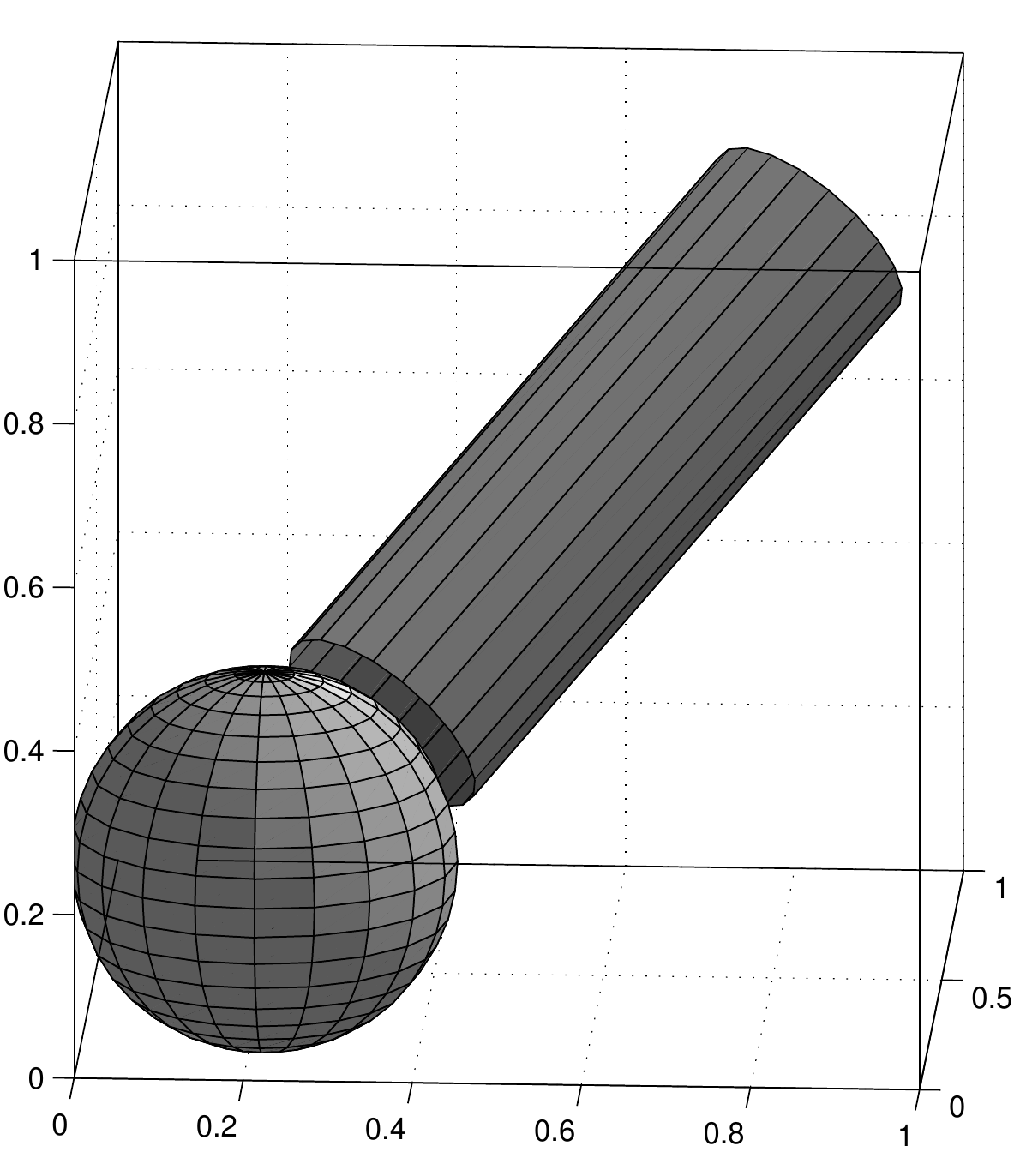}  
 }
\caption{A sphere and a cylinder, axially symmetric intersection with the base, type $sc2$
then $sc3$.}
\label{fig:sp-cyl3}
\end{figure}

\begin{figure}[H] 
\centering
\subfigure[\, Intersecting]{    
\includegraphics[trim = 5cm 7cm 4.5cm 7cm, clip, width=0.45\linewidth]{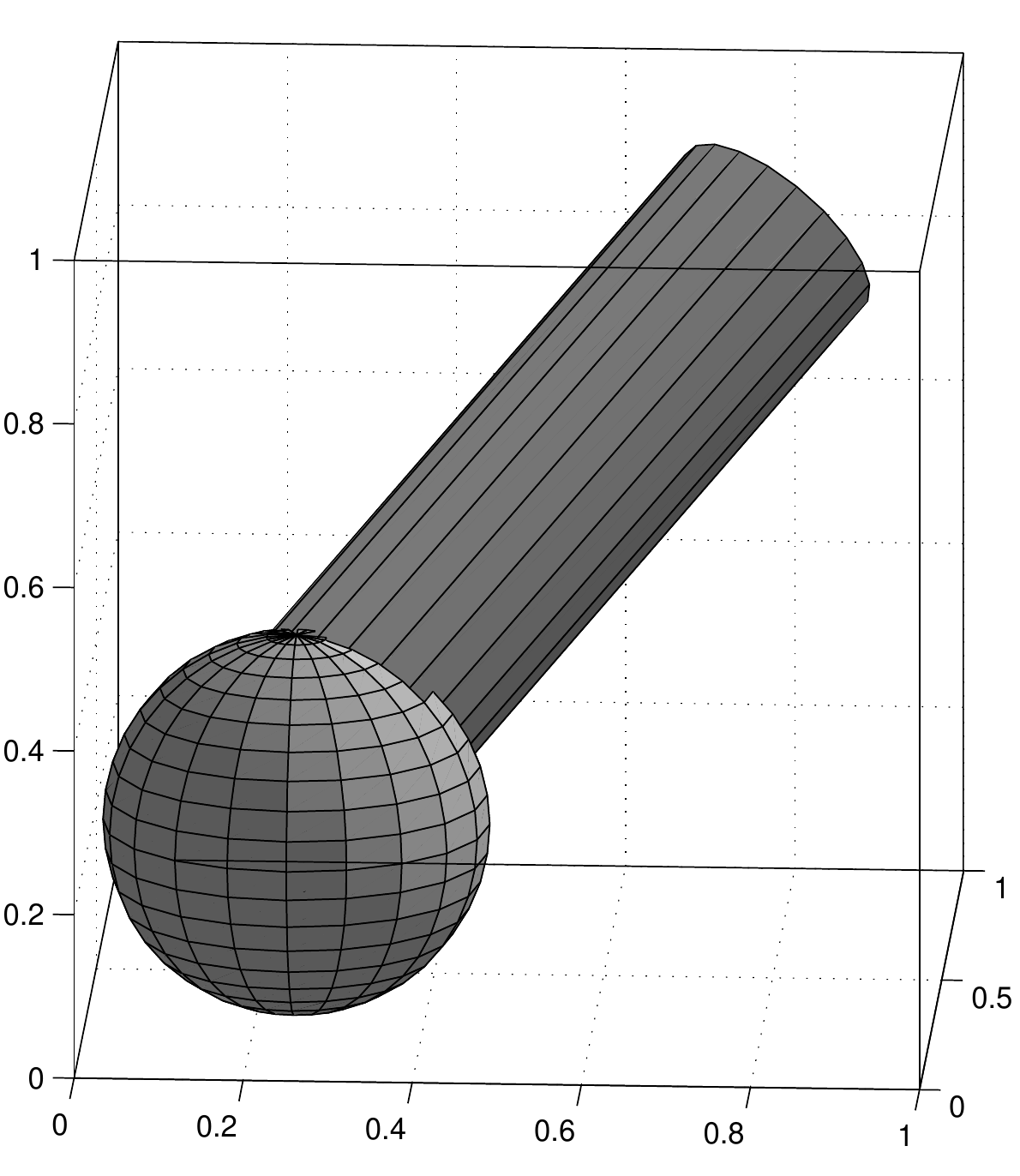}  
}
\subfigure[\, After relaxation]{
    \includegraphics[trim = 5cm 7cm 4.5cm 7cm, clip, width=0.45\linewidth]{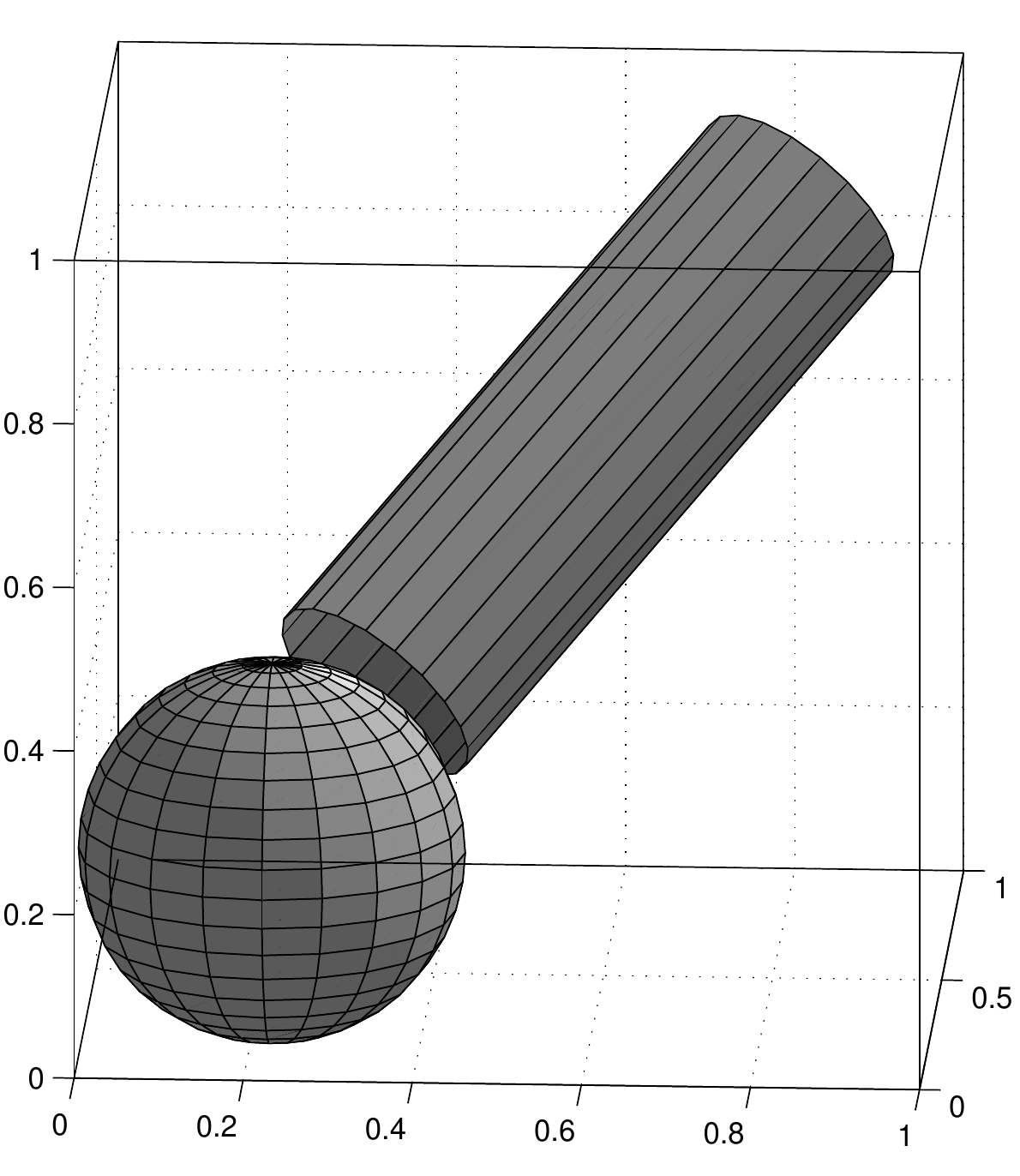}  
 }
\caption{A sphere and a cylinder, intersection with the base, type $sc3$.}
\label{fig:sp-cyl4}
\end{figure}

\begin{figure}[H] 
\centering
\subfigure[\, Intersecting]{    
\includegraphics[trim = 5cm 7cm 4.5cm 7cm, clip, width=0.45\linewidth]{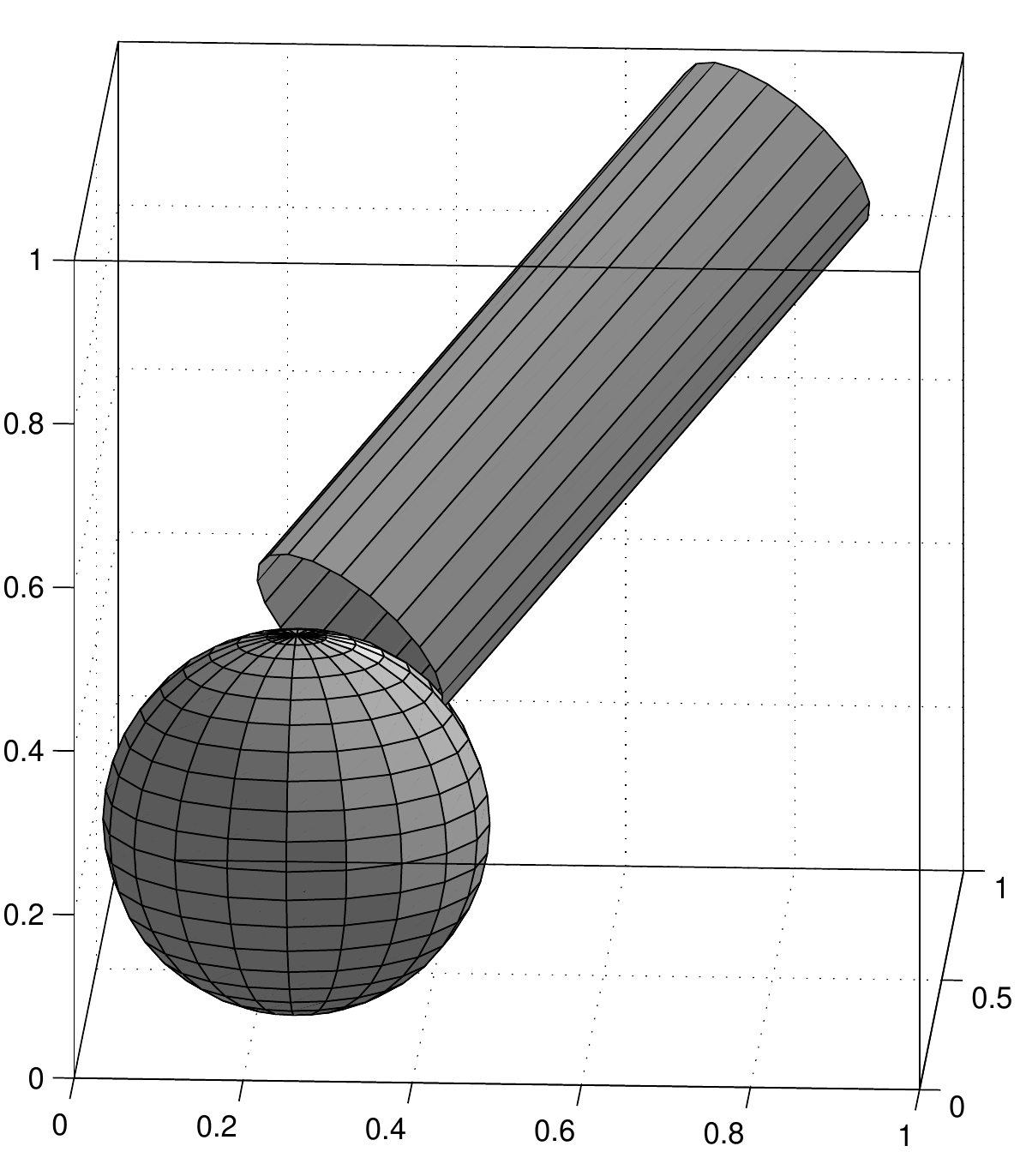}  
}
\subfigure[\, After relaxation]{
    \includegraphics[trim = 5cm 7cm 4.5cm 7cm, clip, width=0.45\linewidth]{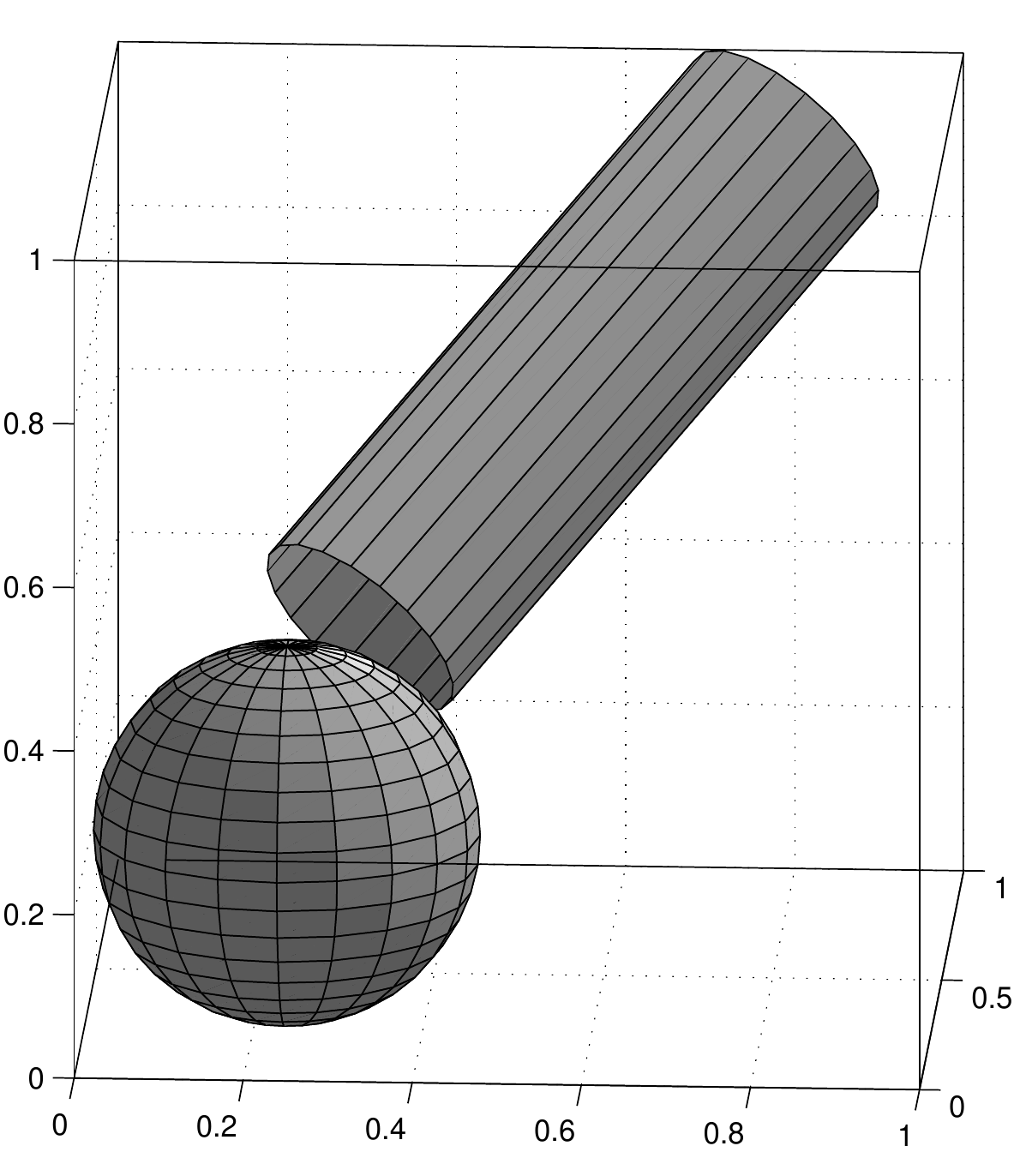}  
 }
\caption{A sphere and a cylinder, intersection with the base boundary, type $sc4$.}
\label{fig:sp-cyl5}
\end{figure}

\begin{figure}[H] 
\centering
\subfigure[\, Intersecting]{    
\includegraphics[trim = 5cm 7cm 4.5cm 7cm, clip, width=0.45\linewidth]{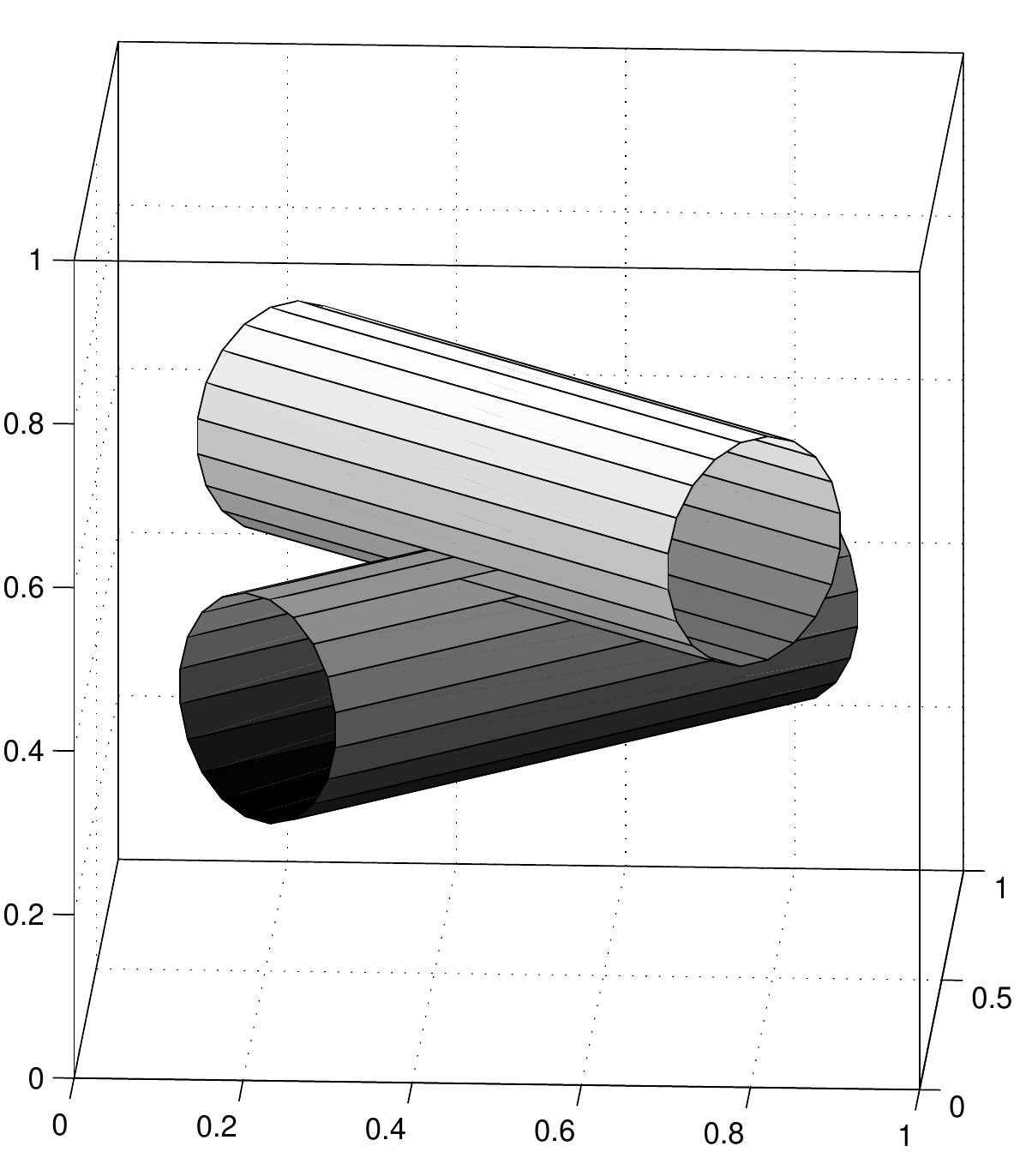}  
}
\subfigure[\, After relaxation]{
    \includegraphics[trim = 5cm 7cm 4.5cm 7cm, clip, width=0.45\linewidth]{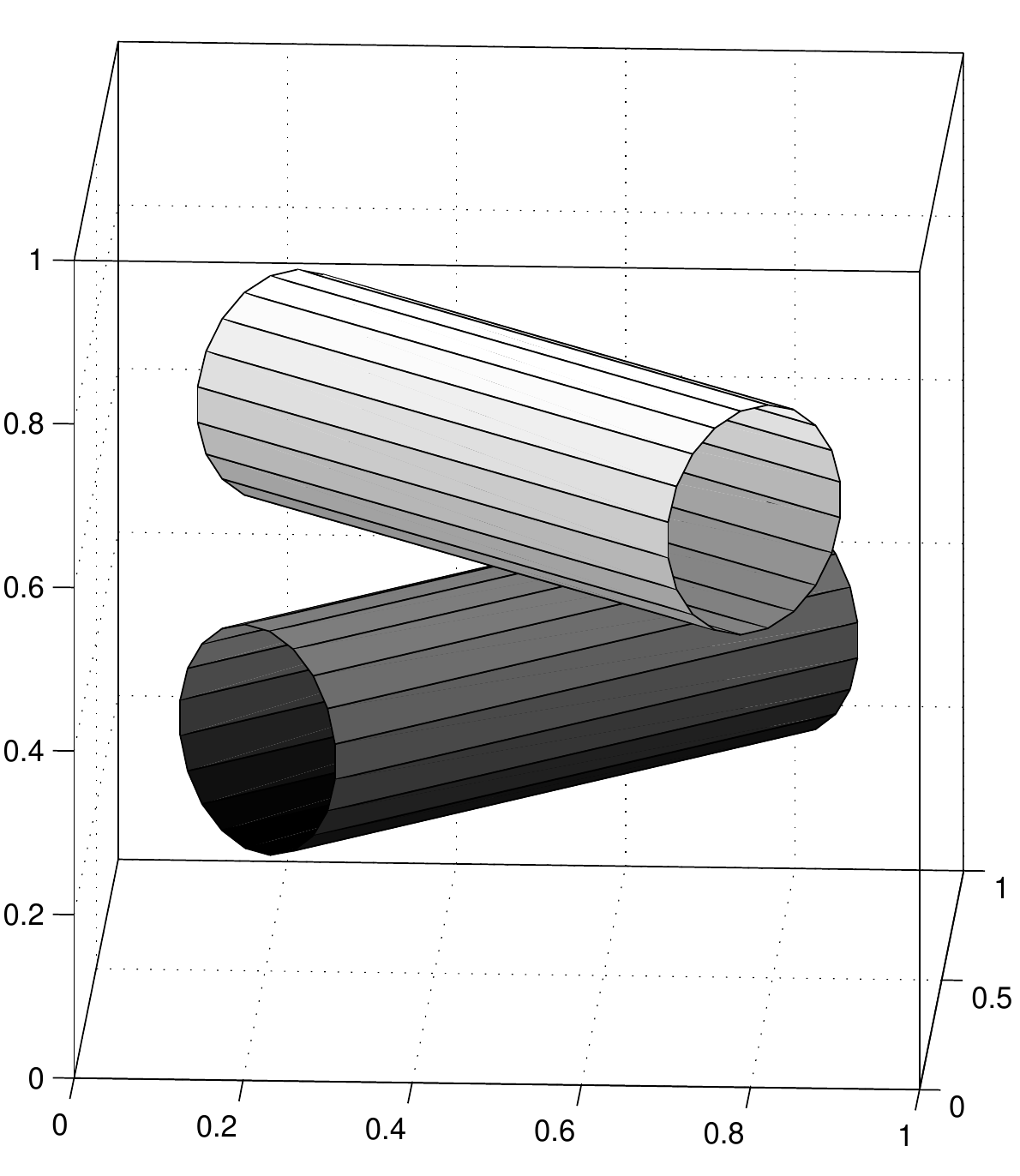}  
 }
\caption{Two cylinders, symmetric intersection, type $cc1$.}
\label{fig:cyl-cyl1}
\end{figure}

\begin{figure}[H] 
\centering
\subfigure[\, Intersecting]{    
\includegraphics[trim = 5cm 7cm 4.5cm 7cm, clip, width=0.45\linewidth]{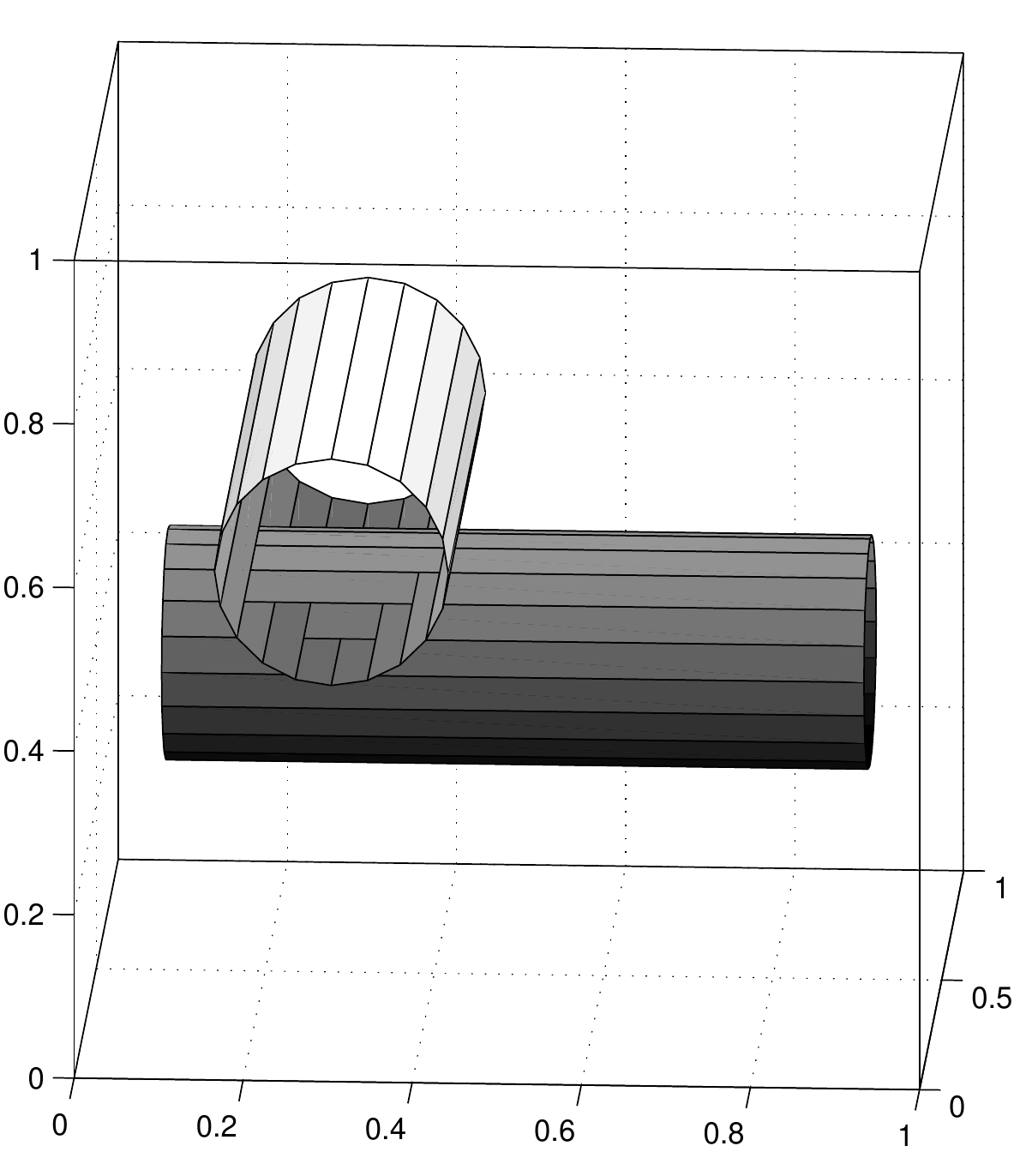}  
}
\subfigure[\, After relaxation]{
    \includegraphics[trim = 5cm 7cm 4.5cm 7cm, clip, width=0.45\linewidth]{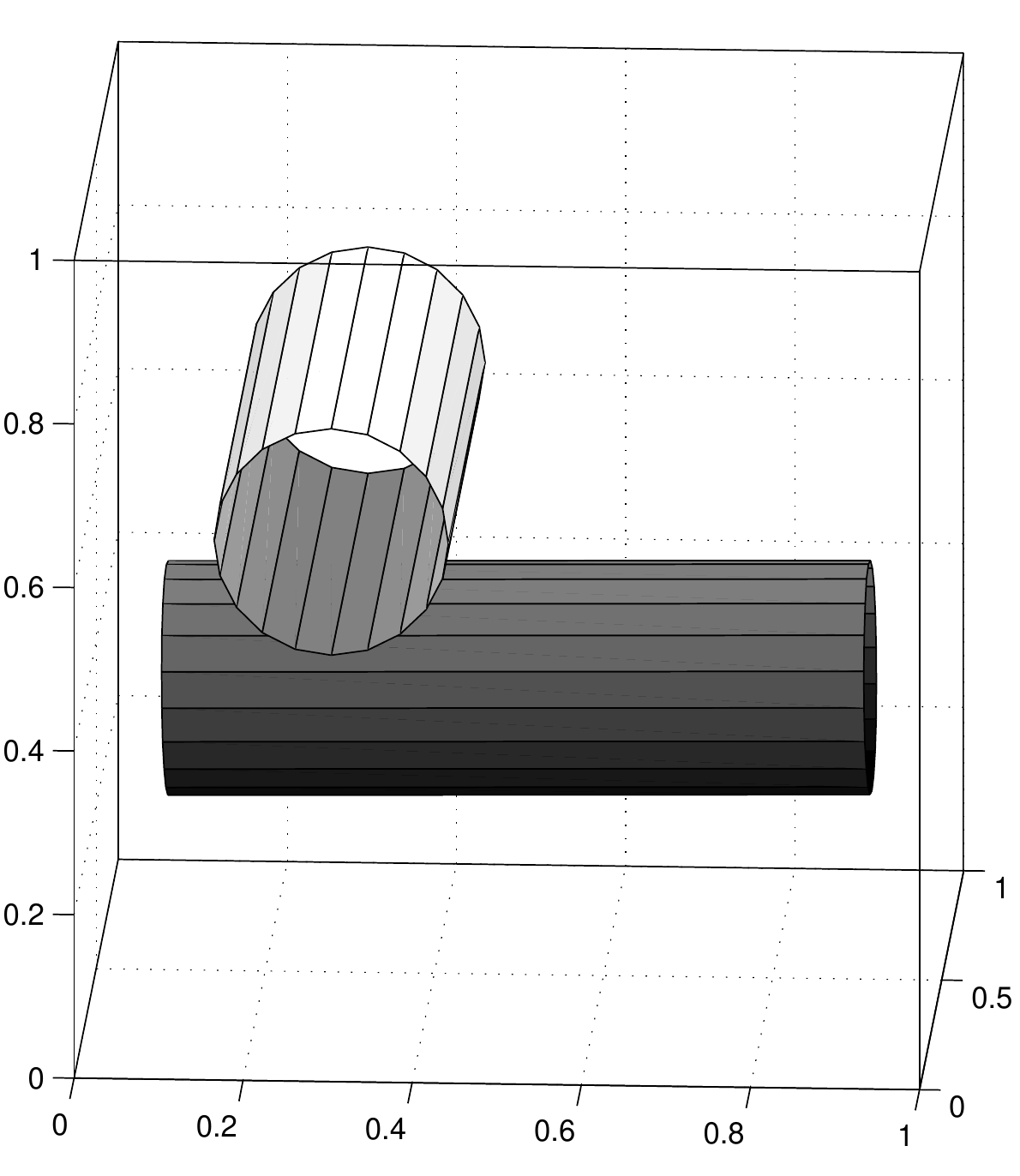}  
 }
\subfigure[\, Intersecting, front view]{    
\includegraphics[trim = 4.5cm 7cm 3.5cm 7cm, clip, width=0.45\linewidth]{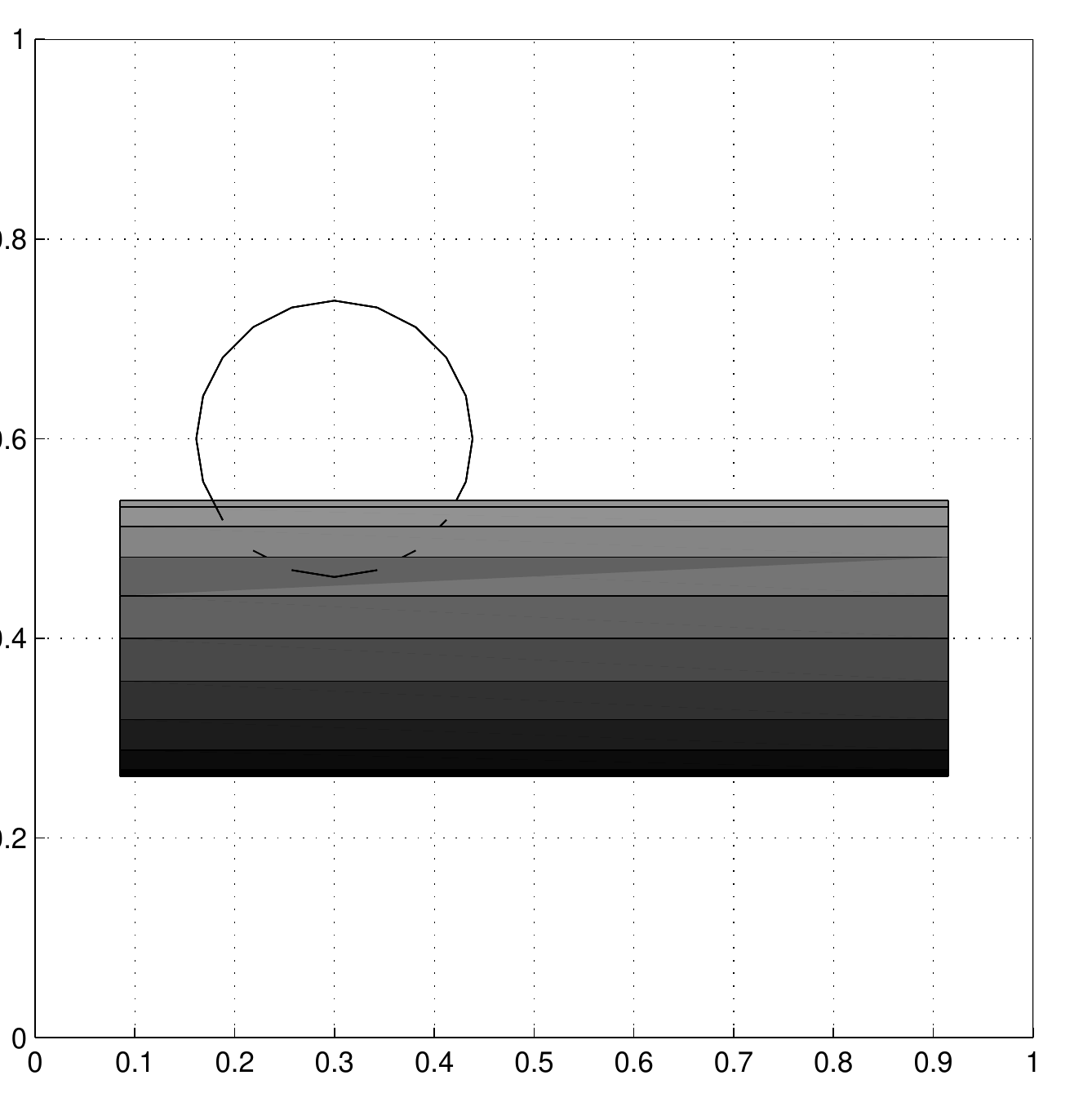}  
}
\subfigure[\, After relaxation, front view]{
    \includegraphics[trim = 4.5cm 7cm 3.5cm 7cm, clip, width=0.45\linewidth]{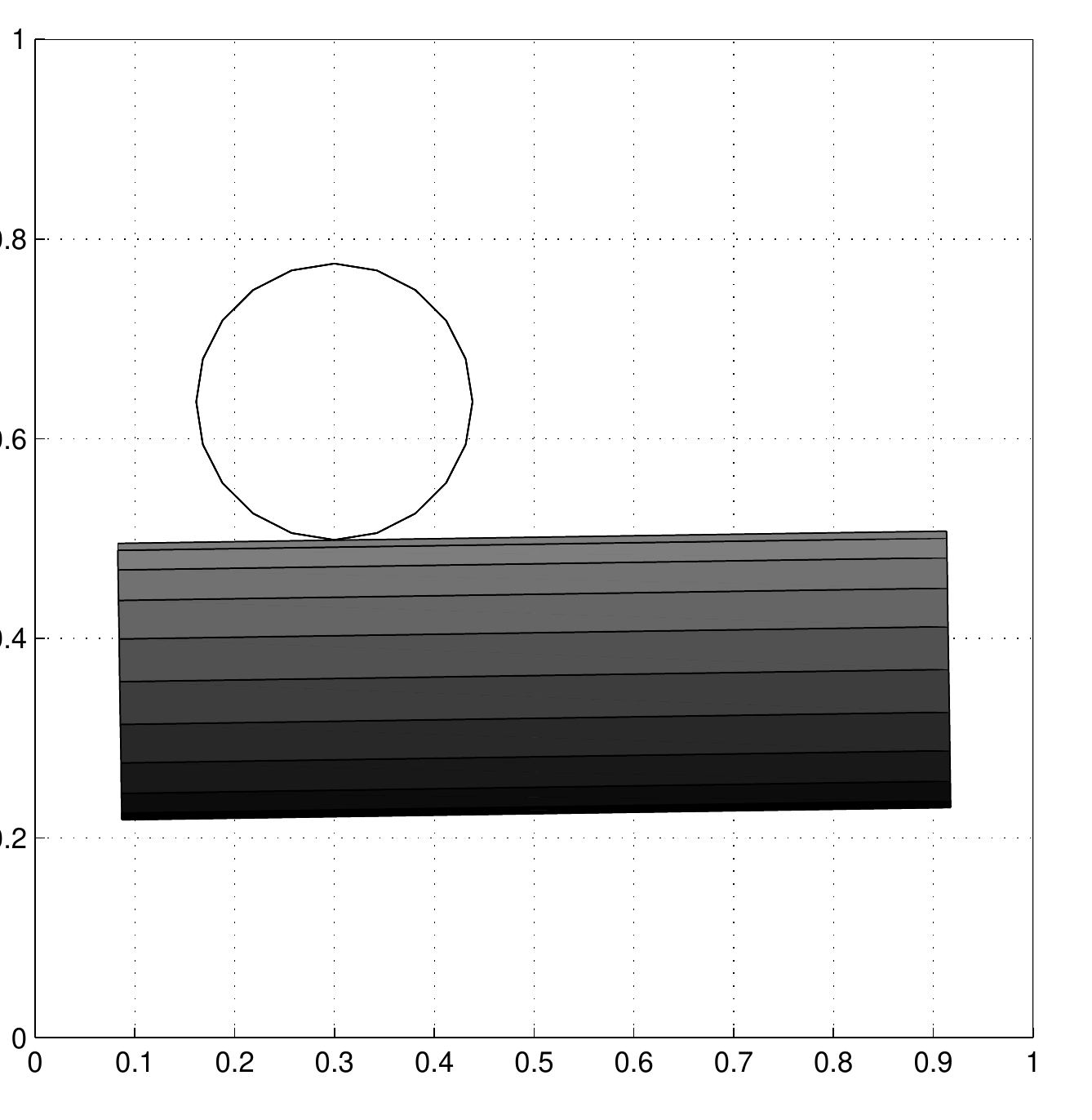}  
 }
\caption{Two cylinders, not symmetric intersection, type $cc1$ -- one cylinder is turning.}
\label{fig:cyl-cyl2}
\end{figure}

\begin{figure}[H] 
\centering
\subfigure[\, Intersecting]{    
\includegraphics[trim = 5cm 7cm 4.5cm 7cm, clip, width=0.45\linewidth]{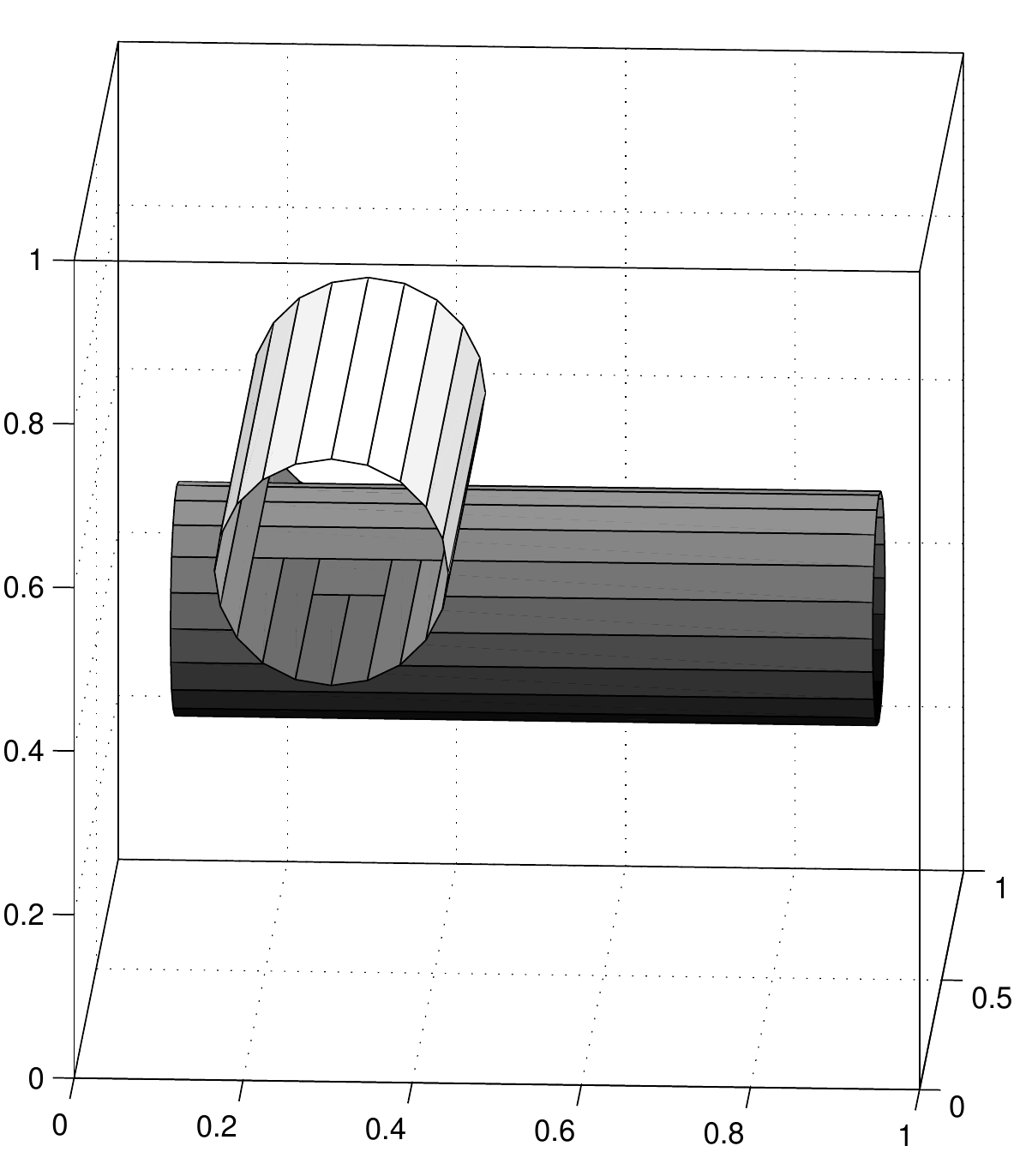}  
}
\subfigure[\, After relaxation]{
    \includegraphics[trim = 5cm 7cm 4.5cm 7cm, clip, width=0.45\linewidth]{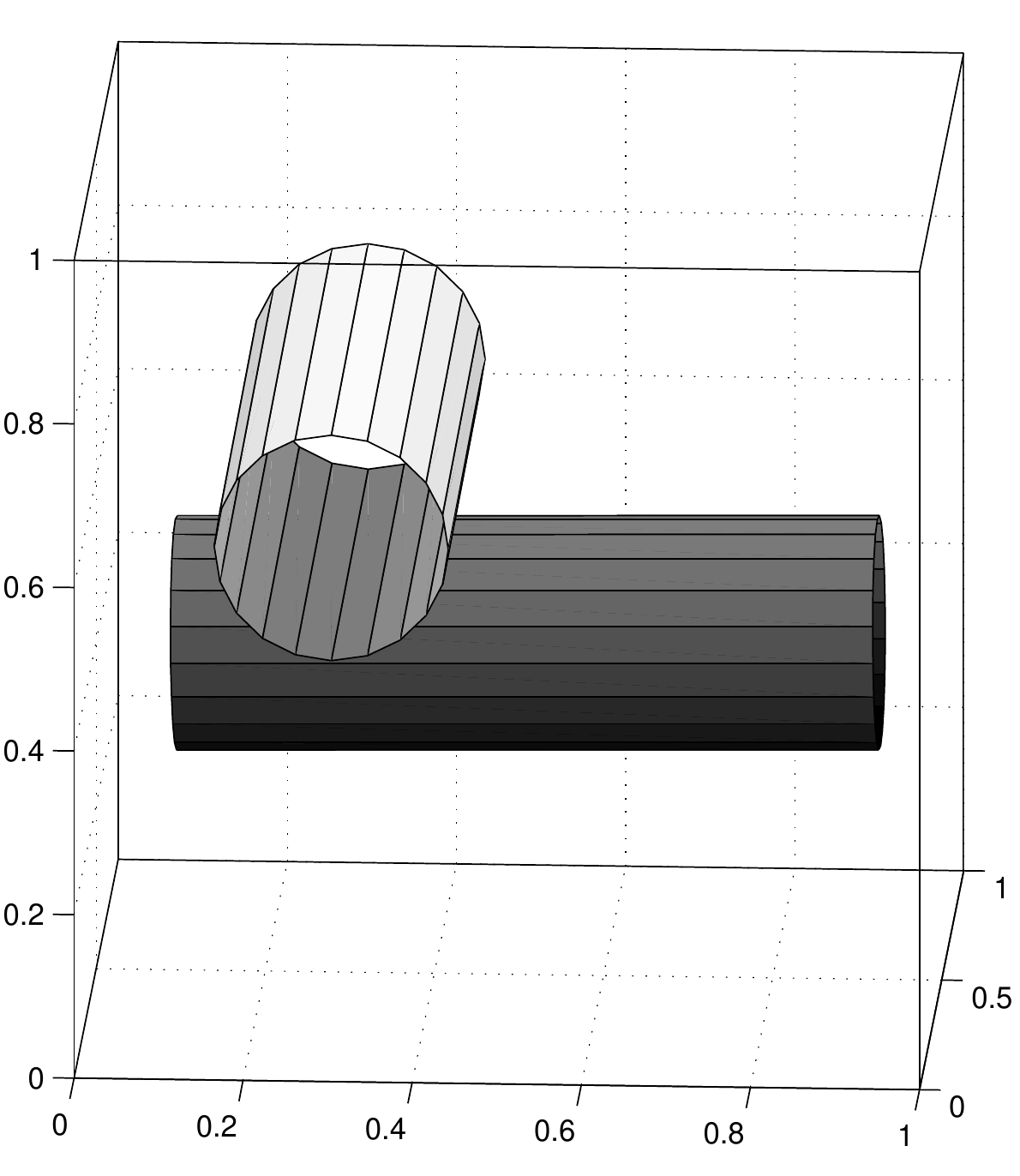}  
 }
\subfigure[\, Intersecting, front view]{    
\includegraphics[trim = 4.5cm 7cm 3.5cm 7cm, clip, width=0.45\linewidth]{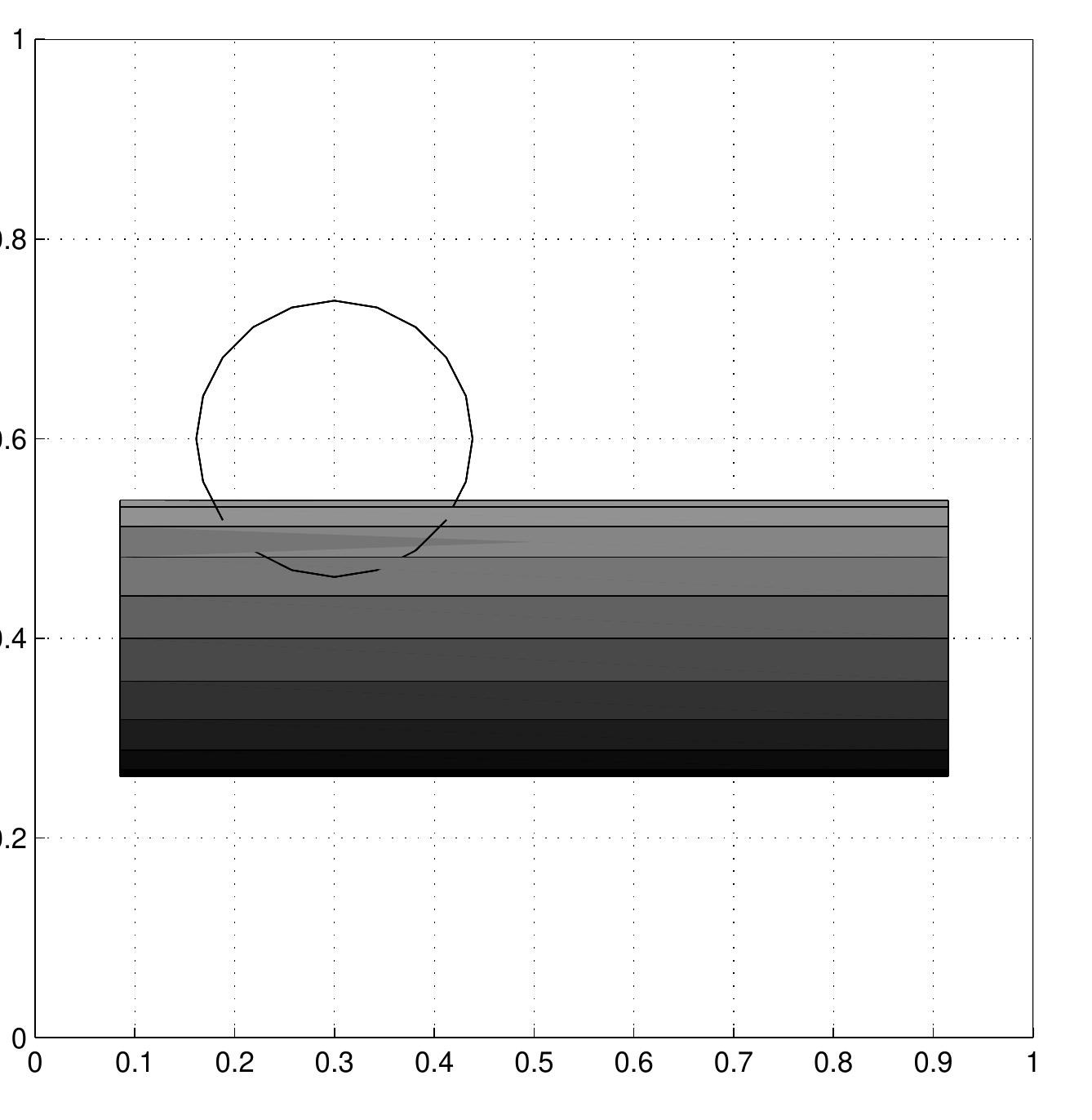}  
}
\subfigure[\, After relaxation, front view]{
    \includegraphics[trim = 4.5cm 7cm 3.5cm 7cm, clip, width=0.45\linewidth]{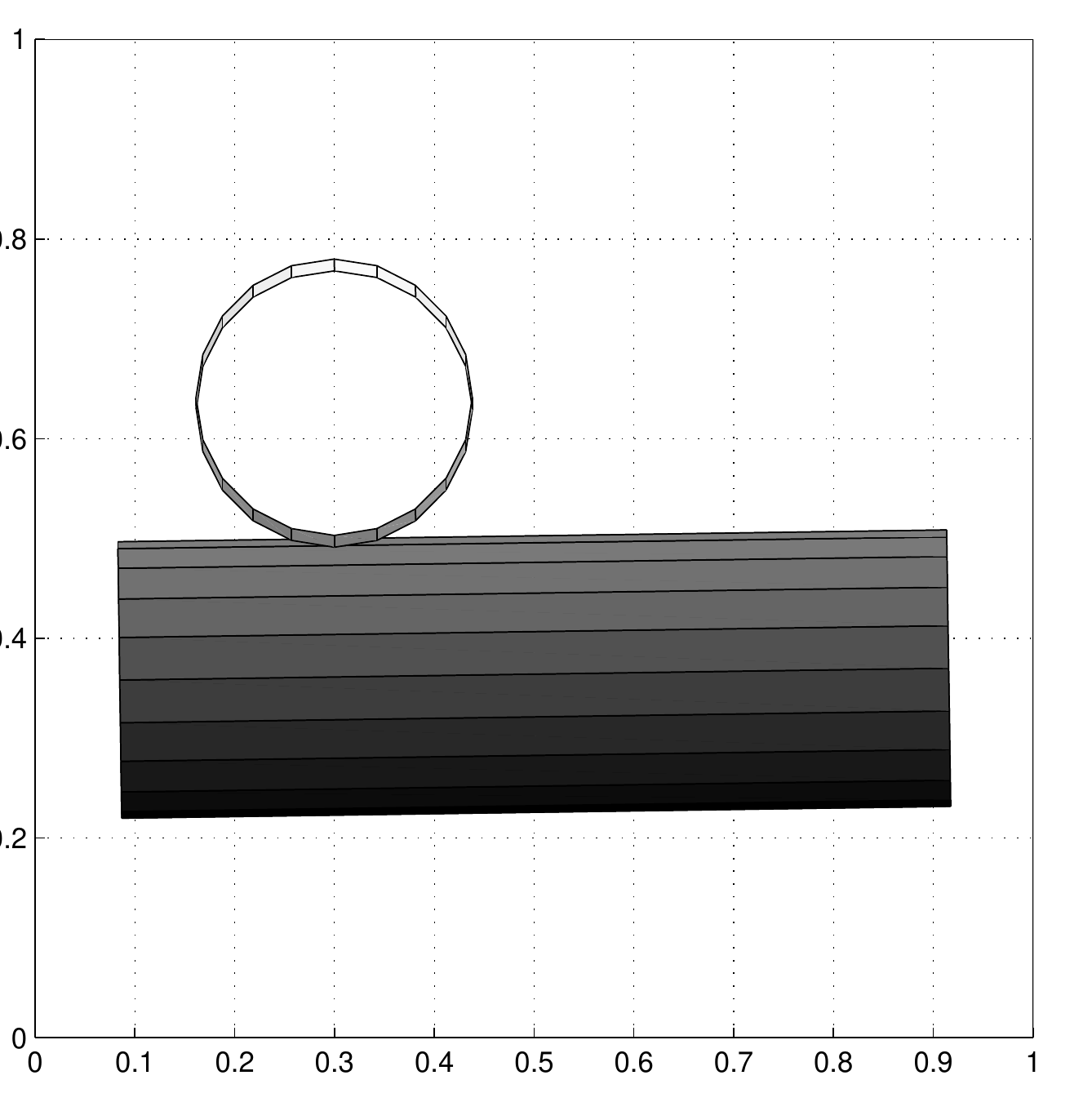}  
 }
\caption{Two cylinders, not symmetric intersection, type $cc1$ -- both cylinders are turning 
(compare with fig. \ref{fig:cyl-cyl2}).}
\label{fig:cyl-cyl3}
\end{figure}

\begin{figure}[H] 
\centering
\subfigure[\, Intersecting]{    
\includegraphics[trim = 5cm 7cm 4.5cm 7cm, clip, width=0.45\linewidth]{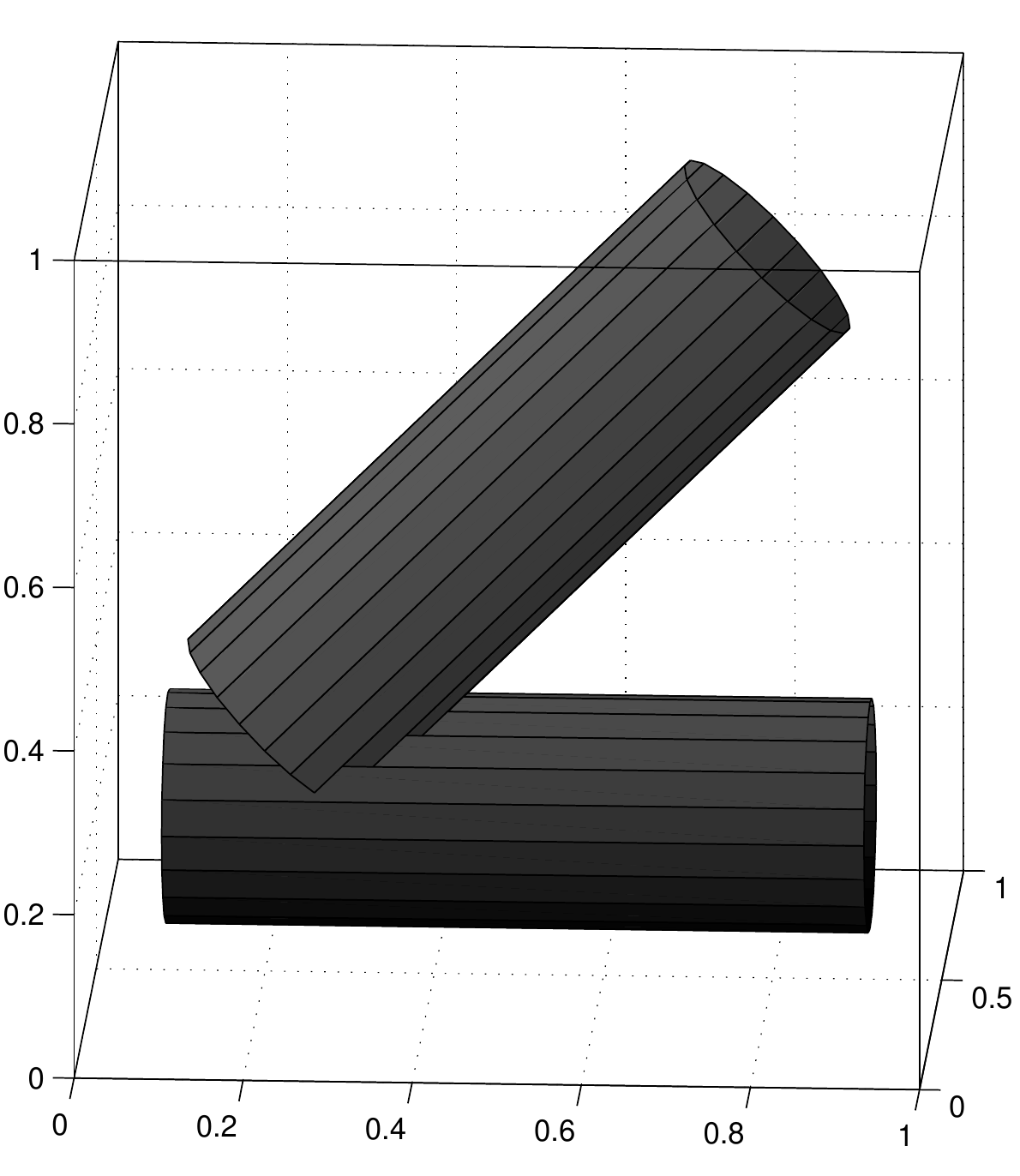}  
}
\subfigure[\, After relaxation]{
    \includegraphics[trim = 5cm 7cm 4.5cm 7cm, clip, width=0.45\linewidth]{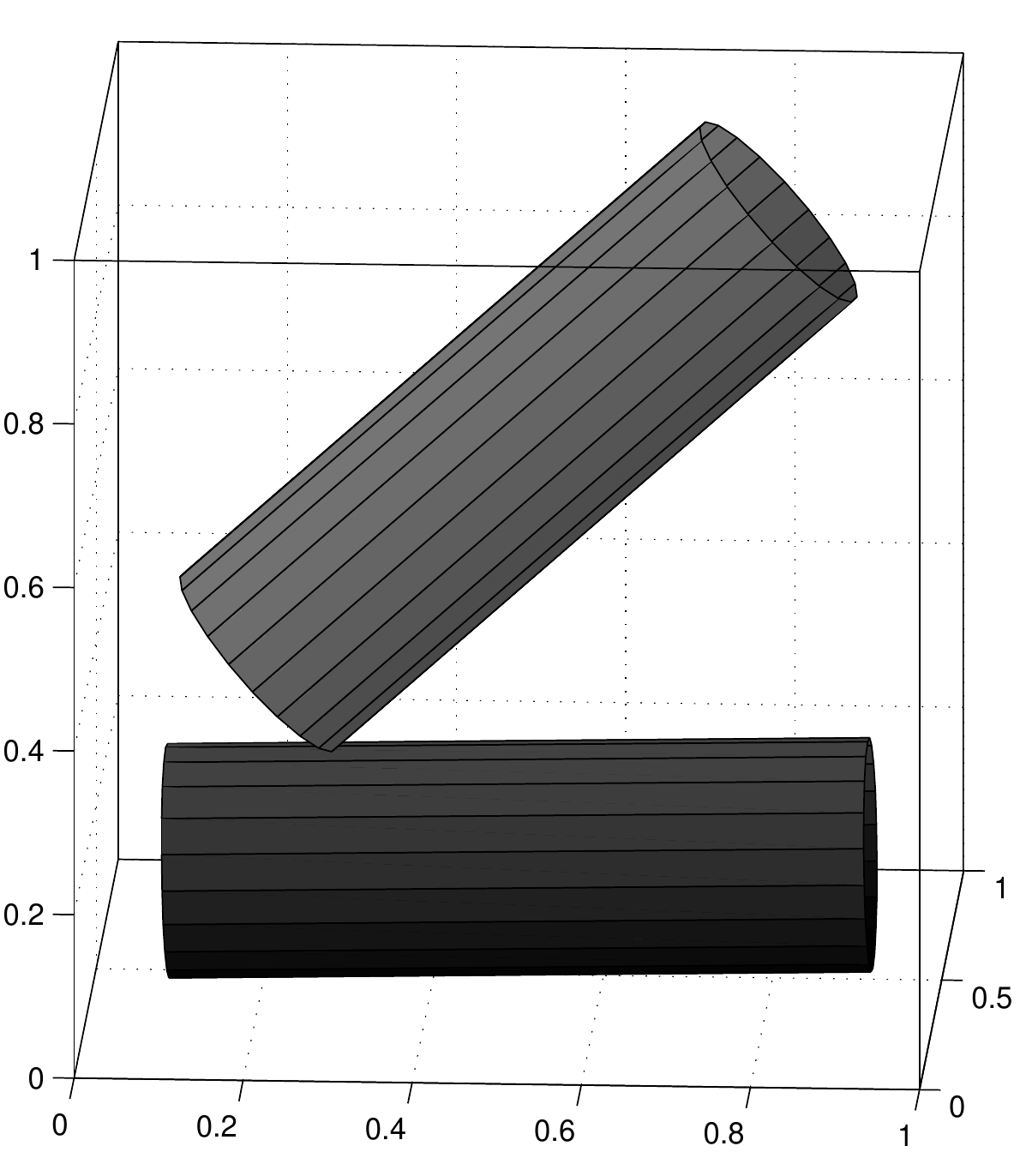}  
 }
\caption{Two cylinders, axes in the same plane, intersection with one base, type $cd1$.}
\label{fig:cyl-cyl4}
\end{figure}

\begin{figure}[H] 
\centering
\subfigure[\, Intersecting]{    
\includegraphics[trim = 5cm 7cm 4.5cm 7cm, clip, width=0.45\linewidth]{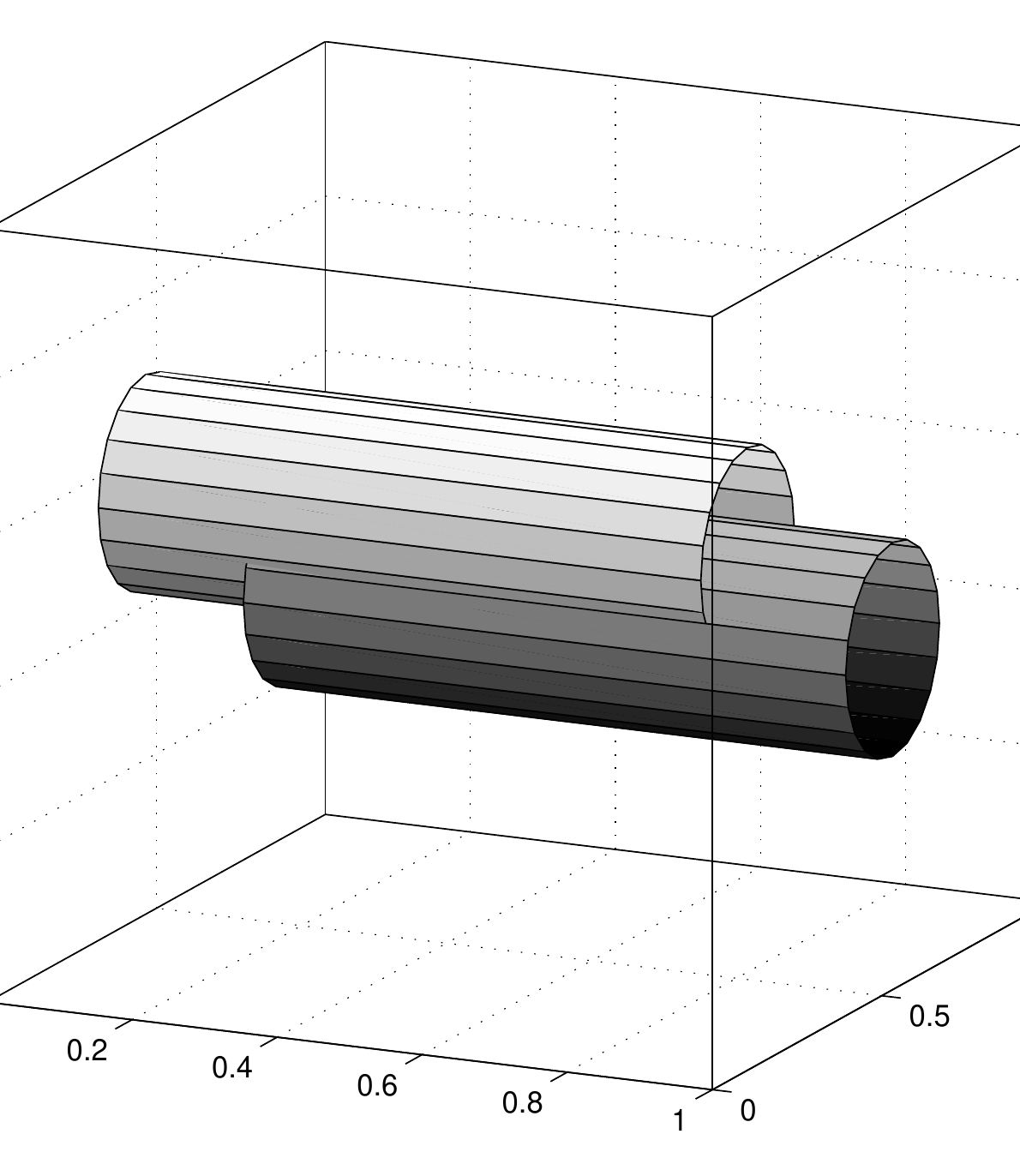}  
}
\subfigure[\, After relaxation]{
    \includegraphics[trim = 5cm 7cm 4.5cm 7cm, clip, width=0.45\linewidth]{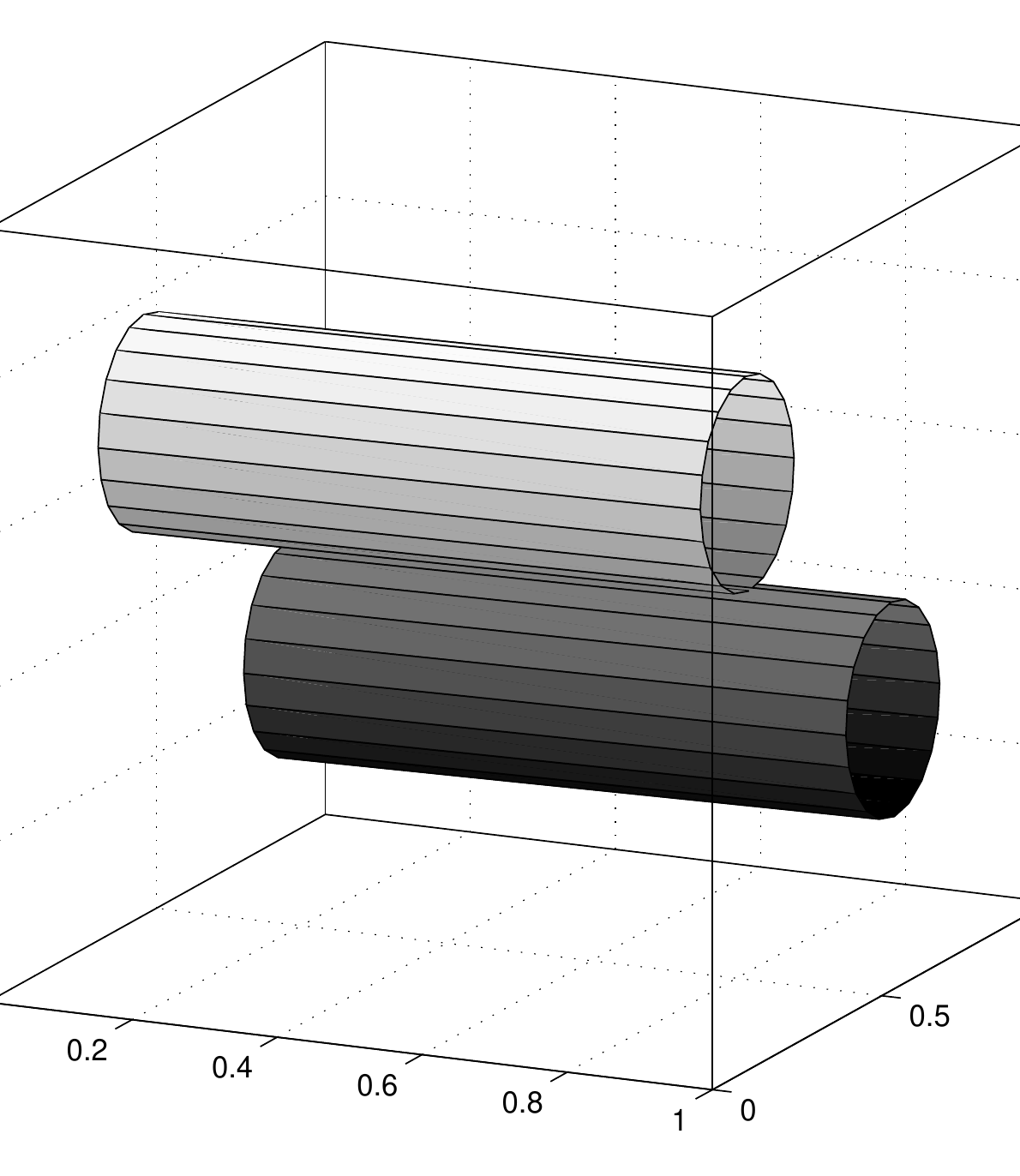}  
 }
\caption{Two cylinders, parallel axes, intersections of type $cc1$ (degenerate case) and $cd2$, then $cd1$.}
\label{fig:cyl-cyl5}
\end{figure}

\begin{figure}[H] 
\centering
\subfigure[\, Intersecting]{    
\includegraphics[trim = 5cm 7cm 4.5cm 7cm, clip, width=0.45\linewidth]{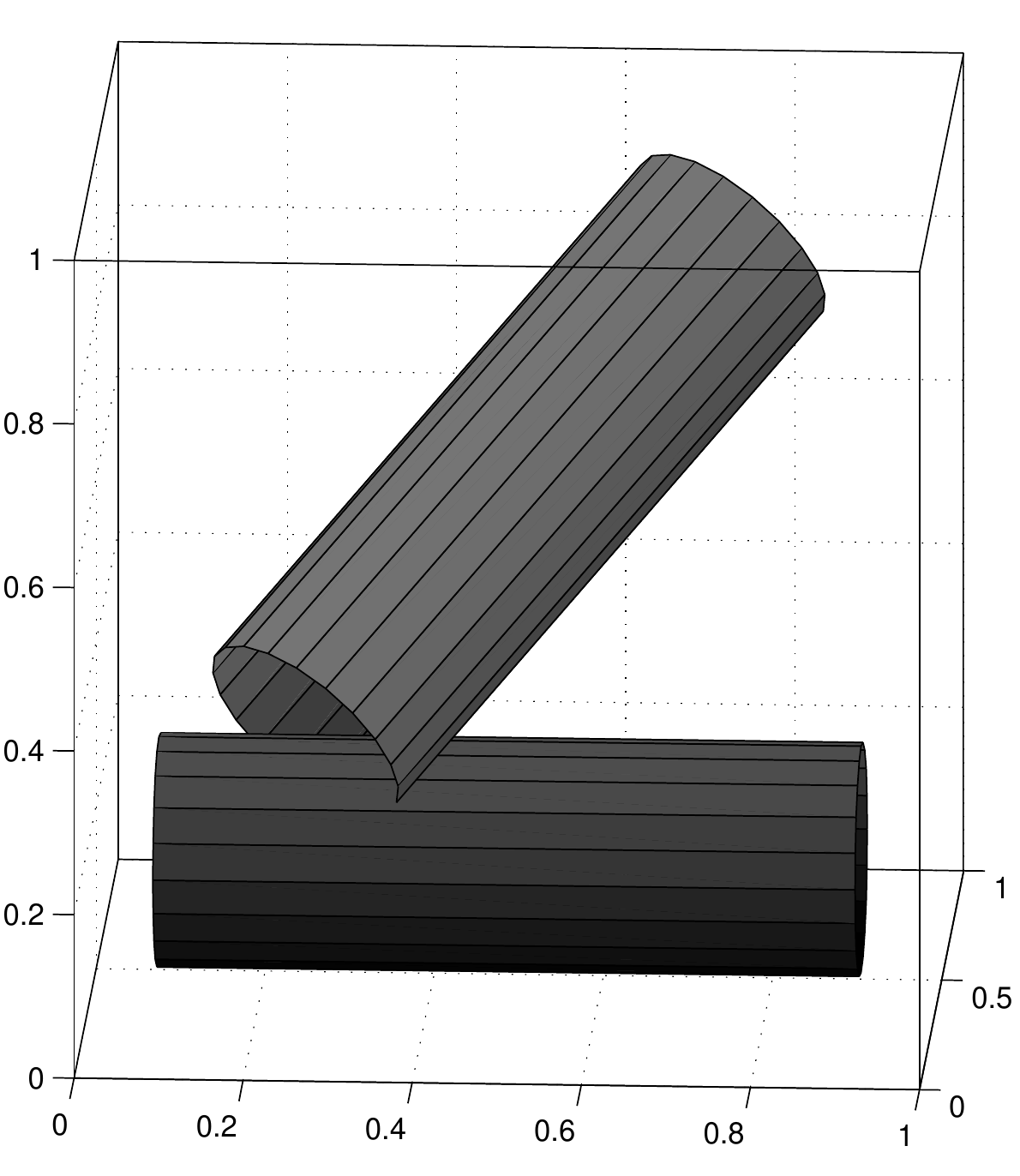}  
}
\subfigure[\, After relaxation]{
    \includegraphics[trim = 5cm 7cm 4.5cm 7cm, clip, width=0.45\linewidth]{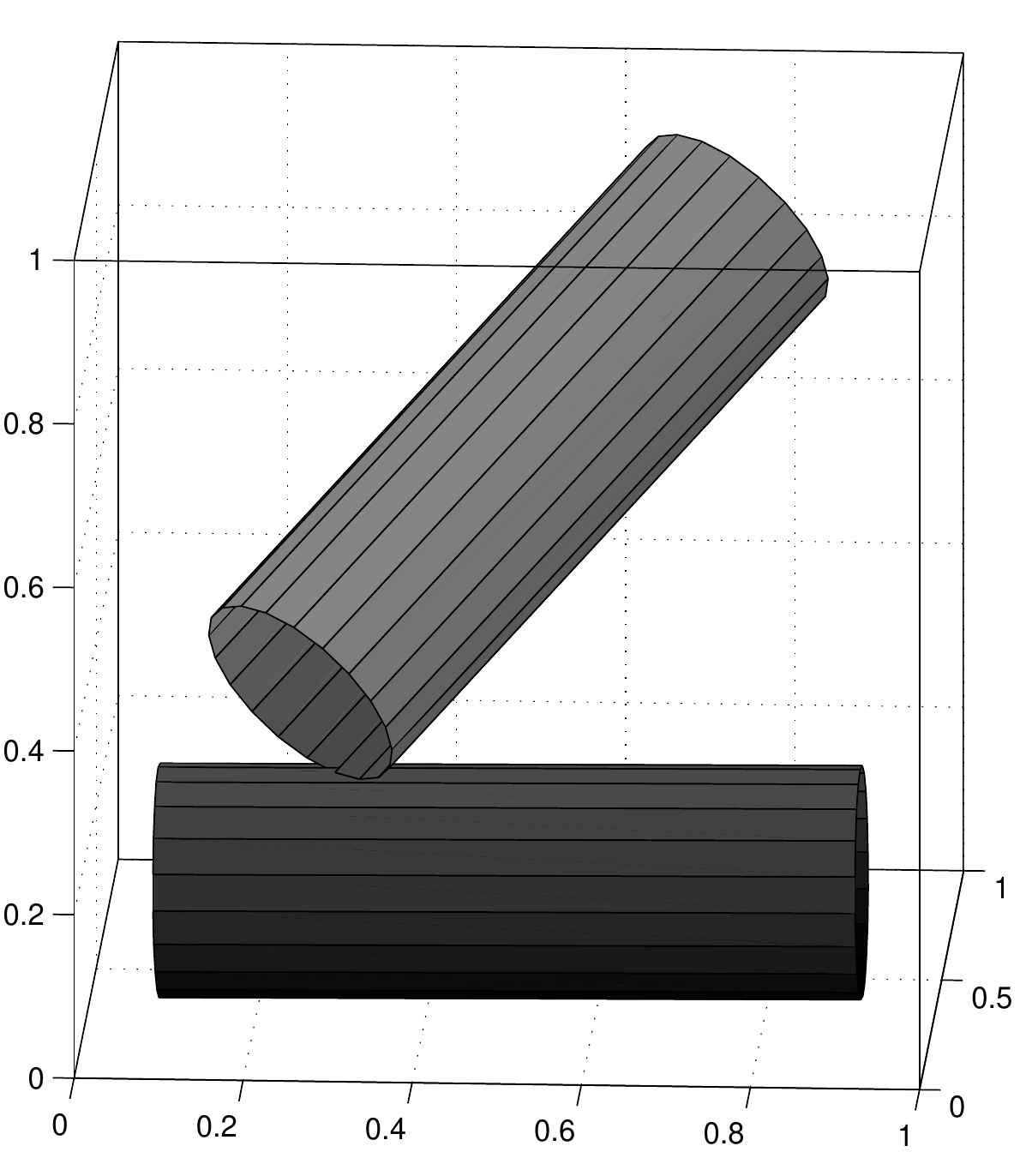}  
 }
\caption{Two cylinders, axes are not coplanar, intersection with one base, type $cd1$.}
\label{fig:cyl-cyl6}
\end{figure}

\begin{figure}[H] 
 \centering
 \subfigure[\, Intersecting]{    
 \includegraphics[trim = 5cm 7cm 4.5cm 7cm, clip, width=0.45\linewidth]{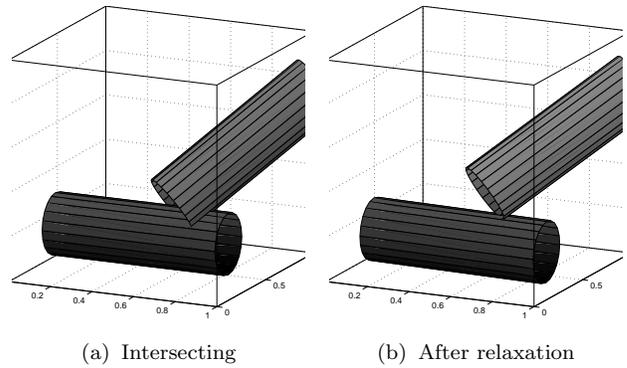}  
 }
 \subfigure[\, After relaxation]{
     \includegraphics[trim = 5cm 7cm 4.5cm 7cm, clip, width=0.45\linewidth]{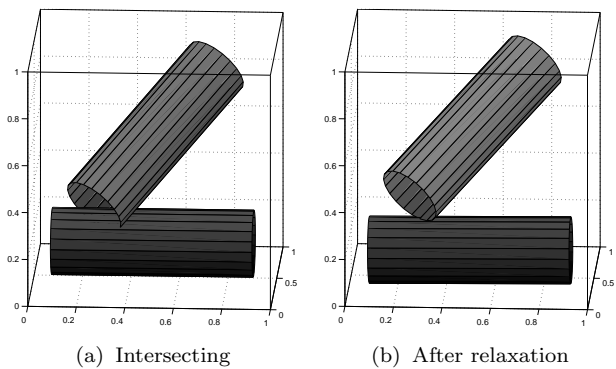}  
  }
 \caption{Two cylinders, axes are not coplanar, intersection with one base, type $cd1$.}
 \label{fig:cyl-cyl7}
 \end{figure}

\begin{figure}[H] 
\centering
\subfigure[\, Intersecting]{    
\includegraphics[trim = 5cm 7cm 4.5cm 7cm, clip, width=0.45\linewidth]{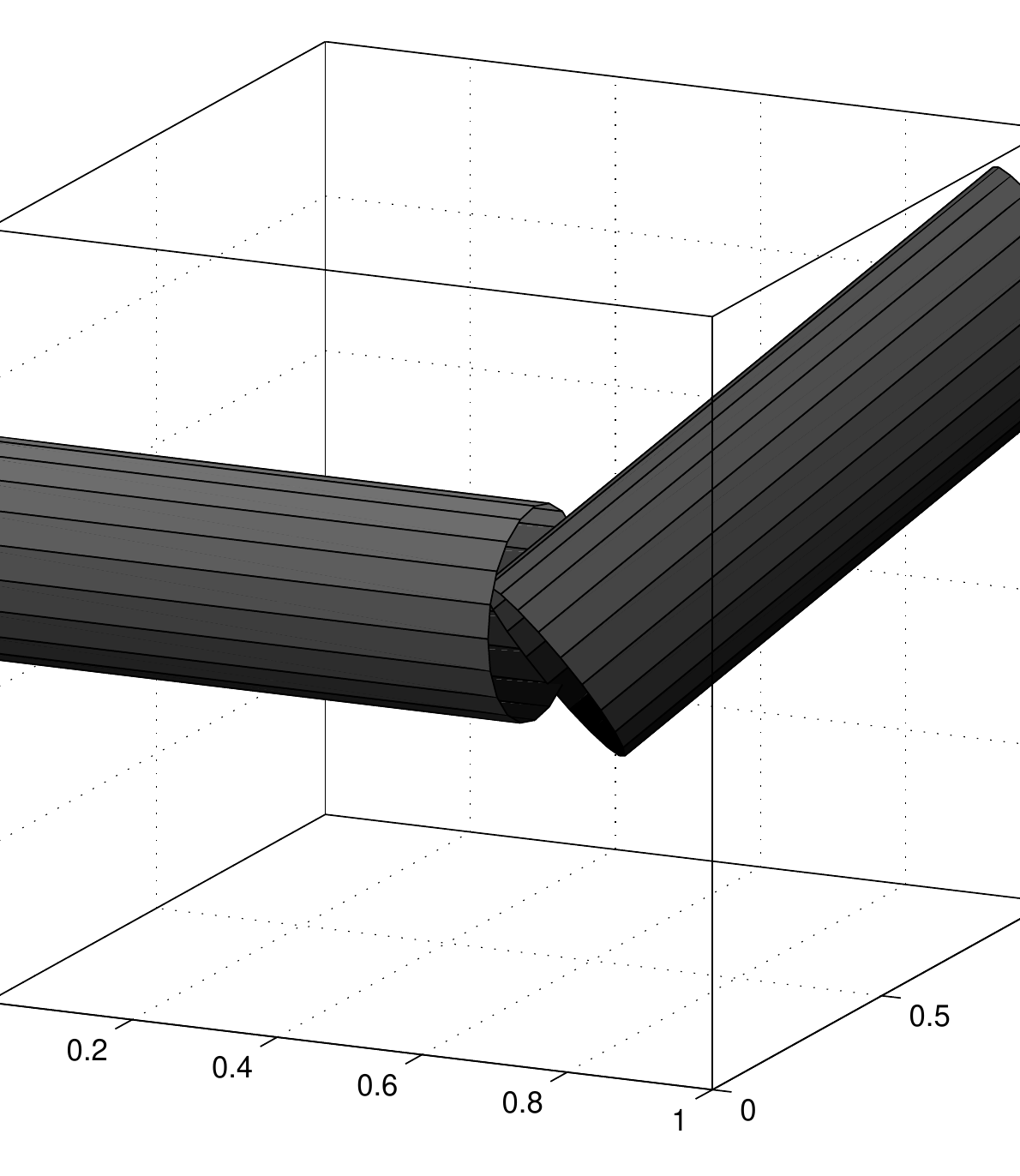}  
}
\subfigure[\, After relaxation]{
    \includegraphics[trim = 5cm 7cm 4.5cm 7cm, clip, width=0.45\linewidth]{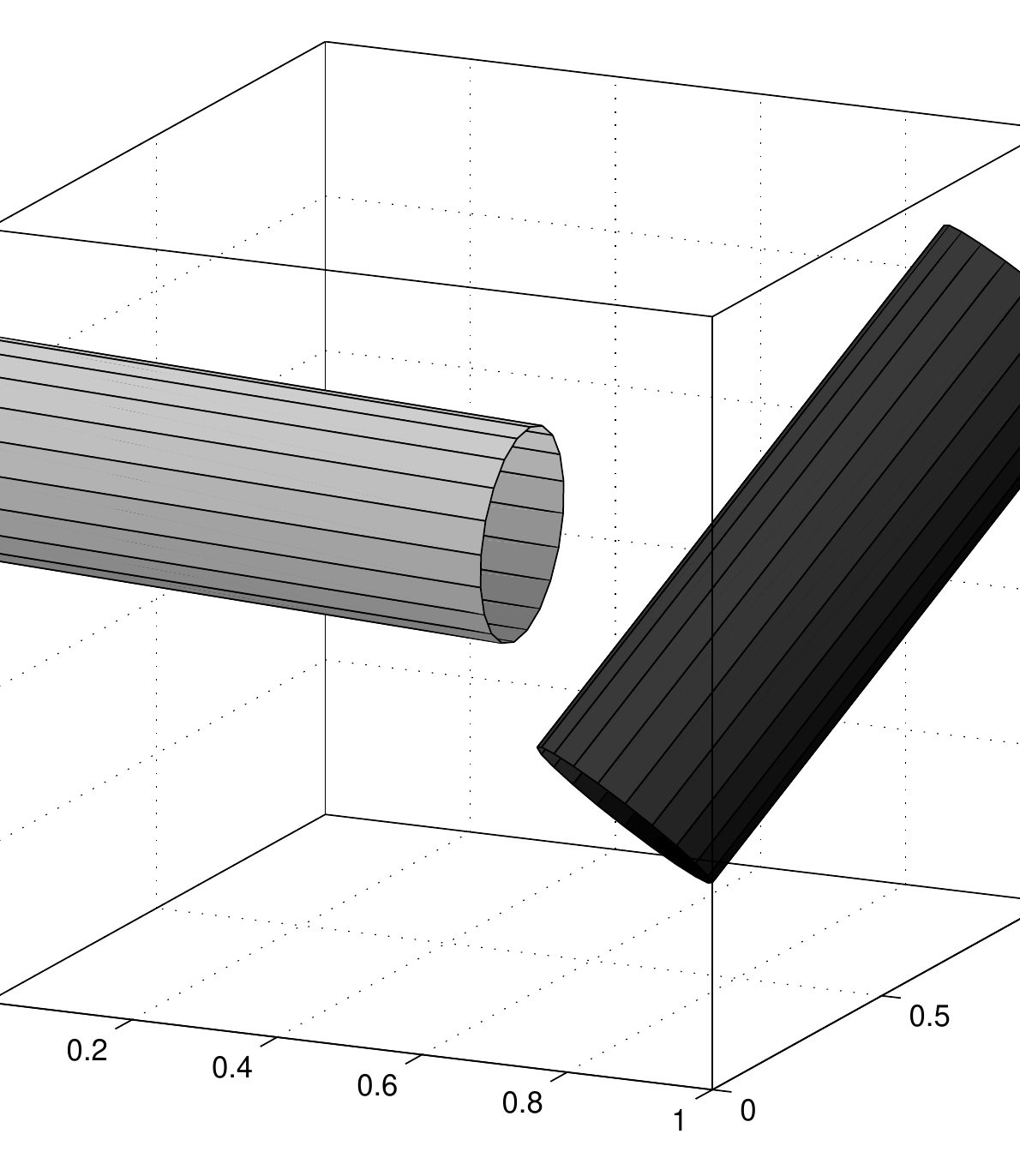}  
 }
\caption{Two cylinders, intersection of bases, type $d1$ or $d2$.}
\label{fig:cyl-cyl8}
\end{figure}

\begin{figure}[H] 
\centering
\subfigure[\, Intersecting]{    
\includegraphics[trim = 5cm 7cm 4.5cm 7cm, clip, width=0.45\linewidth]{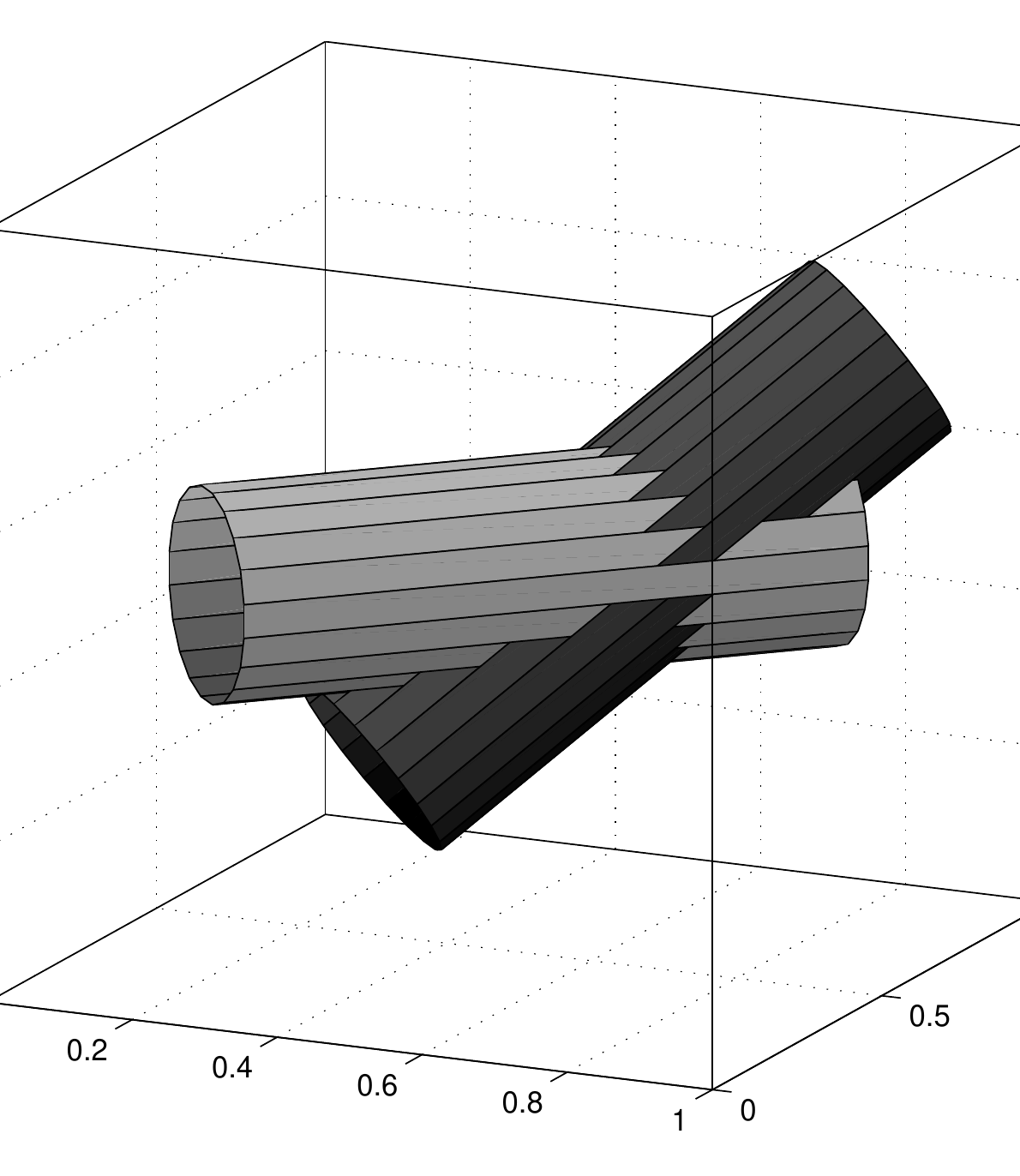}  
}
\subfigure[\, After relaxation]{
    \includegraphics[trim = 5cm 7cm 4.5cm 7cm, clip, width=0.45\linewidth]{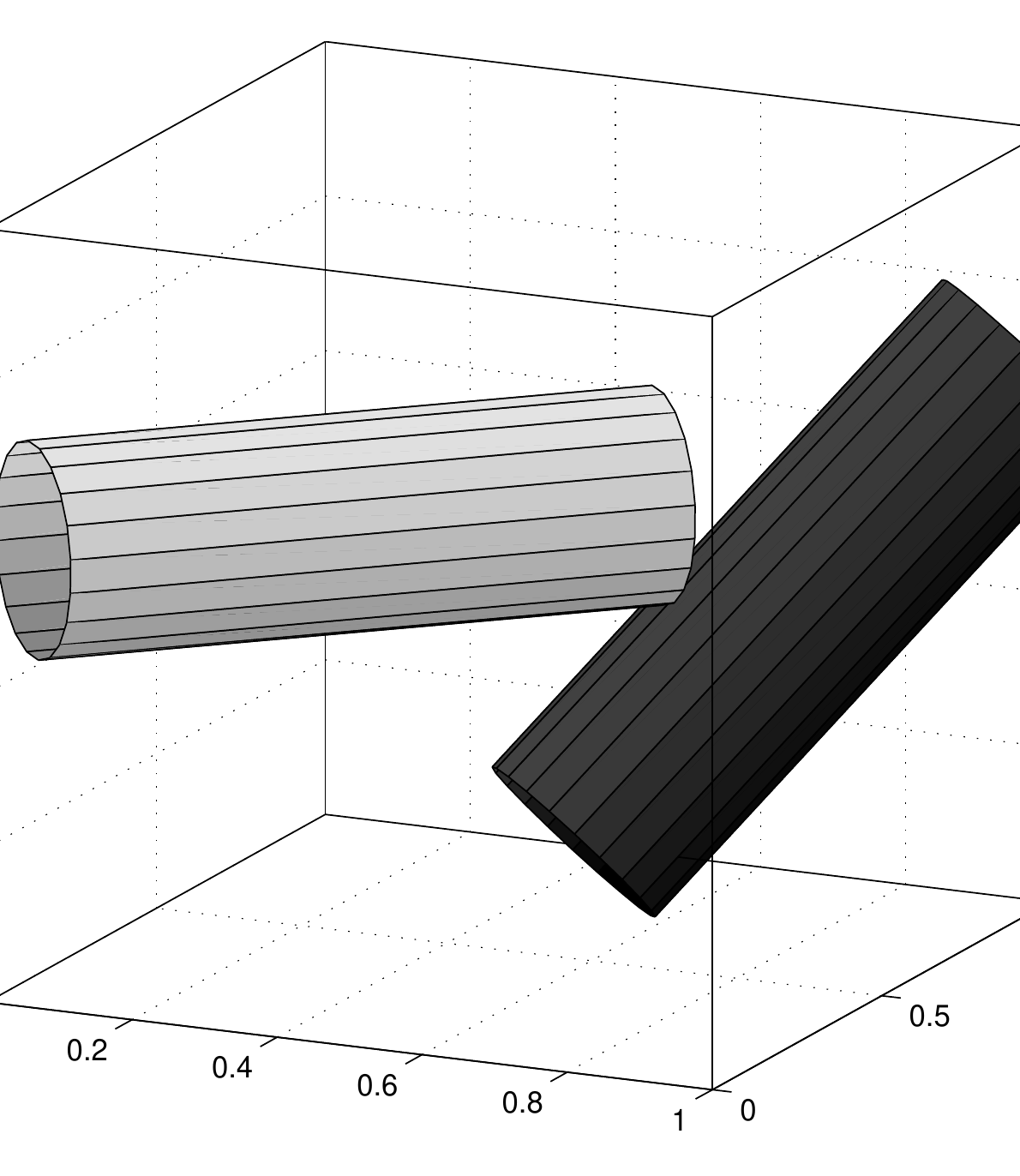}  
 }
\caption{Two cylinders, axes intersect inside the cylinders, type $cc1$ (degenerate case), 
then $cd2$, then $cd1$.}
\label{fig:cyl-cyl9}
\end{figure}

The figures \ref{fig:sp-sp} -- \ref{fig:cyl-cyl9} represent typical intersections of spheres and cylinders 
and the result of application of the relaxation procedure to them.

 \section*{Appendix 2. Time needed for generation using the RSA algorithm}
The tables (\ref{tab:RSA10} -- \ref{tab:RSA50}) show the dependence of time of construction of the RVEs 
following the algorithm \ref{alg:mc_gen}
on various parameters of them: volume fractions $f_s, f_c$, number of inclusions, their geometry.
The couples of numbers in the cells of the tables correspond to two values of the 
aspect ratio $a$ of cylinders (ratio between its length and diameter).
The time estimation (in seconds) is averaged over 20 runs.

\section*{Appendix 3. Time needed for generation using the MD method}

The tables (\ref{tab:MD10} -- \ref{tab:MD50}) show the dependence of time of construction of the RVEs 
using the time-driven MD relaxation method
on various parameters of them: volume fractions $f_s, f_c$, number of inclusions, their geometry.
The couples of numbers in the cells of the tables correspond to two values of the 
aspect ratio $a$ of cylinders (ratio between its length and diameter).
The time estimation (in seconds) is averaged over 20 runs.


\clearpage

\newpage

\begin{table*}[ht] \hspace{-1.5cm} \centering
 \begin{tabular}{|c|c|c|c|c|c|c|c|c|c|} 
 \hline
   $f_s \backslash f_c$  & $0.05$ & $0.1$ & $0.15$ & $0.2$  & $0.25$ & $0.3$   \\
 \hline $0.05 $ & $0.0002|0.0003$  & $0.0003|0.0005$  & $0.0005|0.0012$  & $0.0011|0.0044$ & $0.005|0.021$ & $0.015|1.02$ \\
 \hline $0.1$ &  $0.0002|0.0003$  &  $0.0004|0.0006$  & $0.0007|0.0026$  & $0.0018|0.0059$ & $0.01|0.057$ &  \\
 \hline $0.15$& $0.0003|0.0004$   & $0.0009|0.0027$   & $0.0034|0.023$  & $0.008|0.21$ & &  \\
 \hline $0.2$ & $0.0013|0.0019$   & $0.017|0.088$   & $0.22|$  &  & &   \\
 \hline $0.25$  & $0.022|0.029$  &  $8.49|$  &  &  & & \\
 \hline $0.3$  & $0.015|1.49$   & & & & &  \\
 \hline 
 \end{tabular}
 \caption{\label{tab:RSA10} \modif{RSA: average time of RVE generation (in seconds) for} 10 spheres, 10 cylinders, values of aspect ratio $a = 3|5$ }
 \end{table*}

\begin{table*}[ht]\hspace{-1.5cm} \centering
 \begin{tabular}{|c|c|c|c|c|c|c|c|c|c|} 
 \hline
   $f_s \backslash f_c$  & $0.05$ & $0.1$ & $0.15$ & $0.2$  & $0.25$ & $0.3$   \\
 \hline $0.05 $ & $0.0004|0.0006$  & $0.0007|0.0013$  & $0.0012|0.0031$  & $0.0027|0.012$ & $0.012|0.068$ & $0.062|2.4$ \\
 \hline $0.1$ &  $0.0004|0.0007$  &  $0.001|0.0017$  & $0.0021|0.0042$  & $0.0052|0.017$ & $0.023|0.14$ &  \\
 \hline $0.15$&  $0.0008|0.0011$   & $0.0024|0.005$   & $0.0094|0.077$  & $2.7|7.7$ & &  \\
 \hline $0.2$ & $0.0017|0.003$   & $0.12|0.82$   & $5.8|16.55$  &  & &   \\
 \hline $0.25$  & $0.086|0.17$  &  $2.86|$  &  &  & & \\
 \hline 
 \end{tabular}
 \caption{\label{tab:RSA20} \modif{RSA: average time of RVE generation (in seconds) for} 20 spheres, 20 cylinders, values of aspect ratio $a = 3|5$ }
 \end{table*}

\begin{table*}[ht]\hspace{-1.5cm} \centering
 \begin{tabular}{|c|c|c|c|c|c|c|c|c|c|} 
 \hline
   $f_s \backslash f_c$  & $0.05$ & $0.1$ & $0.15$ & $0.2$  & $0.25$ & $0.3$   \\
 \hline $0.05 $ & $0.0007|0.0014$  & $0.0015|0.0023$  & $0.0026|0.0056$  & $0.0047|0.018$ & $0.017|0.14$ & $0.17|3.36$ \\
 \hline $0.1$ &  $0.001|0.0014$  &  $0.0018|0.0034$  & $0.0041|0.0096$  & $0.01|0.037$ & $0.037|0.68$ &  \\
 \hline $0.15$& $0.0014|0.0021$   & $0.0044|0.0095$   & $0.021|0.099$  & $2.6|13$ & &  \\
 \hline $0.2$ & $0.0027|0.0061$   & $0.21|0.99$   & $6.49|$  &  & &   \\
 \hline $0.25$  & $0.14|1.49$  &    &  &  & & \\
 \hline 
 \end{tabular}
 \caption{\label{tab:RSA30} \modif{RSA: average time of RVE generation (in seconds) for } 30 spheres, 30 cylinders, values of aspect ratio $a = 3|5$ }
 \end{table*}

\begin{table*}[ht]\hspace{-1.5cm} \centering
 \begin{tabular}{|c|c|c|c|c|c|c|c|c|c|} 
 \hline
   $f_s \backslash f_c$  & $0.05$ & $0.1$ & $0.15$ & $0.2$  & $0.25$ & $0.3$   \\
 \hline $0.05 $ & $0.0012|0.0016$  & $0.002|0.0037$  & $0.0034|0.0079$  & $0.0077|0.029$ & $0.034|0.24$ & $0.26|6.74$ \\
 \hline $0.1$ &  $0.0015|0.0021$  &  $0.0027|0.005$  & $0.0058|0.014$  & $0.016|0.057$ & $0.079|0.58$ &  \\
 \hline $0.15$&  $0.0024|0.0035$   & $0.0064|0.013$   & $0.035|0.2$  & $1.78|17.2$ & &  \\
 \hline $0.2$ &  $0.0064|0.011$   & $0.14|1.18$   & $10.6|$  &  & &   \\
 \hline $0.25$  & $0.5|1.6$  &    &  &  & & \\
 \hline 
 \end{tabular}
 \caption{\label{tab:RSA40} \modif{RSA: average time of RVE generation (in seconds) for }40 spheres, 40 cylinders, values of aspect ratio $a = 3|5$ }
 \end{table*}

\begin{table*}[ht]\hspace{-1.5cm} \centering
 \begin{tabular}{|c|c|c|c|c|c|c|c|c|c|} 
 \hline
   $f_s \backslash f_c$  & $0.05$ & $0.1$ & $0.15$ & $0.2$  & $0.25$ & $0.3$   \\
 \hline $0.05 $ & $0.0017|0.0023$  & $0.0027|0.0045$  & $0.0045|0.012$  & $0.011|0.038$ & $0.043|0.38$ & $0.29|8.38$ \\
 \hline $0.1$ &  $0.0021|0.0029$  &  $0.0043|0.0067$  & $0.0082|0.019$  & $0.025|0.08$ & $0.11|0.69$ & $1.44|32$ \\
 \hline $0.15$&  $0.0035|0.0042$   & $0.0086|0.017$   & $0.057|0.29$  & $3.85|$ & &  \\
 \hline $0.2$ &  $0.0092|0.015$   & $0.21|4.94$   & $21.8|$  &  & &   \\
 \hline $0.25$ & $0.21|4.02$  &    &  &  & & \\
 \hline 
 \end{tabular}
 \caption{\label{tab:RSA50} \modif{RSA: average time of RVE generation (in seconds) for }50 spheres, 50 cylinders, values of aspect ratio $a = 3|5$ }
 \end{table*}

 \clearpage

\begin{table*}[ht] \centering 
\hskip-2.5cm
\begin{tabular}{|c|c|c|c|c|c|c|c|c|c|} 
\hline
  $f_s \backslash f_c$  & $0.05$ & $0.1$& $0.15$ & $0.2$  & $0.25$ & $0.3$ & $0.35$ & $0.4$ & $0.45$ \\
\hline $0.05$ &$0.06|0.07$&$0.09|0.12$&  $0.12|0.16$ & $0.14|0.19$ & $0.19|0.37$ & $0.27|0.52$ & $0.38|0.82$ & $0.66|1.62$ &  $0.74|2.17$\\
\hline $0.1$ &$0.07|0.1$&$0.12|0.12$& $0.11|0.16$  & $0.18|0.22$ & $0.21|0.33$ & $0.29|0.7$ & $0.38|1.0$ & $0.63|2.08$ &  \\
\hline $0.15$&$0.08|0.1$&$0.12|0.13$& $0.13|0.22$  & $0.16|0.3$ & $0.27|0.44$ & $0.31|0.72$ & $0.56|2.0$ & & \\
\hline $0.2$ &$0.09|0.1$&$0.12|0.16$&$0.16|0.21$   & $0.23|0.39$ & $0.28|0.54$ & $0.49|0.9$ &  & &\\
\hline $0.25$ &$0.1|0.14$&$0.14|0.18$& $0.19|0.25$  & $0.29|0.43$ &  $0.5|0.93$  & & & &\\
\hline $0.3$  &$0.12|0.14$&$0.15|0.21$& $0.20|0.35$  & $0.35|0.65$ & & & & &\\
\hline $0.35$ &$0.13|0.14$&$0.2|0.25$&  $0.26|0.48$ & & & & & &\\
\hline $0.4$  &$0.15|0.16$&$0.21|0.3$& & & & & & &\\
\hline $0.45$ &$0.17|0.22$& & & & & & & &\\
\hline 
\end{tabular}
\caption{\label{tab:MD10}\modif{MD: average time of RVE generation (in seconds) for } 10 spheres, 10 cylinders,  values of aspect ratio $a = 3|5$ }
\end{table*}

 \begin{table*}[ht] \hskip-2.5cm \centering
 \begin{tabular}{|c|c|c|c|c|c|c|c|c|c|} 
 \hline
   $f_s \backslash f_c$  & $0.05$ & $0.1$ & $0.15$ & $0.2$  & $0.25$ & $0.3$ & $0.35$ & $0.4$ & $0.45$ \\
 \hline $0.05 $ & $0.2|0.25$  & $0.29|0.35$  & $0.38|0.51$  & $0.49|0.69$ & $0.66|1.12$ & $0.84|1.7$ & $0.14|4.1$ & $2.13|8.07$ &  $4.3|7.07$\\
 \hline $0.1$ &  $0.22|0.27$  &  $0.32|0.4$  & $0.41|0.53$  & $0.53|0.84$ & $0.7|1.37$ & $1.07|1.96$ & $1.67|5.6$ & $2.6|10.5$ &  \\
 \hline $0.15$& $0.23|0.29$   & $0.35|0.4$   & $0.45|0.59$  & $0.64|0.87$ & $0.79|1.52$ & $1.21|3$ & $2.08|8.75$ & & \\
 \hline $0.2$ & $0.31|0.33$   & $0.39|0.47$   & $0.47|0.62$  & $0.65|1.11$ & $0.9|1.59$ & $1.86|5.51$ &  & &\\
 \hline $0.25$  & $0.31|0.39$  &  $0.4|0.49$  & $0.53|0.73$  & $0.79|1.22$ &  $1.18|4.29$  & & & &\\
 \hline $0.3$  & $0.39|0.44$   & $0.5|0.54$   & $0.68|0.89$  & $1.07|2.56$ & & & & &\\
 \hline $0.35$ &  $0.42|0.48$ & $0.62|0.66$    &  $0.85|1.39$ & & & & & &\\
 \hline $0.4$  &  $0.45|0.57$ & $0.72|0.81$   & & & & & & &\\
 \hline $0.45$ &  $0.55|0.66$  & & & & & & & &\\
 \hline 
 \end{tabular}
 \caption{\label{tab:MD20}\modif{MD: average time of RVE generation (in seconds) for } 20 spheres, 20 cylinders, values of aspect ratio $a = 3|5$ }
 \end{table*}

  \begin{table*}[ht] \hskip-2.5cm \centering
 \begin{tabular}{|c|c|c|c|c|c|c|c|c|c|} 
 \hline
   $f_s \backslash f_c$  & $0.05$ & $0.1$ & $0.15$ & $0.2$  & $0.25$ & $0.3$ & $0.35$ & $0.4$ & $0.45$ \\
 \hline $0.05 $ & $0.39|0.5$  & $0.59|0.73$  & $0.75|1.06$  & $0.96|1.33$ & $1.15|1.92$ & $1.62|3.57$ & $2.19|8.85$ & $4.01|16.94$ &  $7.36|$\\
 \hline $0.1$ &  $0.47|0.62$  &  $0.59|0.92$  & $0.81|1.1$  & $0.99|1.57$ & $1.22|2.63$ & $1.77|4.4$ & $3.22|9.6$ & $5.77|17.8$ &  \\
 \hline $0.15$& $0.51|0.61$   & $0.68|0.84$   & $0.88|1.17$  & $1.13|1.63$ & $1.53|3.26$ & $2.3|5.62$ & $4.17|15.2$ & & \\
 \hline $0.2$ & $0.59|0.67$   & $0.75|0.96$   & $0.97|1.3$  & $1.29|2.21$ & $1.94|4.7$ & $3.02|11.2$ &  & &\\
 \hline $0.25$  & $0.66|0.74$  &  $0.86|1.06$  & $1.13|1.72$  & $1.51|2.81$ &  $2.55|7.85$  & & & &\\
 \hline $0.3$  & $0.74|0.81$   & $1.01|1.19$   & $1.29|1.88$  & $2.1|4.03$ & & & & &\\
 \hline $0.35$ &  $0.86|0.96$ & $1.13|1.37$    &  $1.7|2.56$ & & & & & &\\
 \hline $0.4$  &  $0.94|1.01$ & $1.38|1.9$   & & & & & & &\\
 \hline $0.45$ &  $1.1|1.32$  & & & & & & & &\\
 \hline 
 \end{tabular}
 \caption{\label{tab:MD30} \modif{MD: average time of RVE generation (in seconds) for }30 spheres, 30 cylinders, values of aspect ratio $a = 3|5$ }
 \end{table*}

  \begin{table*}[ht]\hskip-2.5cm \centering
 \begin{tabular}{|c|c|c|c|c|c|c|c|c|c|} 
 \hline
   $f_s \backslash f_c$  & $0.05$ & $0.1$ & $0.15$ & $0.2$  & $0.25$ & $0.3$ & $0.35$ & $0.4$ & $0.45$ \\
 \hline $0.05 $ & $0.76|0.81$  & $0.94|1.17$  & $1.23|1.69$  & $1.61|2.55$ & $2.1|3.14$ & $2.52|6.23$ & $4.23|14$ & $6.47|25$ &  $12.2|$\\
 \hline $0.1$ &  $0.81|0.82$  &  $1.09|1.15$  & $1.37|1.57$  & $1.83|2.21$ & $2.24|3.48$ & $3.4|8.18$ & $5.08|20.7$ & $8.94|24.9$ &  \\
 \hline $0.15$&  $0.9|0.98$   & $1.16|1.31$   & $1.42|1.72$  & $1.85|2.83$ & $2.53|4.86$ & $3.96|10.3$ & $6.33|26.7$ & & \\
 \hline $0.2$ & $0.97|1.1$   & $1.33|1.5$   & $1.61|2.06$  & $2.17|3.31$ & $3.29|7.21$ & $5.42|16.3$ &  & &\\
 \hline $0.25$  & $1.1|1.14$  &  $1.39|1.62$  & $1.89|2.44$  & $2.54|4.38$ &  $4.27|6.39$  & & & &\\
 \hline $0.3$  & $1.21|1.27$   & $1.62|1.87$   & $2.21|2.96$  & $3.47|6.39$ & & & & &\\
 \hline $0.35$ &  $1.44|1.43$ & $1.94|2.32$    &  $2.87|3.91$ & & & & & &\\
 \hline $0.4$  &  $1.59|1.67$ & $2.25|2.94$   & & & & & & &\\
 \hline $0.45$ &  $2.01|2.15$  & & & & & & & &\\
 \hline 
 \end{tabular}
 \caption{\label{tab:MD40} \modif{MD: average time of RVE generation (in seconds) for }40 spheres, 40 cylinders, values of aspect ratio $a = 3|5$ }
 \end{table*}

  \begin{table*}[ht]\hskip-2.5cm \centering
 \begin{tabular}{|c|c|c|c|c|c|c|c|c|c|} 
 \hline
   $f_s \backslash f_c$  & $0.05$ & $0.1$ & $0.15$ & $0.2$  & $0.25$ & $0.3$ & $0.35$ & $0.4$ & $0.45$ \\
 \hline $0.05 $ & $1.15|1.18$  & $1.38|1.77$  & $1.93|2.37$  & $2.36|3.24$ & $3.04|5.32$ & $3.99|9.01$ & $6.47|20.2$ & $8.8|31.4$ &  $15.6|$\\
 \hline $0.1$ &  $1.2|1.24$  &  $1.56|1.73$  & $1.87|2.42$  & $2.58|3.22$ & $3.41|5.65$ & $4.7|11.5$ & $7.36|22.7$ & $14.3|39.2$ &  \\
 \hline $0.15$& $1.34|1.38$   & $1.67|1.97$   & $2.18|2.5$  & $2.91|4.37$ & $3.81|7.09$ & $5.54|16.2$ & $10.2|44.4$ & & \\
 \hline $0.2$ & $1.5|1.48$   & $1.92|2.11$   & $2.39|3.06$  & $3.31|4.38$ & $4.6|8.27$ & $8.34|25.2$ &  & &\\
 \hline $0.25$  & $1.61|1.7$  &  $2.14|2.36$  & $2.69|3.5$  & $3.84|6.96$ &  $6.39|17.3$  & & & &\\
 \hline $0.3$  & $1.77|1.85$   & $2.3|2.66$   & $3.3|4.37$  & $4.99|9.89$ & & & & &\\
 \hline $0.35$ &  $2.08|2.14$ & $2.8|3.28$    &  $4.28|5.79$ & & & & & &\\
 \hline $0.4$  &  $2.27|2.64$ & $3.47|4.02$   & & & & & & &\\
 \hline $0.45$ &  $2.69|3.33$  & & & & & & & &\\
 \hline 
 \end{tabular}
 \caption{\label{tab:MD50}\modif{MD: average time of RVE generation (in seconds) for } 50 spheres, 50 cylinders,  values of aspect ratio $a = 3|5$ }
 \end{table*}

\end{document}